\documentclass{cmip-l}							%CMI

\usepackage{pstricks}
 \usepackage{graphicx}
\DeclareGraphicsExtensions{.eps} 
\RequirePackage{amsmath,amsthm,amssymb,amsfonts,bbm}

\newcommand{\R}{\mathbb{R}}
\newcommand{\K}{\mathbb{K}}
\newcommand{\Z}{\mathbb{Z}}

\newcommand{\N}{\mathbb{N}}

\newcommand{\T}{\mathbb{T}}
\def\build#1_#2^#3{\mathrel{
\mathop{\kern 0pt#1}\limits_{#2}^{#3}}}
\def\cq{$\hfill \square$}

\def\ve{\varepsilon}

\def\d{\partial}

\def\e{\mathcal{E}}

\def\t{\mathcal{T}}

\def\v{\mathcal{V}}

\def\s{\mathcal{S}}

\def\ov{\overline}

\def\wt{\widetilde}

\def\la{\longrightarrow}

\def\proof{\noindent{\sc Proof.} }

\def\card{{\rm Card\,}}

\def\diam{{\rm diam\,}}
\def\supp{{\rm supp\,}}
\def\sm{\smallskip}
\def\llbracket{[\hspace{-.10em} [ }
\def\rrbracket{ ] \hspace{-.10em}]}
\def\ind{{\bf 1}_}
\def\bn{{\bf n}}
\def\be{{\bf e}}
\def\bT{{\bf T}}
\def\midd{\!\mid\!}

\newcommand{\bp}{{\bf p}}

\newcommand{\bq}{{\bf q}}

\newcommand{\bm}{{\bf m}}

\newcommand{\bt}{{\bf t}}

\newcommand{\bA}{{\bf A}}

\newcommand{\bQ}{{\bf Q}}

\newcommand{\bM}{\mathbf{M}}
\newcommand{\RR}{\mathcal{R}}

\renewcommand{\S}{\mathbb{S}}
\newcommand{\TT}{\mathcal{T}}

\newcommand{\lf}{\lfloor}
\newcommand{\rf}{\rfloor}

\def\cq{$\hfill \square$}

\def\d{{\rm d}}

\def\eps{\varepsilon}

\def\ee{\mathbbm{e}}
\def\ov{\overline}

\def\wt{\widetilde}

\def\la{\longrightarrow}

\def\cat{{\rm Cat}}

\newtheorem{theorem}{Theorem}[section]
\newtheorem{lemma}[theorem]{Lemma}
\newtheorem{proposition}[theorem]{Proposition}
\newtheorem{corollary}[theorem]{Corollary}
\newtheorem{definition}{Definition}[section]
\newtheorem{exercise}[theorem]{Exercise}
\newtheorem{conjecture}{Conjecture}[section]
\newtheorem{remark}[theorem]{Remark}

\begin{document}

\copyrightinfo{2012}{Jean-Fran\c cois Le Gall and Gr\'egory Miermont}

\title{Scaling Limits of Random Trees and Planar Maps} 
\author{Jean-Fran\c cois Le Gall}
\address{Math\'ematiques, bat.425, Universit\'e Paris-Sud, 91405 ORSAY Cedex
FRANCE}
\email{jean-francois.legall@math.u-psud.fr}
\author{Gr\'egory Miermont}
\address{Math\'ematiques, bat.425, Universit\'e Paris-Sud, 91405 ORSAY Cedex
FRANCE}
\email{gregory.miermont@math.u-psud.fr}

%\medskip														%CMI

%\date{\small Universit\'e Paris-Sud}
\maketitle

%\begin{center}{\it\small Lecture notes for the Clay Mathematical				%CMI
 % Institute Summer School in Buzios,\\
%July 11 -- August 7, 2010}
%\end{center}

\tableofcontents

%\bigskip															%CMI

\section{Introduction}

The main goal of these lectures is to present some of the
recent progress 
in the asymptotics for
large random planar maps. Recall that a planar map is simply a 
graph drawn on the two-dimensional sphere and viewed up to 
direct homeomorphisms of the sphere. The faces of the map are
the connected components of the complement of edges, or in
other words the regions of the sphere delimited by the graph.
Special cases of planar maps are triangulations, respectively
quadrangulations, respectively $p$-angulations, where each face is adjacent to exactly $3$,
respectively $4$, respectively $p$, edges (see Section 4 for more precise definitions).

Planar maps play an important
role in several areas of mathematics and physics. They have been
studied extensively in combinatorics since the pioneering work 
of Tutte (see in particular \cite{Tu}), which was motivated by
the famous four-color theorem. Graphs drawn on surfaces also have
important algebraic and geometric applications; see the book \cite{LaZv04}.
In theoretical physics, the enumeration of planar maps (and of maps on surfaces 
of higher genus) has strong connections with matrix models, as shown by
the work of 't Hooft \cite{thooft} and Br\'ezin et al
\cite{BrItPaZu}. More recently, graphs on surfaces have been used in physics 
as discrete models of random geometry in the so-called two-dimensional quantum
gravity; see in particular the book \cite{ADJ} (a different mathematical approach to quantum 
gravity using the Gaussian free field appears in the work
of Duplantier and Sheffield \cite{DS}). A nice account of the connections 
between planar maps and the statistical physics of random surfaces can be
found in Bouttier's thesis \cite{Bo}. From the probabilistic perspective,
a planar map can be viewed as a discretization of a surface, and finding a
continuous limit for large planar maps chosen at random in a suitable class
should lead to an interesting model of a ``Brownian surface''. This is 
of course analogous to the well-known fact that Brownian motion appears
as the scaling limit of long discrete random paths. In a way similar
to the convergence of rescaled random walks to Brownian motion, one
expects that the scaling limit of large random planar maps is
universal in the sense that it should
not depend on the details of the discrete model one is considering. These ideas 
appeared in the pioneering paper of Chassaing and Schaeffer \cite{CSise}
and in the subsequent work of Markert and Mokkadem \cite{MM05} in the case
of quadrangulations, and a little later in Schramm \cite{Sch}, who gave a precise 
form to the question of the existence of a scaling limit for large random 
triangulations of the sphere.

To formulate the latter question, consider a random planar map $M_n$
which is uniformly distributed over a certain class of planar maps (for instance,
triangulations, or quadrangulations) with $n$ faces. Equip the vertex set
$V(M_n)$ with the graph distance $d_{gr}$. It has been known for some
time that the diameter of the resulting metric space is of order $n^{1/4}$ 
when $n$ is large (see \cite{CSise} for the case of quadrangulations).
One then expects that the rescaled random metric spaces $(V(M_n),n^{-1/4}d_{gr})$
will converge in distribution as $n$ tends to infinity
towards a certain random metric space, which should be the same,
up to trivial scaling factors, independently of the class of planar maps we started from.
For the previous convergence to make sense, we need to say what it means for
a sequence of metric spaces to converge. To this end we use the notion of the
Gromov-Hausdorff distance, as it was suggested in \cite{Sch}. Roughly speaking (see
Section 2 for a more precise definition) a sequence $(E_n)$ of compact
metric spaces converges to a limiting space $E_\infty$ if it is possible to embed
isometrically all spaces $E_n$ and $E_\infty$ in the same ``big'' metric space $E$,
in such a way that the Hausdorff distance between $E_n$ and $E_\infty$
tends to $0$ as $n\to\infty$. 

The preceding question of the existence of the scaling limit of large random
planar maps is still open, but there has been significant progress in this
direction, and our aim is to present some of the results that have
been obtained in recent years.

\begin{figure}
\begin{center}
\includegraphics[scale=.25]{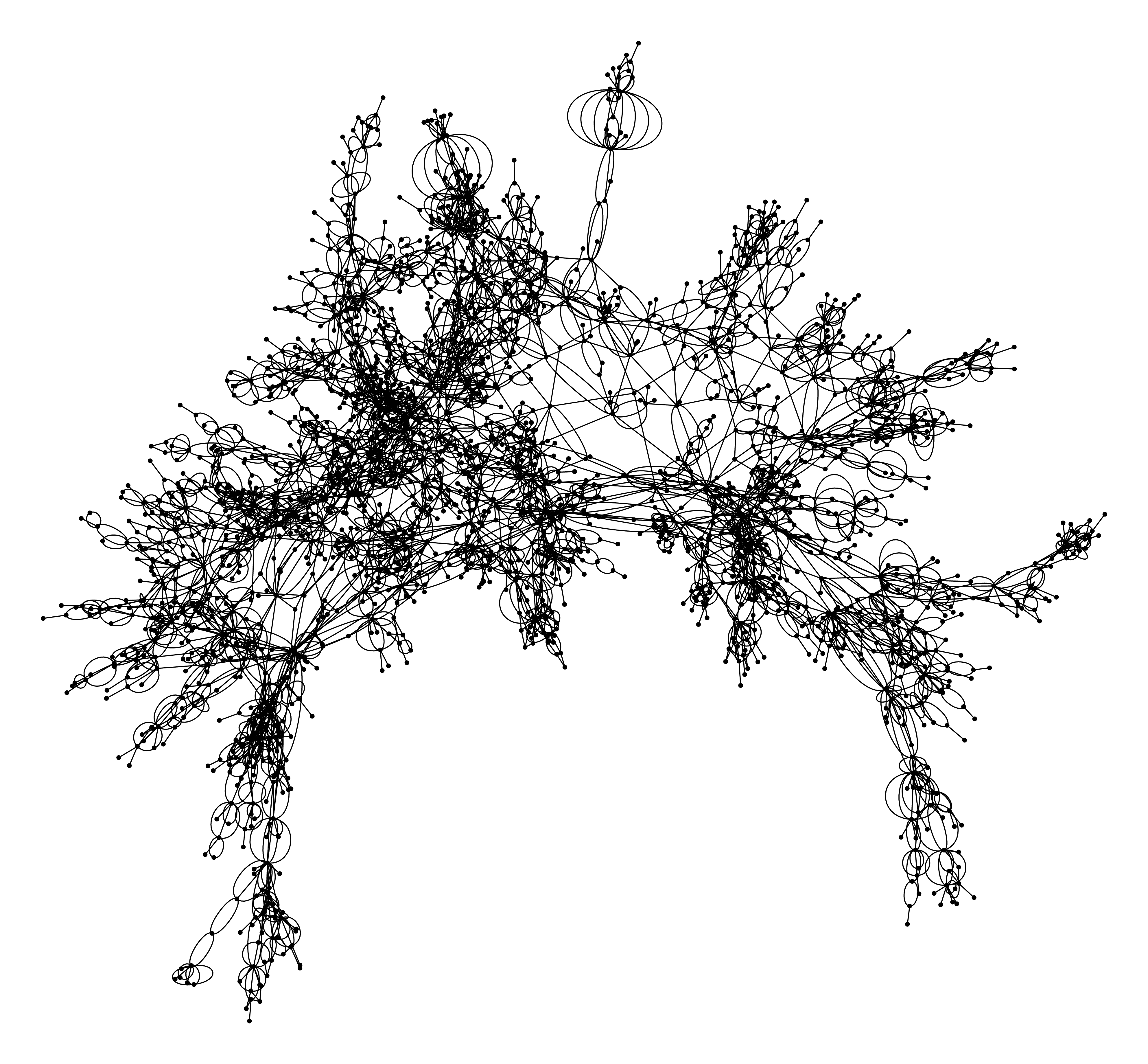}
\includegraphics[scale=.15]{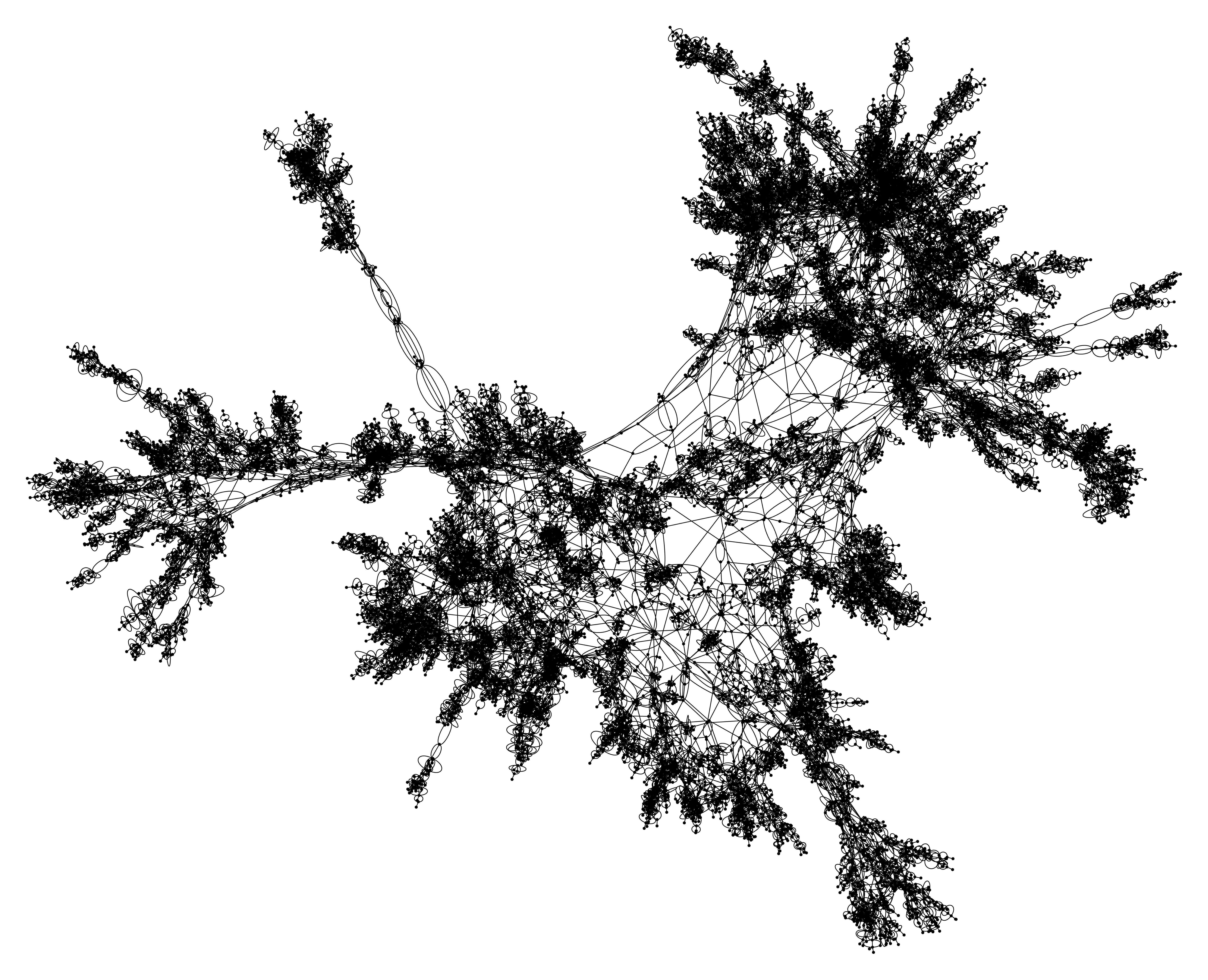}
\caption{Two planar quadrangulations, with respectively 2500 and 20000 vertices. These pictures represent the
quadrangulations as graphs, and do not take account of the embedding in the sphere. Simulations by
J.-F. Marckert.}
\label{quad-simulations}
\end{center}
\end{figure}

Much of the recent progress in the understanding of asymptotic properties of
large random planar maps was made possible by the use of bijections between
different classes of planar maps and certain labeled trees. In the particular case
of quadrangulations, such bijections were discovered by Cori and Vauquelin
\cite{CoVa} and later popularized by Schaeffer \cite{schaeffer98} (see also 
Chassaing and Schaeffer \cite{CSise}). The Cori-Vauquelin-Schaeffer bijection
was extended to much more general planar maps by Bouttier, Di Francesco and 
Guitter \cite{BdFGmobiles}. In the case of bipartite planar maps, this extension
takes a particularly simple form, which explains why some of the recent work 
\cite{jfmgm05,legall06,legall08} concentrates on the bipartite case. The reason why the bijections
between maps and trees are interesting
is the fact that properties of large (labeled) trees are often much easier to
understand than those of large graphs.
Indeed, it has been
known for a long time and in particular since the work of Aldous \cite{Al1,Al3} that one can
often describe the asymptotic properties of large random trees in terms of
``continuous trees'' whose prototype is the so-called CRT or Brownian continuum 
random tree. In the case of trees with labels, the relevant scaling limit
for most of the discrete models of interest is
the CRT equipped with Brownian labels, which can conveniently be constructed 
and studied via 
the path-valued process called the Brownian snake (see e.g. \cite{Zurich}).

A key feature of the bijections between planar maps and labeled trees is the fact that, up
to an appropriate translation, labels on the tree correspond to distances in the map
from a distinguished vertex that plays a special role. Therefore, the known results about scaling
limits of labeled trees immediately give much information about
asymptotics of distances from this distinguished vertex. This idea was
exploited by Chassaing and Schaeffer \cite{CSise} in the case
of quadrangulations and then by
Marckert and Miermont  \cite{jfmgm05} (for bipartite planar maps)
and Miermont \cite{mierinv} (for general planar maps). In view of deriving 
the Gromov-Hausdorff convergence of rescaled planar maps, it is however
not sufficient to control distances from a distinguished vertex. Still, a
simple argument gives an effective bound on the distance between
two arbitrary vertices in terms of quantities depending only on the labels 
on the tree, or equivalently on the distances
from the distinguished vertex (see Proposition \ref{sec:two-useful-bounds-1}(i) below).
This bound was used in \cite{legall06} to show via a compactness argument that the scaling
limit of rescaled uniformly distributed $2p$-angulations with $n$ faces
exists along suitable subsequences. Furthermore,
this scaling limit is a quotient space of the CRT for an equivalence relation defined in
terms of Brownian labels on the CRT: Roughly speaking, two vertices of the CRT
need to be identified if they have the same label and if, when travelling from one 
vertex to the other one along the contour of the CRT, one only encounters
vertices with larger label. The results of \cite{legall06} are not completely
satisfactory, because they require the extraction of suitable subsequences. The reason why this is necessary 
is the fact that the distance on the limiting space (that is, on the quotient of the CRT we have
 just described) has not been fully identified, even though lower and upper bounds are
 available. Still we call Brownian map any random metric space that arises as
 the scaling limit of uniformly distributed $2p$-angulations with $n$ faces. This terminology
 is borrowed from Marckert and Mokkadem \cite{MM05}, who studied a weaker form
 of the convergence of rescaled random quadrangulations. Although the distribution
 of the Brownian map has not been fully characterized, it is possible to derive many 
 properties of this random object (these properties will be common to any of the
 limiting random metric spaces that can arise in the scaling limit). In particular,
 it has been shown that the Brownian map has dimension $4$  \cite{legall06}
 and that it is homeomorphic to the $2$-sphere \cite{lgp,miermont08}. The latter fact is maybe not surprising
 since we started from larger and larger graphs drawn on the sphere: Still it implies that
 large random planar maps will have no ``bottlenecks'', meaning cycles whose
 length is small in comparison with the diameter of the graph but such that both 
 connected components of the complement of the cycle have a macroscopic size.
 
 In the subsequent sections, we discuss most of the preceding results in detail. We
 restrict our attention to the case of quadrangulations, because the bijections
 with trees are simpler in that case: The labeled trees corresponding to quadrangulations
 are just plane trees (rooted ordered trees) equipped with integer labels, such that the
 label of the root is $0$ and the label can change by at most $1$ in absolute value
 along each edge of the tree. 
 
 The first three sections below are devoted to asymptotics for random (labeled) trees, in
 view of our applications to random planar maps.
 In Section 1, we discuss asymptotics for uniformly distributed plane trees with $n$ edges.
 We give a detailed proof of the fact that the suitably rescaled contour function
 of these discrete trees converges in distribution to the normalized Brownian 
 excursion (this is a special case of the results of \cite{Al3}). To this end, we first recall
 the basic facts of excursion theory that we need. In Section 2, we show that
 the convergence of rescaled contour functions can be restated as a convergence
 in the Gromov-Hausdorff sense
 of the trees viewed as random metric spaces for the graph distance. The
 limiting space is then the CRT, which we define precisely as the random 
 real tree coded by a normalized Brownian excursion. Section 2 also contains
 basic facts about the Gromov-Hausdorff distance, and in
 particular its definition in terms of correspondences. In Section 3, we consider
 labeled trees and we give a detailed proof of the fact that rescaled labeled trees converge 
(in a suitable sense) towards the CRT equipped with Brownian labels. 

The last four sections are devoted to planar maps and their scaling limits.
Section 4 presents the combinatorial facts about planar maps that we need.
In particular, we describe the Cori-Vauquelin-Schaeffer bijection between 
(rooted and pointed) quadrangulations and labeled trees.  We also
explain how labels on the tree give access to distances from the distinguished 
vertex in the map, and provide useful upper and lower bounds for
other distances. In Section 5, we give the compactness argument that 
makes it possible
to get sequential limits for rescaled uniformly distributed quadrangulations with
$n$ faces, in the Gromov-Hausdorff sense. The identification of the limit 
(or Brownian map) as 
a quotient space of the CRT for the equivalence relation described above is
explained in Section 6. In that section, we are not able to give the full details of
the proofs, but we try to present the main ideas. As a simple consequence of some
of the estimates needed in the identification of the Brownian map, we also 
compute its Hausdorff dimension. Finally, Section 7 is devoted to the
homeomorphism theorem.  We follow the approach of \cite{miermont08}, which consists in
establishing the absence of  ``bottlenecks'' in the Brownian map
before proving via a theorem of Whyburn that this space is homeomorphic to the sphere.

To conclude this introduction, let us mention that, even though the key problem 
of the uniqueness of the Brownian map remains unsolved, many properties
of this space have been investigated successfully. Often these results 
give insight into the properties of large planar maps. This is in particular the case
for the results of \cite{legall08}, which give a complete description of all geodesics
connecting an arbitrary point of the Brownian map to the distinguished point. Related results
have been obtained in the paper \cite{Mi2}, which deals with maps on surfaces 
of arbitrary genus. Very recently, the homeomorphism theorem of \cite{lgp}
has been extended by Bettinelli \cite{Bet} to higher genus. As a final remark, one expects 
that the Brownian map should be the scaling limit for all random planar maps
subject to some bound on the maximal degree of faces. One may ask what
happens for random planar maps such that the distribution of the degree 
of a typical face has a heavy tail: This problem is discussed in \cite{LGM}, where it
is shown that this case leads to different scaling limits.

\section{Discrete trees and convergence towards the Brownian 
excursion}

\subsection{Plane trees}
\label{SSdiscretetree}

We will be interested in (finite) rooted ordered trees, which
are called plane trees in combinatorics (see e.g. \cite{Sta}). 
We set $\N=\{1,2,\ldots\}$ and by convention $\N^0=\{\varnothing\}$. 
We introduce the set 
$$\mathcal{U}=\bigcup_{n=0}^\infty \N^n.  $$
An
element
of $\mathcal{U}$ is thus a sequence 
$u=(u^1,\ldots,u^n)$ of elements of $\N$, and we set $|u|=n$, so that 
$|u|$ represents the ``generation'' of $u$. If
$u=(u^1,\ldots, u^k)$ and 
$v=(v^1,\ldots, v^\ell)$ belong to $\mathcal{U}$, we write 
$uv=(u^1,\ldots ,u^k,v^1,\ldots ,v^\ell)$
for the concatenation of $u$ and $v$. In particular 
$u\varnothing=\varnothing u=u$.

The mapping $\pi:\mathcal{U}\backslash\{\varnothing\}\la \mathcal{U}$
is defined by $\pi((u^1,\ldots ,u^n))=(u^1,\ldots ,$ \hfill $u^{n-1})$		%CMI
($\pi(u)$ is the ``parent'' of $u$).

\smallskip
A plane tree $\tau$ is a finite subset of
$\mathcal{U}$ such that:
\begin{enumerate}
\item[(i)] $\varnothing\in \tau$.

\item[(ii)] $u\in \tau\backslash\{\varnothing\}\Rightarrow
\pi(u)\in \tau$.

\item[(iii)] For every $u\in\tau$, there exists an integer $k_u(\tau)\geq 0$
such that, for every $j\in\N$,  $uj\in\tau$ if and only if $1\leq j\leq
k_u(\tau)$
\end{enumerate}

The number $k_u(\tau)$ is interpreted as the ``number of children'' of $u$
in
$\tau$.

We denote by ${\bf A}$ the set of all plane trees. In what
follows, we see each vertex of the tree $\tau$ as an individual of a
population  whose $\tau$ is the family tree. By definition, the size $|\tau|$ of $\tau$
is the number of edges of $\tau$, $|\tau|=\#\tau -1$. For every
integer $k\geq 0$,
we put
$${\bf A}_k=\{\tau\in{\bf A}:|\tau|=k\}.$$

\begin{exercise}
Verify that the cardinality of ${\bf A}_k$
is the $k$-th Catalan number
$$\#{\bf A}_k={\rm Cat}_k:= \frac{1}{k+1}{2k\choose k}.$$
\end{exercise}

\smallskip
A plane tree can be coded by its Dyck path or {\bf contour function}. Suppose that the tree is embedded in the half-plane
in  such a way that edges have length one. Informally, we imagine the
motion of a particle that starts at time 
$t=0$ from the root of the tree and then
explores the tree from the left to the right, moving continuously along the edges 
at unit speed (in the way explained by the arrows of Fig.2), until all edges have been
explored and the particle has come back to the root. Since it is
clear that each edge will be crossed twice in this evolution, the total time
needed to explore the tree is $2|\tau|$.
The value $C(s)$ 
of the contour function at time $s\in[0,2|\tau|]$ is the distance (on the tree)
between the position of the particle at time $s$ and the
root. By convention $C(s)=0$ if $s\geq 2|\tau|$.
Fig.2 explains the construction of the contour function better 
than a formal definition.

\begin{figure}

\begin{center}
\unitlength=1pt
\begin{picture}(270,140)

\thicklines \put(50,0){\line(-1,2){20}}
\thicklines \put(50,0){\line(1,2){20}}
\thicklines \put(30,40){\line(0,1){40}}
\thicklines \put(30,40){\line(-1,2){20}}
\thicklines \put(30,40){\line(1,2){20}}
\thicklines \put(30,80){\line(-1,2){20}}
\thicklines \put(30,80){\line(1,2){20}}

\thinlines \put(37,15){\vector(-1,2){6}}
\thinlines \put(17,55){\vector(-1,2){6}}
\thinlines \put(20,67){\vector(1,-2){6}}
\thinlines \put(27,57){\vector(0,1){12}}
\thinlines \put(17,95){\vector(-1,2){6}}
\thinlines \put(22,107){\vector(1,-2){6}}
\thinlines \put(33,95){\vector(1,2){6}}
\thinlines \put(48,107){\vector(-1,-2){6}}
\thinlines \put(33,69){\vector(0,-1){12}}
\thinlines \put(34,55){\vector(1,2){6}}
\thinlines \put(48,67){\vector(-1,-2){6}}
\thinlines \put(40,27){\vector(1,-2){6}}
\thinlines \put(54,15){\vector(1,2){6}}
\thinlines \put(68,27){\vector(-1,-2){6}}

\put(40,-5){$\varnothing$}
\put(20,35){$1$}
\put(60,35){$2$}
\put(12,75){\scriptsize{(1,1)}}
\put(30,75){\scriptsize{(1,2)}}
\put(50,75){\scriptsize{(1,3)}}
\put(12,115){\scriptsize{(1,2,1)}}
\put(50,115){\scriptsize{(1,2,2)}}

\thinlines \put(135,0){\vector(1,0){125}}
\thinlines \put(135,0){\vector(0,1){130}}
\thicklines \put(135,0){\line(1,5){8}}
\thicklines \put(143,40){\line(1,5){8}}
\thicklines \put(151,80){\line(1,-5){8}}
\thicklines \put(159,40){\line(1,5){8}}
\thicklines \put(167,80){\line(1,5){8}}
\thicklines \put(175,120){\line(1,-5){8}}
\thicklines \put(183,80){\line(1,5){8}}
\thicklines \put(191,120){\line(1,-5){8}}
\thicklines \put(199,80){\line(1,-5){8}}
\thicklines \put(207,40){\line(1,5){8}}
\thicklines \put(215,80){\line(1,-5){8}}
\thicklines \put(223,40){\line(1,-5){8}}
\thicklines \put(231,0){\line(1,5){8}}
\thicklines \put(239,40){\line(1,-5){8}}

\thinlines \put(143,0){\line(0,1){2}}
\thinlines \put(151,0){\line(0,1){2}}
\thinlines \put(159,0){\line(0,1){2}}
\thinlines \put(135,40){\line(1,0){2}}
\thinlines \put(135,80){\line(1,0){2}}
\thinlines \put(247,0){\line(0,1){2}}

\put(130,39){\scriptsize 1}
\put(130,79){\scriptsize 2}
\put(142,-6){\scriptsize 1}
\put(150,-6){\scriptsize 2}
\put(158,-6){\scriptsize 3}
\put(240,-8){$\scriptstyle 2|\tau|$}
\put (138,120){$C(s)$}
\put(255,4){$s$}

%\put(30,-20){Tree $\tau$}
%\put(150,-20){Contour function $C(s)$}

\end{picture}

\end{center}
\caption{A tree and its contour function}
\end{figure}

Let $k\geq 0$ be an integer. A Dyck path of length $2k$ is a sequence $(x_0,x_1,x_2,\ldots,$ \hfill $x_{2k})$ %CMI
of nonnegative integers such that $x_0=x_{2k}=0$, and $|x_i-x_{i-1}|=1$ for every 
$i=1,\ldots,2k$. Clearly, if $\tau$ is a plane tree of size $k$, and $(C(s))_{s\geq 0}$
is its contour function, the sequence $(C(0),C(1),\ldots,C({2k}))$ is a Dyck path of
length $2k$. More precisely, we have the following easy result.

\begin{proposition}
\label{Dyck}
The mapping $\tau \mapsto (C(0),C(1),\ldots,C({2k}))$ is a bijection from 
${\bf A}_k$ onto the set of all Dyck paths of length $2k$.
\end{proposition}

\subsection{Galton-Watson trees}

Let $\mu$ be a critical or subcritical offspring distribution. 
This means that
$\mu$ is a probability measure on $\Z_+$ such that
$$\sum_{k=0}^\infty k\mu(k)\leq 1.$$
We exclude the trivial case where $\mu(1)=1$.

To define Galton-Watson trees, we let
$(K_u,u\in\mathcal{U})$ be a collection of independent random variables
with law $\mu$, indexed by the set $\mathcal{U}$. Denote by
$\theta$ the random subset of $\mathcal{U}$
defined by
$$\theta=\{u=(u^1,\ldots ,u^n)\in \mathcal{U}:u^j\leq K_{(u^1,\ldots ,u^{j-1})}
\hbox{ for every }1\leq j\leq n\}.$$

\begin{proposition}
\label{GW}
$\theta$ is a.s. a tree. Moreover, if
$$Z_n=\#\{u\in\theta:|u|=n\},$$
$(Z_n,n\geq 0)$ is a Galton-Watson process with offspring 
distribution $\mu$ and initial value $Z_0=1$.
\end{proposition}

\begin{remark} Clearly $k_u(\theta)=K_u$ for every $u\in \theta$.
\end{remark}

The tree $\theta$, or any random tree
with the same distribution, will be 
called a Galton-Watson tree
with offspring distribution $\mu$, or in short a 
$\mu$-Galton-Watson tree. We also write $\Pi_\mu$
for the distribution of $\theta$ on the space $\bf A$.

\smallskip
We leave the easy proof of the proposition to the reader. 
The finiteness of the tree $\theta$ comes from
the fact that the Galton-Watson process with offspring distribution
$\mu$
becomes extinct a.s., so that $Z_n=0$ for $n$ large.

\smallskip
If $\tau$ is a tree and $1\leq j\leq k_\varnothing(\tau)$, we write
$T_j\tau$ for the tree $\tau$ shifted at $j$:
$$T_j\tau=\{u\in\mathcal{U}:ju\in \tau\}.$$
Note that $T_j\tau$ is a tree.

Then $\Pi _\mu$ may be characterized by the following
two properties (see e.g. \cite{Nev} for more general statements):

\begin{enumerate}
\item[(i)] $\Pi _\mu(k_\varnothing=j)=\mu(j)$,\quad
$j\in \Z_+$.

\item[(ii)] For every $j\geq 1$ with $\mu(j)>0$, the shifted trees
$T_1\tau,\ldots,T_j\tau$
are independent under the conditional probability
$\Pi _\mu(d\tau\mid k_\varnothing =j)$
and their conditional distribution is $\Pi _\mu$.
\end{enumerate}

Property (ii) is often called the branching property of the 
Galton-Watson tree.

We now give an explicit formula for $\Pi_\mu$.

\begin{proposition}
\label{lawGW}
For every $\tau\in{\bf A}$,
$$\Pi_\mu(\tau)=\prod_{u\in\tau} \mu(k_u(\tau)).$$
\end{proposition}

\proof We can easily check that
$$\{\theta=\tau\}=\bigcap_{u\in\tau}\{K_u=k_u(\tau)\},$$
so that
$$\Pi_\mu(\tau)=P(\theta=\tau)=
\prod_{u\in\tau} P(K_u=k_u(\tau))
=\prod_{u\in\tau} \mu(k_u(\tau)).$$
\par\cq

\smallskip
We will be interested in the particular case when $\mu=\mu_0$
is the (critical) geometric offspring distribution, $\mu_0(k)=2^{-k-1}$
for every $k\in \Z_+$. In that case, the proposition gives
$$\Pi_{\mu_0}(\tau)=2^{-2|\tau|-1}$$
(note that $\sum_{u\in\tau} k_u(\tau)= |\tau|$ for every $\tau\in{\bf A}$).

In particular $\Pi_{\mu_0}(\tau)$ only depends on $|\tau|$. 
As a consequence, for every integer $k\geq 0$, the conditional probability distribution 
$\Pi_{\mu_0}(\cdot\mid |\tau|=k)$ is just the uniform
probability measure on ${\bf A}_k$. This fact will be important later.

\subsection{The contour function in the geometric case}

In general, the Dyck path of a Galton-Watson tree does not have a ``nice''
probabilistic structure (see however Section 1 of \cite{probasur}).
In this section we restrict our attention to the case when 
$\mu=\mu_0$ is the critical geometric offspring distribution.

First recall that $(S_n)_{n\geq 0}$ is a simple random walk on $\Z$
(started from $0$) if it can be written as
$$S_n=X_1+X_2+\cdots+X_n$$
where $X_1,X_2,\ldots$ are i.i.d. random variables with distribution
$P(X_n=1)=P(X_n=-1)=\frac{1}{2}$. 

Set $T=\inf\{n\geq 0: S_n=-1\} <\infty$ a.s. The random finite path
$$(S_0,S_1,\ldots,S_{T-1})$$
(or any random path with the same distribution)
is called an excursion of simple random walk. Obviously this
random path is a random Dyck path of length $T-1$. 

\begin{proposition}
\label{Dyck-excursion}
Let $\theta$ be a $\mu_0$-Galton-Watson tree. Then the Dyck path
of $\theta$ is an excursion of simple random walk.
\end{proposition}

\proof Since plane trees are in one-to-one correspondence with 
Dyck paths (Proposition \ref{Dyck}), the statement of the proposition
is equivalent to saying that the random plane tree $\theta$ coded by an 
excursion of simple random walk is a $\mu_0$-Galton-Watson tree. 
To see this, introduce the upcrossing times of the random walk $S$
from $0$ to $1$: 
$$U_1=\inf\{n\geq 0: S_n=1\}\ ,\ V_1=\inf\{n\geq U_1: S_n=0\}$$
and by induction, for every $j\geq 1$,
$$U_{j+1}=\inf\{n\geq V_j: S_n=1\}\ ,\ V_{j+1}=\inf\{n\geq U_{j+1}: S_n=0\}.$$
Let $K=\sup\{j: U_j \leq T\}$ ($\sup\varnothing =0$). From the 
relation between a plane tree and its associated Dyck path, one easily sees
that $k_\varnothing(\theta)=K$, and that for every $i=1,\ldots, K$,
the Dyck path associated with the subtree $T_i\theta$ is the path $\omega_i$, with
$$\omega_i(n):=S_{(U_i+n)\wedge (V_i-1)}-1\quad,\ 0\leq n\leq V_i-U_i-1.$$

A simple application of the Markov property now shows that $K$ is 
distributed according to $\mu_0$ and that conditionally on $K=k$,
the paths $\omega_1,\ldots,\omega_k$ are $k$ independent 
excursions of simple random walk. The characterization of $\Pi_{\mu_0}$ by
properties (i) and (ii) listed before Proposition \ref{lawGW} now shows that
$\theta$ is a $\mu_0$-Galton-Watson-tree. \cq

\subsection{Brownian excursions}

Our goal is to prove that the (suitably rescaled) contour function of a 
tree uniformly distributed over ${\bf A}_k$ converges in distribution
as $k\to\infty$ towards a normalized Brownian excursion. We first need
to recall some basic facts about Brownian excursions.

We consider a standard linear Brownian motion $B=(B_t)_{t\geq 0}$ starting
from the origin. The process $\beta_t=|B_t|$ is called reflected Brownian
motion. We denote by $(L^0_t)_{t\geq 0}$ the local time process of $B$
(or of $\beta$)
at level $0$, which can be defined by the approximation
$$L^0_t = \lim_{\ve \to 0}\frac{1}{2\ve} \int_0^t ds\,\ind{[-\ve,\ve]}(B_s)
= \lim_{\ve \to 0}\frac{1}{2\ve} \int_0^t ds\,\ind{[0,\ve]}(\beta_s),$$
for every $t\geq 0$, a.s. 

Then $(L^0_t)_{t\geq 0}$ is a continuous increasing process, and the set 
of increase points of the function $t\to L^0_t$ coincides with the set 
$${\mathcal Z}=\{t\geq 0: \beta_t=0\}$$
of all zeros of $\beta$. Consequently, if we introduce the right-continuous inverse
of  the local time process,
$$\sigma_\ell:=\inf\{t\geq 0: L^0_t>\ell\}\ , \quad\hbox{for every }\ell\geq 0,$$
we have
$${\mathcal Z}=\{\sigma_\ell:\ell\geq 0\} \cup \{\sigma_{\ell-}:\ell\in D\}$$
where $D$ denotes the countable set of all discontinuity times of the mapping
$\ell\to \sigma_\ell$. 

The connected components  of the open set $\R_+\backslash {\mathcal Z}$ 
are called the {\it excursion intervals} of $\beta$ away from $0$. The preceding discussion shows that,
with probability one,  the
excursion intervals of $\beta$ away from $0$ are exactly the intervals 
$(\sigma_{\ell-},\sigma_\ell)$ for $\ell\in D$. Then, for every $\ell\in D$, we define the 
excursion $e_\ell=(e_\ell(t))_{t\geq 0}$ associated with the interval $(\sigma_{\ell-},\sigma_\ell)$ by setting
$$e_\ell(t)=\left\{\begin{array}{ll}
\beta_{\sigma_{\ell-}+t}\quad&{\rm if}\ 0\leq t\leq \sigma_\ell-\sigma_{\ell-}\,,\\
0&{\rm if}\ t >\sigma_\ell-\sigma_{\ell-}\,.
\end{array}
\right.$$
We view $e_\ell$ as an element of the excursion space $E$, which is defined
by
$$E=\{e\in C(\R_+,\R_+): e(0)=0 \hbox{ and }\zeta(e):=\sup\{s>0:e(s)>0\}\!\in\!(0,\infty)\},$$
where $\sup\varnothing = 0$ by convention. Note that we require $\zeta(e)>0$, so that the zero
function does not belong to $E$.
The space $E$ is equipped with the metric $d$ defined by
$$d(e,e')=\sup_{t\geq 0} |e(t)-e'(t)| + |\zeta(e)-\zeta(e')|$$
and with the associated Borel $\sigma$-field. Notice that $\zeta(e_\ell)=\sigma_\ell-\sigma_{\ell-}$
for every $\ell\in D$.
The following theorem is the basic result of excursion theory in our particular setting. 

\begin{theorem}
\label{Itoth}
The point measure 
$$\sum_{\ell\in D} \delta_{(\ell,e_\ell)}(ds\,de)$$
is a Poisson measure on $\R_+\times E$, with intensity
$$2 ds\otimes \bn(de)$$
where $\bn(de)$ is a $\sigma$-finite measure on $E$.
\end{theorem}

The measure $\bn(de)$ is called the It\^o measure of positive excursions
of linear Brownian motion, or simply the It\^o excursion measure (our measure 
$\bn$ corresponds to the measure $n_+$ in Chapter XII of \cite{RY}). The next corollary follows from standard properties of Poisson measures.

\begin{corollary}
\label{hittingPoisson}
Let $A$ be a measurable subset of $E$ such that $0<\bn(A)<\infty$, and let
$T_A=\inf\{\ell\in D : e_\ell\in A\}$. Then, $T_A$ is exponentially distributed with
parameter $\bn(A)$, and the distribution of $e_{T_A}$ is the conditional
measure
$$\bn(\cdot\midd A)= \frac{\bn(\cdot \cap A)}{\bn(A)}.$$
Moreover, $T_A$ and $e_{T_A}$ are independent.
\end{corollary}

This corollary can be used to calculate various distributions under the 
It\^o excursion measure. The distribution of the height and
the length of the excursion are given as follows:
For every $\ve >0$,
$$\bn\Big(\max_{t\geq 0} e(t)>\ve\Big) =\frac{1}{2\ve}$$
and 
$$\bn(\zeta(e) >\ve) = \frac{1}{\sqrt{2\pi \ve}}.$$
The It\^o excursion measure enjoys the following
scaling property. For every 
$\lambda>0$, define a mapping $\Phi_\lambda:E\la E$
by setting $\Phi_\lambda(e)(t)= \sqrt{\lambda}\,e(t/\lambda)$,
for every $e\in E$ and $t\geq 0$. Then we have $\Phi_\lambda(\bn)=\sqrt{\lambda}
\,\bn$.

This scaling property is useful when defining conditional versions 
of the It\^o excursion measure. We discuss the conditioning of $\bn(de)$ 
with respect to the length $\zeta(e)$.
There exists
a unique collection $(\bn_{(s)},s>0)$
of probability measures on $E$ such that the
following properties hold:

\begin{enumerate}
\item[(i)] For every $s>0$, $\bn_{(s)}(\zeta=s)=1$.

\item[(ii)] For every $\lambda>0$ and $s>0$, 
we have $\Phi_\lambda(\bn_{(s)})=\bn_{(\lambda s)}$.

\item[(iii)] For every measurable subset $A$ of $E$,

$$\bn(A)=
\int_0^\infty \bn_{(s)}(A)\,\frac{ds}{2\sqrt{2\pi s^3}}.$$
\end{enumerate}

We may and will write
$\bn_{(s)}=\bn(\cdot\midd \zeta =s)$.
The measure $\bn_{(1)}=\bn(\cdot\midd \zeta =1)$ is called the law of the
normalized Brownian excursion.  

There are many different descriptions of the It\^o excursion measure: See in
particular \cite[Chapter XII]{RY}. We state the following proposition,
which emphasizes the Markovian properties of $\bn$.
For every $t>0$ and $x> 0$, we set
$$q_t(x)=\frac{x}{\sqrt{2\pi t^3}} \exp(-\frac{x^2}{2t}).$$
Note that the function $t\mapsto q_t(x)$ is the density of the first hitting time of $x$ by $B$.
For $t>0$ and $x,y\in \R$, we also let 
$$p_t(x,y)= \frac{1}{\sqrt{2\pi t}} \exp(-\frac{(y-x)^2}{2t})$$
be the usual Brownian transition density.

\begin{proposition}
\label{ItoMarkov}
The It\^o excursion measure $\bn$ is the only $\sigma$-finite measure on $E$
that satisfies the following two properties:
\begin{enumerate}
\item[(i)] For every $t>0$, and every $f\in C(\R_+,\R_+)$,
$$\bn(f(e(t))\,\ind{\{\zeta>t\}}) = \int_0^\infty f(x)\,q_t(x)\,dx.$$
\item[(ii)] Let $t>0$. Under the conditional probability measure $\bn(\cdot\midd \zeta>t)$, the process $(e(t+r))_{r\geq 0}$ is Markov with the transition 
kernels of Brownian motion stopped upon hitting $0$. 
\end{enumerate}
\end{proposition}

This proposition can be used to establish absolute continuity properties
of the conditional measures $\bn_{(s)}$ with respect to $\bn$. For every $t\geq 0$, let ${\mathcal F}_t$ denote the 
$\sigma$-field on $E$ generated by the mappings $r\mapsto e(r)$, for $0\leq r\leq t$.
Then, if $0<t<1$, the measure $\bn_{(1)}$ is absolutely continuous
with respect to $\bn$
on the $\sigma$-field ${\mathcal F}_t$, with Radon-Nikod\'ym density
$$\frac{d\bn_{(1)}}{d\bn}\Big|_{{\mathcal F}_t}(e)= 2\sqrt{2\pi} \,q_{1-t}(e(t)).
$$
This formula provides a simple derivation of the finite-dimensional marginals under $\bn_{(1)}$, 
noting that the finite-dimensional marginals under $\bn$ are easily obtained from Proposition \ref{ItoMarkov}.
More precisely, for every integer $p\geq 1$, and every 
choice of $0<t_1<t_2<\cdots <t_p<1$, we get that the distribution of $(e(t_1),\ldots,e(t_p))$
under $\bn_{(1)}(de)$ has density
\begin{equation}
\label{finite-dim-excu}
2\sqrt{2\pi}\, q_{t_1}(x_1)\,p^*_{t_2-t_1}(x_1,x_2)\,p^*_{t_3-t_2}(x_2,x_3)\cdots
p^*_{t_p-t_1}(x_{p-1},x_p)\,q_{1-t_p}(x_p)
\end{equation}
where 
$$p^*_t(x,y)= p_t(x,y)-p_t(x,-y)\ ,t>0\;,\ x,y>0$$
is the transition density of Brownian motion killed when it hits $0$. As a side remark, formula (\ref{finite-dim-excu})
shows that the law of $(e(t))_{0\leq t\leq 1}$ under $\bn_{(1)}$ is invariant under time-reversal.

\subsection{Convergence of contour functions to the Brownian excursion}
\label{SSconvcontour}

The following theorem can be viewed as a special case of the results in
Aldous \cite{Al3}. The space of all continuous functions from $[0,1]$
into $\R_+$ is denoted by $C([0,1],\R_+)$, and is equipped with the
topology of uniform convergence.

\begin{theorem}
\label{Aldous}
For every integer $k\geq 1$, let $\theta_k$ be a random tree that is uniformly 
distributed over ${\bf A}_k$, and let $(C_k(t))_{t\geq 0}$
be its contour function. Then
$$\Big(\frac{1}{\sqrt{2k}} C_k(2k\,t)\Big)_{0\leq t\leq 1}
\build{\longrightarrow}_{k\to\infty}^{\rm(d)} (\ee_t)_{0\leq t\leq 1}$$
where $\ee$ is distributed according to $\bn_{(1)}$ (i.e.  
$\ee$ is a normalized Brownian excursion) and 
the convergence holds in the sense of weak convergence of the laws on
the space $C([0,1],\R_+)$.
\end{theorem}

\proof We already noticed that $\Pi_{\mu_0}(\cdot\mid |\tau|=k)$
coincides with the uniform distribution over ${\bf A}_k$. By combining
this with Proposition \ref{Dyck-excursion}, we get that $(C_k(0),C_k(1),$ \hfill $\ldots, C_k({2k}))$	%CMI
 is distributed as an excursion of simple random walk conditioned to have 
length $2k$. Recall our notation $(S_n)_{n\geq 0}$
for simple random walk on $\Z$ starting from $0$, and $T=\inf\{n\geq 0: S_n=-1\}$. To get the desired result, we need to verify that the law
of 
$$\Big( \frac{1}{\sqrt{2k}} S_{\lfloor 2k t\rfloor}\Big)_{0\leq t\leq 1}$$
under $P(\cdot\mid T=2k+1)$ converges to $\bn_{(1)}$ as $k\to\infty$.
This result can be seen as a conditional version of Donsker's theorem
(see Kaigh \cite{Kai} for similar statements). We will provide a detailed proof, because
this result plays a major role in what follows, and because some of the
ingredients of the proof will be needed again in Section 3 below. As usual,
the proof is divided into two parts: We first check the convergence of finite-dimensional
marginals, and then establish the tightness of the sequence of laws.

\medskip
\noindent{\it Finite-dimensional marginals.} We first consider one-dimensional
marginals. So we fix $t\in(0,1)$, and we will verify that
\begin{eqnarray}
\label{1-dim-margi}
\lim_{k\to\infty} \sqrt{2k}\,P\Big(
  S_{\lfloor2kt\rfloor}=\lfloor x\sqrt{2k}\rfloor \hbox{ or }  \lfloor
  x\sqrt{2k}\rfloor+1\,\Big|\, T= 2k+1\Big)
\\= 4\sqrt{2\pi}\,q_t(x)\,q_{1-t}(x),\nonumber
\end{eqnarray}
uniformly when $x$ varies over a compact subset of $(0,\infty)$.
Comparing with the case $p=1$ of formula (\ref{finite-dim-excu}), we see that the law of $(2k)^{-1/2}S_{\lfloor2kt\rfloor}$ under $P(\cdot\mid T=2k+1)$
converges to the law of $e(t)$ under $n_{(1)}(de)$ (we even get a local version of this
convergence).

In order to prove (\ref{1-dim-margi}), we will use two lemmas. The first one is 
a very special case of classical local limit theorems (see e.g. Chapter 2 of Spitzer \cite{Spi}).

\begin{lemma}
\label{local-limit}
For every $\eps>0$,
$$\lim_{n\to\infty} \sup_{x\in\R} \sup_{s\geq \eps}\,\Big|\sqrt{n} P\Big(S_{\lfloor ns\rfloor}=\lfloor x\sqrt{n}\rfloor\hbox{ or }
\lfloor x\sqrt{n}\rfloor+1\Big) - 2\,p_s(0,x)\Big| = 0.
$$
\end{lemma}

In our special situation, the result of the lemma is easily obtained by direct calculations using
the explicit form of the law of $S_n$ and Stirling's formula.

The next lemma is (a special case of) a famous formula of Kemperman (see e.g. \cite{Pit} Chapter 6).
For every integer $\ell\in\Z$, we use $P_\ell$ for a probability measure under which the simple
random walk $S$ starts from $\ell$.

\begin{lemma}
\label{Kemp-form}
For every $\ell\in\Z_+$ and every integer $n\geq 1$,
$$P_\ell( T=n)=\frac{\ell+1}{n}\;P_\ell(S_n=-1).$$
\end{lemma}

\proof It is easy to see that
$$P_\ell( T=n) =\frac{1}{2}P_\ell(S_{n-1}=0,\,T>n-1).$$
On the other hand,
\begin{align*}
P_\ell(S_{n-1}=0,\,T>n-1)&=P_\ell(S_{n-1}=0)-P_\ell(S_{n-1}=0,\,T\leq n-1)\\
&=P_\ell(S_{n-1}=0)-P_\ell(S_{n-1}=-2,\,T\leq n-1)\\
&=P_\ell(S_{n-1}=0)-P_\ell(S_{n-1}=-2),
\end{align*}
where the second equality is a simple application of the reflection principle.
So we have
$$P_\ell( T=n) =\frac{1}{2}\Big(P_\ell(S_{n-1}=0)-P_\ell(S_{n-1}=-2)\Big)$$
and an elementary calculation shows that this is equivalent to the statement
of the lemma. \cq

\smallskip
Let us turn to the proof of (\ref{1-dim-margi}). We first write for $i\in\{1,\ldots,2k\}$
and $\ell\in\Z_+$,
$$P( S_i=\ell\mid T=2k+1)= \frac{ P( \{S_i=\ell\}\cap\{T=2k+1\})}{P(T=2k+1)}.$$
By an application of the Markov property of $S$,
$$P( \{S_i=\ell\}\cap\{T=2k+1\})=P(S_i=\ell, T>i)\,P_\ell(T=2k+1-i).$$
Furthermore, a simple time-reversal argument (we leave the details to the
reader) shows that
$$P(S_i=\ell, T>i) =2 \,P_\ell(T=i+1).$$
Summarizing, we have obtained
\begin{eqnarray}
\label{tech-margi}
\lefteqn{P( S_i=\ell\mid
  T=2k+1)=\frac{2P_\ell(T=i+1)P_\ell(T=2k+1-i)}{P(T=2k+1)}}\\
&&=\frac{2(2k+1)(\ell+1)^2}{(i+1)(2k+1-i)}  \; \frac{P_\ell(S_{i+1}=-1)P_\ell(S_{2k+1-i}=-1)}{P(S_{2k+1}=-1)}\nonumber
\end{eqnarray}
using Lemma \ref{Kemp-form} in the second equality.

We apply this identity with $i=\lfloor 2k t\rfloor $ and $\ell=\lfloor x\sqrt{2k}\rfloor $ or $\ell=\lfloor x\sqrt{2k}\rfloor +1$.
Using Lemma \ref{local-limit}, we have first
$$\frac{2(2k+1)(\lfloor x\sqrt{2k}\rfloor +1)^2}{(\lfloor 2kt\rfloor +1)(2k+1-\lfloor 2kt\rfloor )}\times \frac{1}{P(S_{2k+1}=-1)}
\approx2 \sqrt{2\pi}\,(k/2)^{1/2}\, \frac{x^2}{t(1-t)}$$
and, using Lemma \ref{local-limit} once again,
\begin{eqnarray*}
\lefteqn{P_{\lfloor x\sqrt{2k}\rfloor }(S_{\lfloor 2k t\rfloor +1}=-1)P_{\lfloor x\sqrt{2k}\rfloor }(S_{2k+1-\lfloor 2k t\rfloor }=-1) }\\
&+& 
P_{\lfloor x\sqrt{2k}\rfloor +1}(S_{\lfloor 2k t\rfloor +1}=-1)P_{\lfloor x\sqrt{2k}\rfloor +1}(S_{2k+1-\lfloor 2k t\rfloor }=-1)\\&& \approx 2\, k^{-1}\,p_t(0,x)p_{1-t}(0,x).
\end{eqnarray*}
Putting these estimates together, and noting that $q_t(x)=(x/t)p_t(0,x)$, we arrive at (\ref{1-dim-margi}).

Higher order marginals can be treated in a similar way. Let us sketch the argument in the
case of two-dimensional marginals. We observe that, if $0<i<j<2k$ and if $\ell,m\in\Z_+$, 
we have, by the same arguments as above,
\begin{eqnarray*}
\lefteqn{P(S_i=\ell,S_j=m,T=2k+1)}\\
&=& 2\,P_\ell(T=i+1)\,P_\ell(S_{j-i}=m,T>j-i)\,P_m(T=k+1-j).
\end{eqnarray*}
Only the middle term $P_\ell(S_{j-i}=m,T>j-i)$ requires a different treatment than in
the case of one-dimensional marginals. However, by an application of the reflection
principle, one has 
$$P_\ell(S_{j-i}=m,T>j-i) = P_\ell(S_{j-i}=m)- P_\ell(S_{j-i}=-m-2).$$
Hence, using Lemma \ref{local-limit}, we easily obtain that for $x,y>0$ and 
$0<s<t<1$,
\begin{eqnarray*}
P_{\lfloor x\sqrt{2k}\rfloor }(S_{\lfloor 2k t\rfloor -\lfloor 2k s\rfloor }=\lfloor y\sqrt{2k}\rfloor )
+ P_{\lfloor x\sqrt{2k}\rfloor +1}(S_{\lfloor 2k t\rfloor -\lfloor 2k
  s\rfloor }=\lfloor y\sqrt{2k}\rfloor )\\
\approx (2k)^{-1/2}\,p^*_{t-s}(x,y),
\end{eqnarray*}
and the result for two-dimensional marginals follows in a straightforward way. 

\medskip
\noindent{\it Tightness.} We start with some combinatorial considerations. 
We fix $k\geq 1$. Let  $(x_0,x_1,\ldots,x_{2k})$be a Dyck path with length $2k$,
and let $i\in\{0,1,\ldots,2k-1\}$. We set, for every $j\in\{0,1,\ldots,2k\}$,
$$x^{(i)}_j= x_i + x_{i\oplus j} - 2\,\min_{i\wedge(i\oplus j)\leq n\leq i\vee (i \oplus j)} x_n$$
with the notation $i\oplus j= i+j$ if $i+j\leq 2k$, and
 $i\oplus j=i+j-2k$ if $i+j>2k$. It is elementary to see that 
 $(x^{(i)}_0,x^{(i)}_1,\ldots,x^{(i)}_{2k})$ is again a Dyck path with length $2k$.
 Moreover, the mapping $\Phi_i:(x_0,x_1,\ldots,x_{2k})\la (x^{(i)}_0,x^{(i)}_1,\ldots,x^{(i)}_{2k})$ 
 is a bijection from the set of all Dyck paths with length $2k$ onto itself. To see this,
 one may check that the composition $\Phi_{2k-i}\circ \Phi_i$ is the identity mapping.
 This property is easily verified by viewing $\Phi_i$ as a mapping defined on
 plane trees with $2k$ edges (using Proposition \ref{Dyck}): The plane tree corresponding
 to the image under $\Phi_i$ of the Dyck path associated with a tree $\tau$ is
 the ``same'' tree $\tau$ re-rooted at the corner corresponding to the $i$-th step of the contour exploration
 of $\tau$. From this observation it is obvious that the composition $\Phi_{2k-i}\circ\Phi_i$
leads us back to the original plane tree.
 
 To simplify notation,
we set for every $i,j\in\{0,1,\ldots,2k\}$,
$$\check C_k^{i,j}= \min_{i\wedge j\leq n\leq i\vee j} C_k(n).$$
The preceding discussion then gives the identity in distribution
\begin{equation}
\label{discrete-reroot}
\Big(C_k(i)+ C_k(i\oplus j) - 2 \check C_k^{i,i\oplus j}\Big)_{0\leq j\leq 2k}\build{=}_{}^{(d)}
(C_k(j))_{0\leq j\leq 2k}. 
\end{equation}

\begin{lemma}
\label{moment-contour}
For every integer $p\geq 1$, there exists a constant $K_p$
such that, for every $k\geq 1$ and every $i\in\{0,1,\ldots,2k\}$,
$$E[C_k(i)^{2p}]\leq K_p\,i^p.$$
\end{lemma}

Assuming that the lemma holds, the proof of tightness is easily completed. Using the
identity (\ref{discrete-reroot}), we get for $0\leq i<j\leq 2k$,
\begin{eqnarray*}
E[(C_k(j)-C_k(i))^{2p}]&\leq& E[(C_k(i)+ C_k(j) - 2 \check
C_k^{i,j})^{2p}]\\
&=&E[C_k(j-i)^{2p}]\\
&\leq & K_p(j-i)^p.
\end{eqnarray*}
It readily follows that the bound
$$E\Big[\Big(\frac{C_k(2kt)-C_k(2ks)}{\sqrt{2k}}\Big)^{2p}\Big]
\leq K_p\,(t-s)^p.$$
holds at least if $s$ and $t$ are of the form $s=i/2k$, $t=j/2k$, with $0\leq i<j\leq 2k$. Since the function $C_k$
is $1$-Lipschitz, a simple argument shows that the same bound holds (possibly
with a different constant $K_p$) whenever $0\leq s<t\leq 1$. This gives the
desired tightness, but we still have to prove the lemma.

\smallskip
\noindent{\bf Proof of Lemma \ref{moment-contour}.}
Clearly, we may restrict our attention to the case $1\leq i\leq k$ (note that
$(C_k(2k-i))_{0\leq i\leq 2k}$ has the same distribution as $(C_k(i))_{0\leq i\leq 2k}$).
Recall that $C_k(i)$ has the same distribution as $S_i$ under
$P(\cdot \mid T=2k+1)$. By formula (\ref{tech-margi}), we have thus,
for every integer $\ell\geq 0$,
$$P(C_k(i)=\ell)=\frac{2(2k+1)(\ell+1)^2}{(i+1)(2k+1-i)}  \; \frac{P_\ell(S_{i+1}=-1)P_\ell(S_{2k+1-i}=-1)}{P(S_{2k+1}=-1)}.$$
From Lemma \ref{local-limit} (and our assumption $i\leq k$), we can find two positive
constants $c_0$ and $c_1$ such that
$$P(S_{2k+1}=-1)\geq c_0(2k)^{-1/2}\ ,\quad P_\ell(S_{2k+1-i}=-1)\leq c_1(2k)^{-1/2}.$$
It then follows that
\begin{eqnarray*}
P(C_k(i)=\ell)&\leq& 4c_1(c_0)^{-1}\,\frac{(\ell +1)^2}{i+1}\, P_\ell(S_{i+1}=-1) 
\\&=& 4c_1(c_0)^{-1}\,\frac{(\ell +1)^2}{i+1}\, P(S_{i+1}=\ell+1).
\end{eqnarray*}
Consequently,
\begin{align*}
E[C_k(i)^{2p}]&=\sum_{\ell=0}^\infty \ell^{2p} P(C_k(i)=\ell)\\
&\leq \frac{4c_1(c_0)^{-1}}{i+1} \sum_{\ell=0}^\infty \ell^{2p}(\ell+1)^2\,P(S_{i+1}=\ell+1)\\
&\leq  \frac{4c_1(c_0)^{-1}}{i+1} \,E[(S_{i+1})^{2p+2}].
\end{align*}
However, it is well known and easy to prove that $E[(S_{i+1})^{2p+2}]\leq K'_p(i+1)^{p+1}$,
with some constant $K'_p$ independent of $i$. This completes the proof of the
lemma and of Theorem \ref{Aldous}.
\cq

\smallskip
Extensions and variants of Theorem \ref{Aldous} can be found in \cite{Al3}, \cite{Duq} and \cite{DuL0}.
To illustrate the power of this theorem, let us give a typical application.
The height $H(\tau)$ of a plane tree $\tau$ is the maximal generation of
a vertex of $\tau$.

\begin{corollary}
\label{height-tree}
Let $\theta_k$ be uniformly 
distributed over ${\bf A}_k$. Then 
$$\frac{1}{\sqrt{2k}} H(\theta_k) \build{\longrightarrow}_{k\to\infty}^{\rm(d)}
\max_{0\leq t\leq 1} \ee_t.$$
\end{corollary}

Since 
$$\frac{1}{\sqrt{2k}} H(\theta_k)= \max_{0\leq t\leq 1} \Big(\frac{1}{\sqrt{2k}} C_k(2k\,t)\Big)$$
the result of the corollary is immediate from Theorem \ref{Aldous}.

The limiting distribution in Corollary \ref{height-tree} is known in the form
of a series: For every $x>0$,
$$P\Big(\max_{0\leq t\leq 1} \ee_t >x\Big)
=2 \sum_{k=1}^\infty (4k^2x^2-1)\,\exp(-2k^2x^2).$$
See Chung \cite{Ch}.

\section{Real trees and the Gromov-Hausdorff convergence}

Our main goal in this section is to interpret the convergence of contour functions
in Theorem \ref{Aldous} as a convergence of discrete random trees towards 
a ``continuous random tree'' which is coded by the Brownian excursion in the 
same sense as a plane tree is coded by its contour function. We need to introduce a 
suitable notion of a continuous tree, and then to explain in which sense the
convergence takes place.

\subsection{Real trees}

We start with a formal definition. In these notes, we consider only
{\it compact} real trees, and so we include this compactness property
in the definition.

\begin{definition}\label{sec:real-trees}
A compact metric space $(\t,d)$ is a real tree if the following two
properties hold for every $a,b\in \t$.

\begin{enumerate}
\item[(i)] There is a unique
isometric map
$f_{a,b}$ from $[0,d(a,b)]$ into $\t$ such
that $f_{a,b}(0)=a$ and $f_{a,b}(
d(a,b))=b$.
\item[(ii)] If $q$ is a continuous injective map from $[0,1]$ into
$\t$, such that $q(0)=a$ and $q(1)=b$, we have
$$q([0,1])=f_{a,b}([0,d(a,b)]).$$
\end{enumerate}

A rooted real tree is a real tree $(\t,d)$
with a distinguished vertex $\rho=\rho(\t)$ called the root.
In what follows, real trees will always be rooted, even if this
is not mentioned explicitly. 
\end{definition}

Informally, one should think of a (compact) real tree as a 
connected union of line segments
in the plane with  no loops. Asssume for simplicity that there are finitely many
segments in the union. Then, for any two
points $a$ and $b$ in the tree, there is a unique path going from $a$ to $b$
in the tree, which is the concatentation of finitely many line segments. The distance between $a$ and
$b$ is then the length of this path.

Let us consider a rooted real tree $(\t,d)$.
The range of the mapping $f_{a,b}$ in (i) is denoted by
$\llbracket a,b\rrbracket$ (this is the ``line segment'' between $a$
and $b$ in the tree). 
In particular, $\llbracket \rho,a\rrbracket$ is the path 
going from the root to $a$, which we will interpret as the ancestral
line of vertex $a$. More precisely, we can define a partial order on the
tree by setting $a\preccurlyeq b$
($a$ is an ancestor of $b$) if and only if $a\in\llbracket \rho,b
\rrbracket$.

If $a,b\in\t$, there is a unique $c\in\t$ such that
$\llbracket \rho,a
\rrbracket\cap \llbracket \rho,b
\rrbracket=\llbracket \rho,c
\rrbracket$. We write $c=a\wedge b$ and call $c$ the most recent
common ancestor to $a$ and $b$.

By definition, the multiplicity of a vertex $a\in\t$ is the number of
connected components of $\t\backslash \{a\}$. Vertices
of $\t$ which have multiplicity
$1$ are called leaves.

\subsection{Coding real trees}
\label{SScodingrealtree}

In this subsection, we describe a method for constructing real trees,
which is well-suited to our forthcoming applications
to random trees.  This method is nothing but a continuous analog
of the coding of discrete trees by contour functions. 

We
consider a	 (deterministic) continuous function
$g:[0,1]\longrightarrow[0,\infty)$ 
such that $g(0)=g(1)=0$. To avoid
trivialities, we will also assume that $g$ is not identically zero.
For every $s,t\in[0,1]$, we set
$$m_g(s,t)=\inf_{r\in[s\wedge t,s\vee t]}g(r),$$
and
$$d_g(s,t)=g(s)+g(t)-2m_g(s,t).$$
Clearly $d_g(s,t)=d_g(t,s)$ and it is also easy to verify the triangle
inequality
$$d_g(s,u)\leq d_g(s,t)+d_g(t,u)$$
for every $s,t,u\in[0,1]$. We then introduce the equivalence relation
$s\sim t$ iff $d_g(s,t)=0$ (or equivalently iff $g(s)=g(t)=m_g(s,t)$). Let
$\t_g$ be the quotient space
$$\t_g=[0,1]/ \sim.$$
Obviously the function $d_g$ induces a distance on $\t_g$, and we keep the
notation $d_g$ for this distance. We denote by
$p_g:[0,1]\longrightarrow
\t_g$ the canonical projection. Clearly $p_g$ is continuous (when
$[0,1]$ is equipped with the Euclidean metric and $\t_g$ with the
metric $d_g$), and the metric space $(\t_g,d_g)$ is thus compact.

\begin{theorem}
\label{tree-deterministic}
The metric space $(\t_g,d_g)$ is a real tree.
We will view $(\t_g,d_g)$ as a rooted tree with root $\rho=p_g(0)=p_g(1)$.
\end{theorem}

\begin{remark} It is also possible to prove that any (rooted) real
tree can be represented in the form $\t_g$. We will leave this
as an exercise for the reader. 
\end{remark}

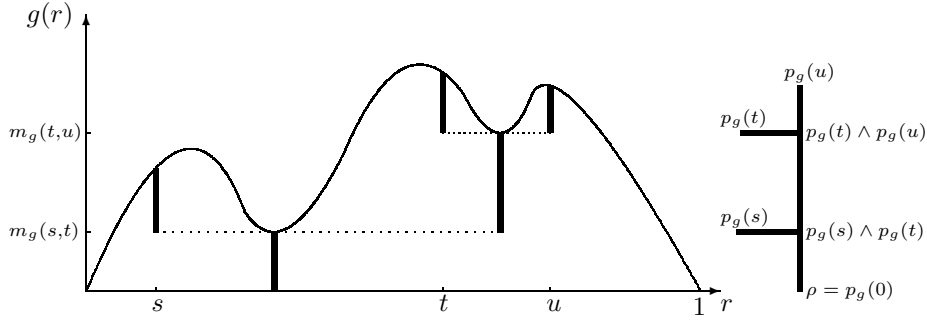
\begin{figure}
\begin{center}
\unitlength0.75pt
\begin{picture}(400,160)
\thinlines \put(10,10){\vector(1,0){320}}
\thinlines \put(10,10){\vector(0,1){140}}
\bezier{400}(10,10)(60,130)(90,50)
\bezier{400}(90,50)(110,20)(140,80)
\bezier{400}(140,80)(170,150)(200,110)
\bezier{400}(200,110)(220,70)(235,110)
\bezier{400}(235,110)(250,135)(320,10)
\bezier{40}(45,40)(133,40)(219,40)
\bezier{20}(190,90)(217,90)(244,90)
\linethickness{2pt}\put(45,40){\line(0,1){32}}
\put(105,10){\line(0,1){30}}
\put(219,40){\line(0,1){50}}
\put(190,90){\line(0,1){30}}
\put(244,90){\line(0,1){24}}
\thinlines\put(45,10){\line(0,1){2}}
\thinlines\put(190,10){\line(0,1){2}}
\thinlines\put(244,10){\line(0,1){2}}
\put(43,0){$s$}
\put(188,0){$t$}
\put(242,0){$u$}
\put(316,-2){$1$}
\put(330,0){$r$}
\put(-20,146){$g(r)$}
\put(-29,38){$\scriptstyle m_g(s,t)$}
\put(-29,88){$\scriptstyle m_g(t,u)$}
\thinlines\put(10,40){\line(1,0){2}}
\thinlines\put(10,90){\line(1,0){2}}
\linethickness{2pt}\put(370,10){\line(0,1){30}}
\linethickness{2pt}\put(370,40){\line(-1,0){32}}
\put(370,40){\line(0,1){50}}
\put(370,90){\line(-1,0){30}}
\put(370,90){\line(0,1){24}}
\put(373,38){\scriptsize$p_g(s)\wedge p_g(t)$}
\put(373,88){\scriptsize$p_g(t)\wedge p_g(u)$}
\put(373,8){\scriptsize$\rho=p_g(0)$}
\put(330,45){\scriptsize$p_g(s)$}
\put(330,95){\scriptsize$p_g(t)$}
\put(362,118){\scriptsize$p_g(u)$}

\end{picture}

\caption{Coding a tree by a continuous function}
\end{center}
\end{figure}

To get an intuitive understanding of Theorem \ref{tree-deterministic}, the 
reader should have a look at Fig.3. This figure shows how to construct a
simple subtree of $\t_g$, namely the
``reduced tree'' consisting of the union of the ancestral lines in $\t_g$ of three
vertices $p_g(s),p_g(t),p_g(u)$ corresponding to three (given)
times $s,t,u\in[0,1]$. This reduced tree is the union of the five bold line
segments that are constructed from the graph of $g$ in the way explained on
the left part of the figure. Notice that the lengths of the horizontal dotted lines play no
role in the construction, and that the reduced tree should be viewed as pictured on the
right part of Fig.3. The ancestral line 
of $p_g(s)$ (resp. $p_g(t),p_g(u)$) is a line segment
of length $g(s)$ (resp. $g(t),g(u)$). The ancestral lines of $p_g(s)$
and $p_g(t)$ share a common part, which has length $m_g(s,t)$ (the line
segment at the bottom in the left or the right part of Fig.3), and of course 
a similar property holds for the ancestral lines of $p_g(s)$
and $p_g(u)$, or of $p_g(t)$
and $p_g(u)$.

The following re-rooting lemma, which is of independent interest, is a useful
ingredient of the proof of Theorem \ref{tree-deterministic} (a discrete version of 
this lemma already appeared at the beginning of the proof of tightness
in Theorem \ref{Aldous}).

\begin{lemma}
\label{root-change}
Let $s_0\in[0,1)$. For any real $r\geq 0$, denote the fractional part of $r$
by
$\ov r=r-\lfloor r\rfloor$. Set
$$g'(s)=g(s_0)+g(\ov{s_0+s})-2m_g(s_0,\ov{s_0+s}),$$
for every $s\in [0,1]$. Then,
the function $g'$ is continuous and satisfies
$g'(0)=g'(1)=0$, so that we can define
$\t_{g'}$. Furthermore, for
every
$s,t\in[0,1]$, we have
\begin{equation}
\label{iso-root}
d_{g'}(s,t)=d_g(\ov{s_0+s},\ov{s_0+t})
\end{equation}
and there exists a unique isometry $R$ from $\t_{g'}$ onto $\t_{g}$
such that, for every $s\in [0,1]$,
\begin{equation}
\label{def-changeroot}
R(p_{g'}(s))=p_g(\ov{s_0+s}).
\end{equation}
\end{lemma}

Assuming that Theorem \ref{tree-deterministic} is proved, we see that
$\t_{g'}$ coincides with the real tree $\t_g$ re-rooted at $p_g(s_0)$.
Thus the lemma tells us which function codes the tree $\t_g$ re-rooted
at an arbitrary vertex.

\smallskip
\noindent{\bf Proof.} It is immediately checked that $g'$ satisfies
the same assumptions as $g$, so that we
can make sense of $\t_{g'}$. Then the key step is to verify
the relation (\ref{iso-root}). Consider first the case where
$s,t\in[0,1-s_0)$. Then two possibilities may occur.

If $m_g(s_0+s,s_0+t)\geq m_g(s_0,s_0+s)$, then
$m_g(s_0,s_0+r)=m_g(s_0,s_0+s)=m_g(s_0,s_0+t)$ for every $r\in[s,t]$, and
so
$$m_{g'}(s,t)=g(s_0)+m_g(s_0+s,s_0+t)
-2m_g(s_0,s_0+s).$$
It follows that
\begin{align*}
d_{g'}(s,t)&=g'(s)+g'(t)-2m_{g'}(s,t)\\
&=g(s_0+s)-2m_g(s_0,s_0+s)+g(s_0+t)\\
&\qquad-2m_g(s_0,s_0+t)
-2(m_g(s_0+s,s_0+t)-2m_g(s_0,s_0+s))\\
&=g(s_0+s)+g(s_0+t)
-2m_g(s_0+s,s_0+t)\\
&=d_g(s_0+s,s_0+t).
\end{align*}

If $m_g(s_0+s,s_0+t)< m_g(s_0,s_0+s)$, then the minimum in the definition
of $m_{g'}(s,t)$ is attained at $r_1$ defined as the first $r\in[s,t]$
such that $g(s_0+r)=m_g(s_0,s_0+s)$ (because for $r\in[r_1,t]$ we will
have $g(s_0+r)-2m_g(s_0,s_0+r)\geq -m_g(s_0,s_0+r)\geq
-m_g(s_0,s_0+r_1)$). Therefore,
$$m_{g'}(s,t)=g(s_0)-m_g(s_0,s_0+s),$$
and
\begin{eqnarray*}
d_{g'}(s,t)
&=&g(s_0+s)-2m_g(s_0,s_0+s)+g(s_0+t)
\\ & &-2m_g(s_0,s_0+t)
+2m_g(s_0,s_0+s)\\
&=&d_g(s_0+s,s_0+t).
\end{eqnarray*}

The other cases are treated in a similar way and are left to the reader.

By (\ref{iso-root}), if $s,t\in[0,1]$ are such that $d_{g'}(s,t)=0$,
then
$d_g(\ov{s_0+s},\ov{s_0+t})=0$ so that $p_g(\ov{s_0+s})=p_g(\ov{s_0+t})$.
Noting that $\t_{g'}=p_{g'}([0,1])$, we can define
$R$ in a unique way by the relation (\ref{def-changeroot}). From
(\ref{iso-root}),
$R$ is an isometry, and it is also immediate that
$R$ takes $\t_{g'}$ onto $\t_g$. \cq

\medskip

Thanks to the lemma, the fact that $\t_g$  verifies property (i) in the 
definition of a real tree is obtained from the particular case when
$a=\rho$ and $b =p_g(s)$ for some $s\in[0,1]$.
In that case however, the isometric mapping $f_{\rho,b}$ is easily constructed 
by setting
$$f_{\rho,b}(t)=p_g(\sup\{r\leq s:g(r)=t\})\ ,\quad\hbox{for every }
0\leq t\leq g(s)=d_g(\rho,b).$$
The remaining part of the argument is straightforward: See Section 2 in \cite{DuL}.

\begin{remark}
A short proof of Theorem \ref{tree-deterministic} using the characterization of real trees
via the so-called four-point condition can be found in \cite{EW}.
\end{remark}

The following simple observation will be useful in Section 7: If $s,t\in[0,1]$, the 
line segment $\llbracket p_g(s),p_g(t )\rrbracket$ in the tree $\t_g$ coincides with the
collection of the vertices $p_g(r)$, for all $r\in[0,1]$ such that either
$g(r)=m_g(r,s)\geq m_g(s,t)$ or $g(r)=m_g(r,t)\geq m_g(s,t)$. This easily
follows from the construction of the distance $d_g$. 

\subsection{The Gromov-Hausdorff convergence}

In order to make sense of the convergence of 
discrete trees towards real trees, we will use the Gromov-Hausdorff distance
between compact metric spaces, which
has been introduced by Gromov (see e.g. \cite{Gro}) in view of
geometric applications.

\sm

If $(E,\delta)$ is a metric space, the
notation $\delta_{Haus}(K,K')$ stands for the usual Hausdorff metric between
compact subsets of $E$ :
$$\delta_{Haus}(K,K')=
\inf\{\varepsilon>0:K\subset U_\varepsilon(K')\hbox{ and }
K'\subset U_\varepsilon(K)\},$$
where $U_\varepsilon(K):=\{x\in E:\delta(x,K)\leq \varepsilon\}$.

A pointed metric space is just a pair  consisting of a metric space $E$ and a distinguished 
point $\rho$ of $E$. We often write $E$ instead of $(E,\rho)$ to simplify notation.

Then, if $(E_1,\rho_1)$ and $(E_2,\rho_2)$ are two 
pointed compact metric spaces, we define the distance
$d_{GH}(E_1,E_2)$ by
$$d_{GH}(E_1,E_2)=\inf\{\delta_{Haus}(\varphi_1(E_1),\varphi_2(E_2))\vee
\delta(\varphi_1(\rho_1),
\varphi_2(\rho_2))\}$$ 
where the infimum is over all possible choices of the metric space
$(E,\delta)$ and the isometric
embeddings
$\varphi_1:E_1\la E$ and
$\varphi_2:E_2\la E$ of $E_1$ and $E_2$ into $E$. 

Two pointed compact metric spaces $E_{1}$ and $E_{2}$ are called equivalent if
there is an isometry that maps 
$E_{1}$ onto $E_{2}$ and preserves the distinguished points. Obviously
$d_{GH}(E_1,E_2)$ only depends on the equivalence classes of $E_1$ and
$E_2$. We denote by $\K$ the space of all equivalence classes of pointed compact metric spaces.

\begin{theorem}
\label{Burago}
$d_{GH}$
defines  a metric on the set $\K$.
Furthermore the metric space 
$(\K,d_{GH})$
is separable and complete.
\end{theorem}
 
A proof of the fact that $d_{GH}$
is  a metric on the set $\K$ can be found in \cite[Theorem 7.3.30]{BBI}. This proof is 
in fact concerned with the non-pointed case, but the argument is easily adapted to
our setting. The separability of the space $(\K,d_{GH})$ follows from the fact that
finite metric spaces are dense in $\K$. Finally the completeness of $(\K,d_{GH})$
can be obtained as a consequence of the compactness theorem 
in \cite[Theorem 7.4.15]{BBI}. 

In our applications, it will be important to have the following alternative definition of $d_{GH}$.
First recall that if
$(E_1,d_1)$ and $(E_2,d_2)$ are two compact metric spaces, a
correspondence between $E_1$ and $E_2$ is a subset $\mathcal{R}$ of
$E_1\times E_2$ such that for every $x_1\in E_1$ there exists at least one
$x_2\in E_2$ such that $(x_1,x_2)\in\mathcal{R}$ and conversely
for every $y_2\in E_2$ there exists at least one
$y_1\in E_1$ such that $(y_1,y_2)\in\mathcal{R}$. The distortion of
the correspondence $\mathcal{R}$ is defined by
$${\rm dis}(\mathcal{R})=\sup\{|d_1(x_1,y_1)-d_2(x_2,y_2)|:
(x_1,x_2),(y_1,y_2)\in\mathcal{R}\}.$$

\begin{proposition}
\label{corresp}
Let $(E_1,\rho_1)$ and $(E_2,\rho_2)$ be two pointed compact metric spaces. Then,
\begin{equation}
\label{bounddist}
d_{GH}(E_1,E_2)=\frac{1}{2}\ \inf_{\mathcal{R}\in\mathcal{C}(E_1,E_2),\,(\rho_1,\rho_2)\in\mathcal{R}}\,{\rm dis}(\mathcal{R}),
\end{equation}
where $\mathcal{C}(E_1,E_2)$ denotes the set of all correspondences 
between $E_1$ and $E_2$.
\end{proposition}

See \cite[Theorem 7.3.25]{BBI} for a proof of this proposition in the non-pointed case,
which is easily adapted.

The following consequence of Proposition \ref{corresp} will be very useful. Notice
that a rooted real tree can be viewed as a pointed compact metric space, whose
distinguished point is the root.

\begin{corollary}
\label{dist-trees}
Let $g$ and $g'$ be two continuous functions 
from $[0,1]$ into $\R_+$, such that $g(0)=g(1)=g'(0)=g'(1)=0$.
Then,
$$d_{GH}(\t_g,\t_{g'})\leq 2\|g-g'\|,$$
where $\|g-g'\|=\sup_{t\in[0,1]}|g(t)-g'(t)|$ is the supremum norm of $g-g'$.
\end{corollary}

\noindent{\bf Proof.} We rely 
on formula (\ref{bounddist}). We can construct a correspondence
between 
$\t_g$ and $\t_{g'}$ by setting
$$\mathcal{R}=\{(a,a'):\exists t\in[0,1]\hbox{ such that }a=p_g(t)\hbox{ and }
a'=p_{g'}(t)\}.$$
Note that $(\rho,\rho')\in \mathcal{R}$, if $\rho=p_g(0)$, resp. $\rho'=p_{g'}(0)$,
is the root of $\t_g$, resp. the root of $\t_{g'}$.
In order to bound the distortion of $\mathcal{R}$, let
$(a,a')\in\mathcal{R}$ and $(b,b')\in\mathcal{R}$. By
the definition of $\mathcal{R}$ we can find $s,t\geq 0$ such that
$p_g(s)=a$, $p_{g'}(s)=a'$ and $p_g(t)=b$, $p_{g'}(t)=
b'$. Now recall that
\begin{align*}
d_g(a,b)&=g(s)+g(t)-2m_g(s,t),\\
d_{g'}(a',b')&=g'(s)+g'(t)-2m_{g'}(s,t),
\end{align*}
so that
$$|d_g(a,b)-d_{g'}(a',b')|\leq 4\|g-g'\|.$$
Thus we have ${\rm dis}(\mathcal{R})\leq 4\|g-g'\|$ and the desired result
follows from (\ref{bounddist}). \cq

\subsection{Convergence towards the CRT}
\label{convCRT}

As in  subsection \ref{SSconvcontour}, we use the notation $\ee$ for a
normalized Brownian excursion. We view $\ee=(\ee_t)_{0\leq t\leq 1}$ as a 
(random) continuous function
over the interval $[0,1]$, which satisfies the same assumptions as the
function $g$ in subsection \ref{SScodingrealtree}.

\begin{definition}
The Brownian continuum random tree, also called the CRT, is the random real tree $\t_{\ee}$
coded by the normalized Brownian excursion.
\end{definition}

The CRT $\t_{\ee}$ is thus a random variable taking values in the set $\K$. 
Note that the measurability of this random variable follows from 
Corollary \ref{dist-trees}.

\begin{remark} Aldous \cite{Al1},\cite{Al3} uses a different method to define the
CRT. The preceding definition then corresponds to Corollary 22 in \cite{Al3}.
Note that our normalization differs by an unimportant scaling factor $2$
from the one in Aldous' papers: The CRT there is the tree $\t_{2\ee}$
instead of $\t_{\ee}$.
\end{remark}

We will now restate Theorem \ref{Aldous} as a convergence in distribution
of discrete random trees towards the CRT in the space $(\K,d_{GH})$.

\begin{theorem}
\label{approxCRT}
For every $k\geq 1$, let $\theta_k$ be uniformly distributed 
over ${\bf A}_k$, and equip $\theta_k$ with the usual graph distance $d_{gr}$. 
Then 
$$(\theta_k,(2k)^{-1/2}d_{gr})\build{\la}_{k\to\infty}^{\rm (d)} (\t_\ee,d_{\ee})$$
in the sense of convergence  in distribution for random variables with values in $(\K,d_{GH})$. 
\end{theorem}

\proof As in Theorem \ref{Aldous}, let $C_k$ be the contour function of 
$\theta_k$, and define a rescaled version of $C_k$ by setting
$$\wt C_k(t)= (2k)^{-1/2} C_k(2k\,t)$$
for every $t\in[0,1]$. Note that the function $\wt C_k$ is continuous and nonnegative 
over $[0,1]$ and vanishes at $0$ and at $1$. Therefore we can define the
real tree $\t_{\wt C_k}$.

Now observe that this real tree is very closely related to the (rescaled) discrete tree $\theta_k$.
Indeed $\t_{\wt C_k}$ is (isometric to) a finite union of 
line segments of length $(2k)^{-1/2}$ in the plane, with genealogical structure prescribed by
$\theta_k$, in the way suggested
in the left part of Fig.2.  From this observation, and the definition of the Gromov-Hausdorff
distance, we easily get
\begin{equation}
\label{compare-discrete}
d_{GH}\Big((\theta_k,(2k)^{-1/2}d_{gr}), (\t_{\wt C_k},d_{\wt C_k})\Big)\leq (2k)^{-1/2}.
\end{equation}
On the other hand, by combining Theorem \ref{Aldous} and Corollary \ref{dist-trees}, we have
$$
(\t_{\wt C_k},d_{\wt C_k})\build{\la}_{k\to\infty}^{\rm (d)} (\t_\ee,d_{\ee}).$$
The statement of Theorem \ref{approxCRT} now follows from the latter convergence
and (\ref{compare-discrete}). \cq

\begin{remark} 
Theorem \ref{approxCRT} contains in fact less information than Theorem \ref{Aldous},
because the lexicographical ordering that is inherent to the notion of a plane tree 
(and also to the coding of real trees by functions) disappears when we look
at a plane tree as a metric space. Still, Theorem \ref{approxCRT} is important
from the conceptual viewpoint: It is crucial to think of the CRT  as a
continuous limit of rescaled discrete random trees.
\end{remark}

There are analogs of Theorem \ref{approxCRT} for  
other classes of combinatorial trees. For instance, if $\tau_n$ is distributed uniformly
among all rooted Cayley trees with $n$ vertices, then 
$(\tau_n,(4n)^{-1/2}d_{gr})$ converges in distribution
to the CRT $\t_{\ee}$, in the space $\K$. Similarly, discrete random trees
that are uniformly distributed over binary trees with $2k$ edges converge
in distribution (modulo a suitable rescaling) towards the CRT. All
these results can be derived from a general statement of convergence
of 
conditioned Galton-Watson trees due to Aldous \cite{Al3} (see also
\cite{probasur}). A recent work of Haas and Miermont \cite{HM} provides
further extensions of Theorem \ref{approxCRT} to P\'olya trees (unordered
rooted trees).

\section{Labeled trees and the Brownian snake}

\subsection{Labeled trees}
\label{SSlabeledtree}

In view of forthcoming applications to random planar maps, we now introduce 
labeled trees. A labeled tree is a pair $(\tau, (\ell(v))_{v\in\tau})$ that consists
of a plane tree $\tau$ (see subsection \ref{SSdiscretetree}) and a collection
$(\ell(v))_{v\in\tau}$ of integer labels assigned to the vertices of $\tau$ -- in our 
formalism for plane trees, the tree $\tau$ coincides with the set of all its vertices. 
We assume that labels satisfy the following three properties:
\begin{enumerate}
\item[(i)] for every $v\in\tau$, $\ell(v)\in \Z\;$;
\item[(ii)] $\ell(\varnothing)=0\;$;
\item[(iii)] for every $v\in\tau\backslash\{\varnothing\}$, $\ell(v) - \ell({\pi(v)})=1,0,\hbox{ or }-1$,
\end{enumerate}
where we recall that $\pi(v)$ denotes the parent of $v$. Condition (iii) just means that when 
crossing an edge
of $\tau$ the label can change by at most $1$ in absolute value.

The motivation for introducing labeled trees comes from the fact that
(rooted and pointed) planar quadrangulations can be coded by such trees (see Section 4 below).
Our goal in the present section is to derive asymptotics for large labeled trees
chosen uniformly at random, in the same way as Theorem \ref{Aldous}, or
Theorem \ref{approxCRT}, provides asymptotics for large plane trees.
For every integer $k\geq 0$, we denote by $\bT_k$ the set of all
labeled trees with $k$ edges. It is immediate that
$$\#\bT_k= 3^k \#{\bf A}_k= \frac{3^k}{k+1}{2k\choose k}$$
simply because for each edge of the tree there are three possible choices 
for the label increment along this edge.

Let $(\tau, (\ell(v))_{v\in\tau})$ be a labeled tree with $k$ edges. As we saw in subsection \ref{SSdiscretetree},
the plane tree $\tau$ is coded by its contour function $(C_t)_{t\geq 0}$. We can similarly
encode the labels by another function $(V_t)_{t\geq 0}$, which is defined as follows.
If we explore the tree $\tau$ by following its contour, in the way suggested by the 
arrows of Fig.2, we visit successively all vertices of $\tau$ (vertices that are not leaves
are visited more than once). Write $v_0=\varnothing, v_1,v_2,\ldots, v_{2k}=\varnothing$ for the
successive vertices visited in this exploration. For instance, in the particular example 
of Fig.1 we have 
$$v_0=\varnothing, v_1=1, v_2=(1,1),v_3=1,v_4=(1,2),v_5=(1,2,1), v_6=(1,2),\ldots$$
The finite sequence $v_0, v_1,v_2,\ldots, v_{2k}$ will be called the contour exploration
of the vertices of $\tau$.

Notice that $C_i=|v_i|$, for every $i=0,1,\ldots,2k$, by the definition of the contour function.
We similarly set
$$V_i=\ell({v_i})\hbox{ \ for every }i=0,1,\ldots,2k.$$
To complete this definition, we set $V_t=0$ for $t\geq 2k$ and, for every $i=1,\ldots,2k$, we define
$V_t$ for $t\in(i-1,i)$ by using linear interpolation. We will call 
$(V_t)_{t\geq 0}$ the ``label contour function'' of the labeled tree $(\tau, (\ell(v))_{v\in\tau})$
Clearly  
$(\tau, (\ell(v))_{v\in\tau})$ is determined by the pair $(C_t,V_t)_{t\geq 0}$.

Our goal is now to describe the scaling limit of this pair when the labeled tree 
$(\tau, (\ell(v))_{v\in\tau})$ is chosen uniformly at random in $\bT_k$ and $k\to\infty$.
As an immediate consequence of Theorem \ref{Aldous} (and the fact that 
the number of possible labelings is the same for every plane tree with
$k$ edges), the scaling limit of $(C_t)_{t\geq 0}$ is the normalized 
Brownian excursion. To describe the scaling limit of $(V_t)_{t\geq 0}$ we need
to introduce the Brownian snake.

\subsection{The snake driven by a deterministic function}
\label{snakedeter}

Consider a continuous function $g:[0,1]\la \R_+$ 
such that $g(0)=g(1)=0$ (as 
in subsection \ref{SScodingrealtree}). We also assume that $g$ is H\"older
continuous: There exist two positive constants $K$ and $\gamma$
such that, for every $s,t\in[0,1]$,
$$|g(s)-g(t)|\leq K\,|s-t|^\gamma.$$

As in subsection \ref{SScodingrealtree}, we also set, 
for every $s,t\in[0,1]$,
$$m_g(s,t)=\min_{r\in[s\wedge t,s\vee t]}g(r).$$

\begin{lemma}
\label{posidefinite}
The function $(m_g(s,t))_{s,t\in[0,1]}$ is nonnegative definite in the sense
that, for every integer $n\geq 1$, for every $s_1,\ldots,s_n\in[0,1]$
and every $\lambda_1,\ldots,\lambda_n\in\R$, we have
$$\sum_{i=1}^n\sum_{j=1}^n \lambda_i \lambda_j \,m_g(s_i,s_j) \geq 0.$$
\end{lemma}

\proof Fix $s_1,\ldots,s_n\in[0,1]$, and let $t\geq 0$. For $i,j\in\{1,\ldots,n\}$, put $i\approx j$ 
if $m_g(s_i,s_j)\geq t$. Then $\approx$ is an equivalence relation on
$\{i:g(s_i)\geq t\}\subset \{1,\ldots,n\}$. By summing over the different classes
of this equivalence relation, we get that
$$\sum_{i=1}^n\sum_{j=1}^n \lambda_i \lambda_j \ind{\{t\leq m_g(s_i,s_j)\}} 
=\sum_{{\mathcal C}\;{\rm class\;of\;}\approx} \;\Big(\sum_{i\in{\mathcal C}} \lambda_i\Big)^2\geq 0.$$
Now integrate with respect to $dt$ to get the desired result. \cq

\smallskip
By Lemma \ref{posidefinite} and a standard application of the
Kolmogorov extension theorem, there exists a centered Gaussian process
$(Z^g_s)_{s\in[0,1]}$ whose covariance is
$$E[Z^g_sZ^g_t]=m_g(s,t)$$
for every $s,t\in[0,1]$. Consequently we have
\begin{eqnarray*}
E[(Z^g_s-Z^g_t)^2]&=&
  E[(Z^g_s)^2]+E[(Z^g_t)^2]-2E[Z^g_sZ^g_t]\\
&=& g(s)+g(t)-2m_g(s,t)\\
&\leq& 2K\,|s-t|^\gamma,
\end{eqnarray*}
where the last bound follows from our H\"older continuity assumption on $g$
(this calculation also shows that $E[(Z^g_s-Z^g_t)^2]=d_g(s,t)$, in the notation
of subsection \ref{SScodingrealtree}). From the previous bound and an
application of the Kolmogorov continuity criterion, the process $(Z^g_s)_{s\in[0,1]}$
has a modification with continuous sample paths. This leads us
to the following definition.

\begin{definition}
\label{detersnake}
The snake driven by the function $g$ is the centered Gaussian process
$(Z^g_s)_{s\in[0,1]}$ with continuous sample paths and covariance 
$$E[Z^g_sZ^g_t]=m_g(s,t)\ ,\quad s,t\in[0,1].$$
\end{definition}

Notice that we have in particular $Z^g_0=Z^g_1=0$. More generally, for every $t\in[0,1]$,
$Z^g_t$ is normal with mean $0$ and variance $g(t)$.

\begin{remark} Recall from subsection \ref{SScodingrealtree} the definition of the equivalence
relation $\sim$ associated with $g$: $s\sim t$ iff $d_g(s,t)=0$. Since we have $E[(Z^g_s-Z^g_t)^2]=d_g(s,t)$,
a simple argument shows that almost surely for every $s,t\in[0,1]$, the condition $s\sim t$
implies that $Z^g_s=Z^g_t$. In other words we may view $Z^g$ as a process indexed 
by the quotient $[0,1]\,/\!\sim$, that is by the tree $\t_g$. Indeed, it is then very natural to
interpret $Z^g$ as Brownian motion indexed by the tree $\t_g$: In the particular case
when $\t_g$ is a finite union of segments (which holds if $g$ is piecewise monotone),
$Z^g$ can be constructed by running independent Brownian motions along the
branches of $\t_g$.
It is however more convenient to view $Z^g$ as a process indexed by $[0,1]$
because later the function $g$ (and thus the tree $\t_g$) will be random and 
we avoid considering
a random process indexed by a random set. 
\end{remark}

\subsection{Convergence towards the Brownian snake}
\label{convSnake}

Let $\ee$ be as previously a normalized Brownian excursion. By standard properties 
of Brownian paths, the function $t\mapsto \ee_t$ is a.s. H\"older continuous (with
exponent $\frac{1}{2}-\ve$ for any $\ve >0$), and so we can apply the construction
of the previous subsection to (almost) every realization of $\ee$.

In other words, we can construct a pair $(\ee_t,Z_t)_{t\in [0,1]}$ of continuous random processes,
whose distribution is characterized by the following two properties:
\begin{enumerate}
\item[(i)] $\ee$ is a normalized Brownian excursion;
\item[(ii)] conditionally given $\ee$, $Z$ is distributed as the snake driven by $\ee$. 
\end{enumerate}

The process $Z$ will be called the Brownian snake (driven by $\ee$). This terminology is a little different
from the usual one: Usually, the Brownian snake is viewed as a path-valued process
(see e.g. \cite{Zurich}) and $Z_t$ would correspond only to the
terminal point of the value at time $t$ of this path-valued process. 

We can now answer the question raised at the end of subsection \ref{SSlabeledtree}. The following theorem
is due to Chassaing and Schaeffer \cite{CSise}. More general results can be found in \cite{JM}.

\begin{theorem}
\label{convsnake}
For every integer $k\geq 1$, let $(\theta_k,(\ell^k(v))_{v\in\theta_k})$ be distributed 
uniformly over the set $\bT_k$ of all labeled trees with $k$ edges. Let 
$(C_k(t))_{t\geq 0}$ and $(V_k(t))_{t\geq 0}$ be respectively the 
contour function and the label contour function of the 
labeled tree $(\theta_k,(\ell^k(v))_{v\in\theta_k})$. Then,
$$\Big( \frac{1}{\sqrt{2k}} C_k(2k\,t), \Big(\frac{9}{8k}\Big)^{1/4} V_k(2k\,t)\Big)_{t\in[0,1]}
\build{\la}_{k\to\infty}^{\rm(d)} (\ee_t,Z_t)_{t\in[0,1]}$$
where the convergence holds in the sense of weak convergence of the
laws on the space $C([0,1],\R_+^2)$. 
\end{theorem}

\proof   From Theorem \ref{Aldous} and the Skorokhod representation theorem, we may assume
without loss of generality that
\begin{equation}
\label{almosure}
\sup_{0\leq t\leq 1} |(2k)^{-1/2} C_k(2kt) - \ee_t| \build{\la}_{k\to\infty}^{\rm a.s.} 0.
\end{equation}

We first discuss the convergence of finite-dimensional marginals: We prove that
for every choice of $0\leq t_1<t_2<\cdots<t_p\leq 1$, we have 
\begin{equation}
\label{finitedim}
\Big( \frac{1}{\sqrt{2k}} C_k(2k\,t_i), \Big(\frac{9}{8k}\Big)^{1/4} V_k(2k\,t_i)\Big)_{1\leq i\leq p}
\build{\la}_{k\to\infty}^{\rm(d)} (\ee_{t_i},Z_{t_i})_{1\leq i\leq p}.
\end{equation}
Since for every $i\in\{1,\ldots,n\}$,
$$|C_k(2kt_i)-C_k(\lfloor 2kt_i\rfloor )|\leq 1\ ,\quad |V_k(2kt_i)-V_k(\lfloor 2kt_i\rfloor )|\leq 1$$
we may replace $2kt_i$ by its integer part $\lfloor 2kt_i\rfloor $ in (\ref{finitedim}).

Consider the case $p=1$. We may assume that $0<t_1<1$, because otherwise the
result is trivial. It is immediate that conditionally on $\theta_k$, the label
increments $\ell^k(v)-\ell^k({\pi(v)})$, $v\in\theta_k\backslash\{\emptyset\}$, are 
independent and uniformly distributed over $\{-1,0,1\}$. Consequently, we may write
$$(C_k(\lfloor 2kt_1\rfloor ),V_k(\lfloor 2kt_1\rfloor )) \build{=}_{}^{\rm(d)} \Big(C_k(\lfloor 2kt_1\rfloor ),
\sum_{i=1}^{C_k(\lfloor 2kt_1\rfloor )}\eta_i\Big)$$
where the variables $\eta_1,\eta_2,\ldots$ are independent and uniformly distributed over
$\{-1,0,$ \hfill $1\}$, and are also independent of the trees $\theta_k$. By the central limit theorem,    %CMI
$$\frac{1}{\sqrt{n}} \sum_{i=1}^n \eta_i \build{\la}_{n\to\infty}^{\rm(d)} \Big(\frac{2}{3}\Big)^{1/2}N$$
where $N$ is a standard normal variable. Thus
if we set for $\lambda\in \R$,
$$\Phi(n,\lambda)= E\Big[\exp\Big(i\frac{\lambda}{\sqrt{n}}\sum_{i=1}^n \eta_i\Big)\Big]$$
we have $\Phi(n,\lambda)\la \exp(-\lambda^2/3)$ as $n\to\infty$.

Then, for every $\lambda,\lambda'\in \R$, we get by conditioning on $\theta_k$
\begin{align*}
&E\Big[ \exp\Big(i\frac{\lambda}{\sqrt{2k}}C_k(\lfloor  2kt_1\rfloor )+i\frac{\lambda'}{\sqrt{C_k(\lfloor  2kt_1\rfloor )}}\sum_{i=1}^
{C_k(\lfloor  2kt_1\rfloor )} \eta_i\Big)\Big]\\
&\qquad 
= E\Big[\exp\Big(i\frac{\lambda}{\sqrt{2k}}C_k(\lfloor  2kt_1\rfloor )\Big)\times \Phi(C_k(\lfloor  2kt_1\rfloor ), \lambda')\Big]\\
&\qquad
\build{\la}_{k\to\infty}^{} E[\exp(i\lambda \ee_{t_1})]\times \exp(-\lambda'^2/3)
\end{align*}
using the (almost sure) convergence of $(2k)^{-1/2}C_k(\lfloor  2kt_1\rfloor )$ towards $\ee_{t_1}>0$. 
In other words we have obtained the joint convergence in distribution
\begin{equation}
\label{convt1}
\Big(\frac{C_k(\lfloor  2kt_1\rfloor )}{\sqrt{2k}}, \frac{1}{\sqrt{C_k(\lfloor  2kt_1\rfloor )}}\sum_{i=1}^{C_k(\lfloor  2kt_1\rfloor )} \eta_i \Big)\build{\la}_{k\to\infty}^{\rm(d)} (\ee_{t_1}, 
(2/3)^{1/2}N),
\end{equation}
where the normal variable $N$ is independent of $\ee$. 

From the preceding observations, we have
\begin{eqnarray*}
\lefteqn{\Big(\frac{C_k(\lfloor  2kt_1\rfloor
    )}{\sqrt{2k}},\!\Big(\frac{9}{8k}\Big)^{1/4}V_k(\lfloor
  2kt_1\rfloor )\Big)}\\
&\!\! \build{=}_{}^{\rm(d)}\!\!&
 \Big(\frac{C_k(\lfloor  2kt_1\rfloor )}
{\sqrt{2k}},\! \Big(\frac{3}{2}\Big)^{1/2} \Big(\frac{C_k(\lfloor  2kt_1\rfloor )}{\sqrt{2k}}\Big)
^{1/2}
\frac{1}{\sqrt{C_k(\lfloor  2kt_1\rfloor )}}\sum_{i=1}^
{C_k(\lfloor  2kt_1\rfloor )} \eta_i \Big)
\end{eqnarray*}
and from (\ref{convt1}) we get
$$\Big(\frac{C_k(\lfloor  2kt_1\rfloor )}{\sqrt{2k}},\Big(\frac{9}{8k}\Big)^{1/4}V_k(\lfloor  2kt_1\rfloor )\Big)
\build{\la}_{k\to\infty}^{\rm(d)} (\ee_{t_1}, 
\sqrt{\ee_{t_1}}\,N).$$
This gives (\ref{finitedim}) in the case $p=1$, since by construction
it holds that $(\ee_{t_1},Z_{t_1})\build{=}_{}^{\rm(d)}(\ee_{t_1}, 
\sqrt{\ee_{t_1}}\,N)$. 

Let us discuss the case $p=2$ of (\ref{finitedim}). We fix $t_1$ and $t_2$ with $0<t_1<t_2<1$. Recall the notation
$$\check C_k^{i,j}= \min_{i\wedge j\leq n\leq i\vee j} C_k(n)\;,\qquad i,j\in\{0,1,\ldots,2k\}$$
introduced in Section 1.
Write $v^k_0=\varnothing,v^k_1,\ldots,v^k_{2k}=\varnothing$
for the contour exploration of vertices of  $\theta_k$ (see the
 end of subsection \ref{SSlabeledtree}). Then we know that
 \begin{align*}
&C_k(\lfloor 2kt_1\rfloor )=|v^k_{\lfloor 2kt_1\rfloor }|,\;
C_k(\lfloor 2kt_2\rfloor )=|v^k_{\lfloor 2kt_2\rfloor }|, \; \\
&V_k(\lfloor 2kt_1\rfloor )=\ell^k({v^k_{\lfloor 2kt_1\rfloor }}),\;
V_k(\lfloor 2kt_2\rfloor )=\ell^k({v^k_{\lfloor 2kt_2\rfloor }}),
\end{align*}
and furthermore $\check C_k^{\lfloor 2kt_1\rfloor ,\lfloor
  2kt_2\rfloor }$ is the generation in $\theta_k$ of the last common
ancestor to $v^k_{\lfloor 2kt_1\rfloor }$ and $v^k_{\lfloor
  2kt_2\rfloor }$.  From the properties of labels on the tree
$\theta_k$, we now see that conditionally on $\theta_k$,
\begin{eqnarray}
\lefteqn{(V_k(\lfloor 2kt_1\rfloor ),V_k(\lfloor 2kt_2\rfloor
  )) \build{=}_{}^{\rm(d)} }\nonumber
\\
&&\Big( \sum_{i=1}^{\check C_k^{\lfloor 2kt_1\rfloor ,\lfloor
    2kt_2\rfloor }}\!\!\!\! \eta_i + \!\!\!\!  \sum_{i=\check C_k^{\lfloor 2kt_1\rfloor ,\lfloor 2kt_2\rfloor }+1}^{C_k(\lfloor 2kt_1\rfloor )} \!\!\!\! \eta'_i
\;,\;
\sum_{i=1}^{\check C_k^{\lfloor 2kt_1\rfloor ,\lfloor 2kt_2\rfloor
  }}\!\!\!\!  \eta_i + \!\!\!\! \sum_{i=\check C_k^{\lfloor
    2kt_1\rfloor ,\lfloor 2kt_2\rfloor }+1}^{C_k(\lfloor 2kt_2\rfloor
  )}\!\!\!\!  \eta''_i\Big) \label{convt2}
\end{eqnarray}
where the variables $\eta_i,\eta'_i,\eta''_i$ are independent and
uniformly distributed over $\{-1,0,$ \hfill $1\}$. 			%CMI

From (\ref{almosure}), we have
\begin{align*}
&\Big((2k)^{-1/2}C_k(\lfloor 2kt_1\rfloor ),(2k)^{-1/2}C_k(\lfloor 2kt_2\rfloor ),(2k)^{-1/2} \check C_k^{\lfloor 2kt_1\rfloor ,\lfloor 2kt_2\rfloor }\Big)\\
&\qquad \build{\la}_{k\to\infty}^{\rm a.s.}
(\ee_{t_1},\ee_{t_2},m_\ee(t_1,t_2)).
\end{align*}
By arguing as in the case $p=1$, we now deduce from (\ref{convt2}) that
\begin{align*}
&\Big(\frac{C_k(\lfloor 2kt_1\rfloor )}{\sqrt{2k}},\frac{C_k(\lfloor 2kt_2\rfloor )}{\sqrt{2k}},\Big(\frac{9}{8k}\Big)^{1/4}V_k(\lfloor 2kt_1\rfloor ),\Big(\frac{9}{8k}\Big)^{1/4}V_k(\lfloor 2kt_2\rfloor )\Big)\\
&\qquad\build{\la}_{k\to\infty}^{\rm(d)}
(\ee_{t_1},\ee_{t_2}, \sqrt{m_\ee(t_1,t_2)}\,N +
\sqrt{\ee_{t_1}-m_\ee(t_1,t_2)}\, N'\;,\\
&\qquad \qquad \qquad \qquad \quad\;\sqrt{m_\ee(t_1,t_2)}\,N + \sqrt{\ee_{t_2}-m_\ee(t_1,t_2)}\, N'')
\end{align*}
where $N,N',N''$ are three independent standard normal variables, which are also independent of $\ee$. The limiting
distribution in the last display is easily identified with that of $(\ee_{t_1},\ee_{t_2},Z_{t_1},Z_{t_2})$, and this gives the
case $p=2$ in (\ref{finitedim}). The general case is proved by similar arguments and we leave details to the reader.

To complete the proof of Theorem \ref{convsnake}, we need a tightness argument. The laws of the processes 
$$\Big( \frac{1}{\sqrt{2k}} C_k(2k\,t)\Big)_{t\in[0,1]}$$
are tight by Theorem \ref{Aldous}, and so we need only verify the tightness of the processes
$$\Big(\Big(\frac{9}{8k}\Big)^{1/4} V_k(2k\,t)\Big)_{t\in[0,1]}.$$
This is a consequence of the following lemma, which therefore 
completes the proof of Theorem \ref{convsnake}. \cq

\begin{lemma}
\label{tightness}
For every integer $p\geq 1$, there exists a constant $K_p<\infty$ such that,
for every $k\geq 1$ and every $s,t\in[0,1]$,
$$E\Big[\Big(\frac{V_k(2kt)- V_k(2ks)}{k^{1/4}}\Big)^{4p}\Big] \leq K_p \,|t-s|^p.$$
\end{lemma}

\proof Simple arguments show that we may restrict our attention to the case when 
$s=i/(2k)$, $t=j/(2k)$, with $i,j\in\{0,1,\ldots,2k\}$. By using the same
decomposition as in (\ref{convt2}), we have
\begin{equation}
\label{convt3}
V_k(j)-V_k(i) \build{=}_{}^{\rm(d)} \sum_{n=1}^{d_{gr}(v^k_i,v^k_j)} \eta_n
\end{equation}
where the random variables $\eta_n$ are independent and uniform over $\{-1,0,1\}$
(and independent of $\theta_k$) and
$$d_{gr}(v^k_i,v^k_j) = C_k(i)+ C_k(j) - 2 \check C_k^{i,j}$$
is the graph distance in the tree $\theta_k$ between vertices $v^k_i$ and $v^k_j$. 
From (\ref{convt3}) and by conditioning with respect to $\theta_k$, we get 
the existence of a constant $K'_p$ such that
$$E[(V_k(i)-V_k(j))^{4p}]\leq K'_p\, E[(d_{gr}(v^k_i,v^k_j))^{2p}].$$
So the lemma will be proved if we can verify the bound
\begin{equation}
\label{convt4}
E[(C_k(i)+ C_k(j) - 2 \check C_k^{i,j})^{2p}] \leq K''_p\,|j-i|^p
\end{equation}
with a constant $K''_p$ independent of $k$. 
By the identity (\ref{discrete-reroot}), it is enough to prove that this bound holds
for $i=0$. However, the case $i=0$ is exactly Lemma \ref{moment-contour}.
This completes the proof.

\section{Planar maps}\label{sec:planar-maps}

\subsection{Definitions}\label{sec:prem-defin-et}

A map is a combinatorial object, which can be best visualized as a
class of graphs embedded in a surface. In these lectures, we will
exclusively focus on the case of {\em plane} (or {\em planar}) maps,
where the surface is the 2-dimensional sphere $\S^2$.

Let us first formalize the notion of map. We will not enter into
details, referring the reader to the book by Mohar and Thomassen
\cite{MoTh01} for a very complete exposition. Another useful
reference, discussing in depth the different equivalent ways to define
maps (in particular through purely algebraic notions) is the book by
Lando and Zvonkin \cite[Chapter 1]{LaZv04}.

An {\em oriented edge} in $\S^2$ is a mapping $e:[0,1]\to \S^2$ that
is continuous, and such that either $e$ is injective, or the
restriction of $e$ to $[0,1)$ is injective and $e(0)=e(1)$. In the
  latter case, $e$ is also called a loop. An oriented edge will always
  be considered up to reparametrization by a continuous increasing
  function from $[0,1]$ to $[0,1]$, and we will always be interested
  in properties of edges that do not depend on a particular
  parameterization. The origin and target of $e$ are the points
  $e^-=e(0)$ and $e^+=e(1)$. The reversal of $e$ is the oriented edge
  $\ov{e}=e(1-\cdot)$. An {\em edge} is a pair $\be=\{e,\ov{e}\}$,
  where $e$ is an oriented edge. The {\em interior} of $\be$ is
  defined as $e((0,1))$.

An {\em embedded graph} in $\S^2$ is a graph\footnote{all the graphs
  considered here are finite, and are multigraphs in which multiple
  edges and loops are allowed} $G=(V,E)$ such that 
\begin{itemize}
\item
$V$ is a (finite) subset of $\S^2$
\item
$E$ is a (finite) set of edges in $\S^2$
\item
the vertices incident to $\be=\{e,\ov{e}\}\in E$ are $e^-,e^+\in V$
\item
the interior of an edge $\be\in E$ does not intersect $V$ nor the edges of
$E$ distinct from $\be$
\end{itemize}

The support of an embedded graph $G=(V,E)$ is 
$$\supp(G)=V\cup \bigcup_{\be=\{e,\ov{e}\}\in E}e([0,1])\, .$$ A {\em
  face} of the embedding is a connected component of the set
$\S^2\setminus\supp(G)$.

\begin{definition}\label{sec:definitions-3}
A (planar) map is a connected embedded graph. Equivalently, a map is
an embedded graph whose faces are all homeomorphic to the Euclidean
unit disk in $\R^2$.
\end{definition}

Topologically, one would say that a map is the $1$-skeleton of a
CW-complex decomposition of $\S^2$. We will denote maps using bold
characters $\bm,\bq,\ldots$

Let $\bm=(V,E)$ be a map, and let $\overrightarrow{E}=\{e\in \be:\be\in
E\}$ be the set of all oriented edges of $\bm$. Since $\S^2$ is oriented, it is
possible to define, for every oriented edge $e\in \overrightarrow{E}$,
a unique face $f_e$ of $\bm$, located to the left of the edge $e$. We
call $f_e$ the face {\em incident} to $e$. Note that the edges
incident to a given face form a closed curve in $\S^2$, but not
necessarily a Jordan curve (it can happen that $f_e=f_{\ov{e}}$ for
some $e$). The degree of a face $f$ is defined as
$$\deg(f)=\#\{e\in \overrightarrow{E}:f_e=f\}\, .$$ The oriented edges
incident to a given face $f$, are arranged cyclically in
counterclockwise order around the face in what we call the {\em facial
  ordering}. With every oriented edge $e$, we can associate a {\em
  corner} incident to $e$, which is a small simply connected
neighborhood of $e^-$ intersected with $f_e$. Then the
corners of two different oriented edges do not intersect.

Of course, the degree of a vertex $u\in V$ is the usual
graph-theoretical notion
$$\deg(u)=\#\{e\in \overrightarrow{E}:e^-=u\}\, .$$ Similarly as for
faces, the outgoing edges from $u$ are organized cyclically in
counterclockwise order around $u$.

\begin{figure}
\begin{center}
\includegraphics[scale=.8]{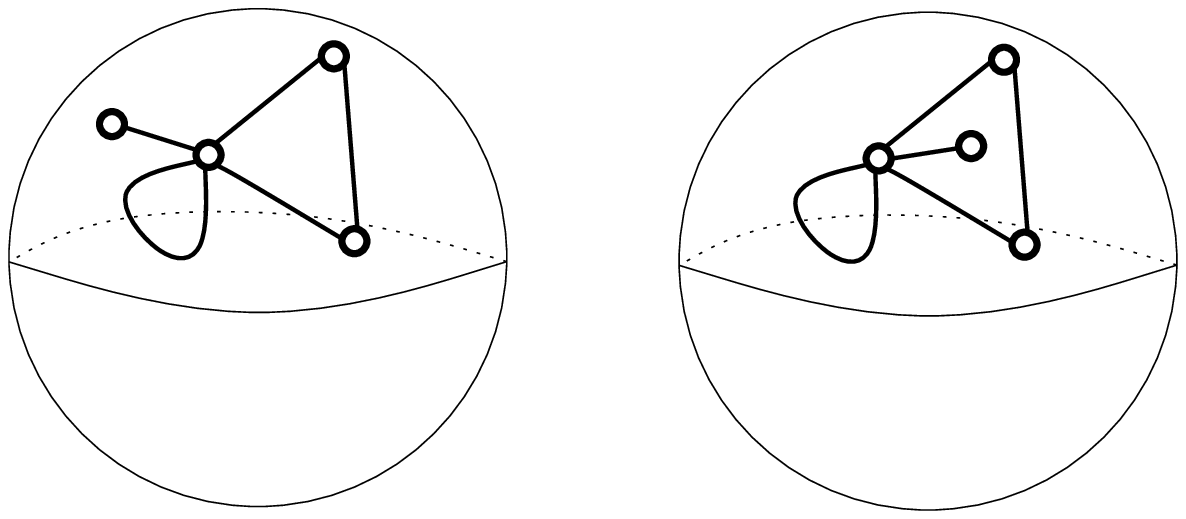}
\caption{Two planar maps, with 4 vertices and 3 faces of degrees 1,3,6
  and 1,4,5 respectively}
\label{fig:cartes}
\end{center}
\end{figure}   

A {\em rooted} map is a pair $(\bm,e)$ where $\bm=(V,E)$ is a map and
$e\in \overrightarrow{E}$ is a distinguished oriented edge, called the
root. We often omit the mention of $e$ in the notation. 

\subsection{Euler's formula}\label{sec:eulers-formula}

An important property of maps is the so-called {\em Euler
  formula}. If $\bm$ is a map, $V(\bm),E(\bm),F(\bm)$ denote respectively  the
sets of all vertices, edges and faces of $\bm$. Then,
\begin{equation}\label{eq:11}
\#V(\bm)-\#E(\bm)+\#F(\bm)=2\, .
\end{equation}
This is a relatively easy result in the case of interest (the planar
case): One can remove the edges of the graph one by one until a
spanning tree $\bt$ of the graph is obtained, for which the result is
trivial (it has one face, and $\#V(\bt)=\#E(\bt)+1$).

\subsection{Isomorphism, automorphism and rooting}\label{sec:isomorphisme}

In the sequel, we will always consider maps ``up to deformation'' in
the following sense. 

\begin{definition}\label{sec:isom-autom-sous-1}
The maps $\bm,\bm'$ on $\S^2$ are isomorphic if there exists an
orientation-preserving homeomorphism $h$ of $\S^2$ onto itself, such
that $h$ induces a graph isomorphism of $\bm$ with $\bm'$.

The rooted maps $(\bm,e)$ and $(\bm',e')$ are isomorphic if $\bm$ and
$\bm'$ are isomorphic through a homeomorphism $h$ that maps $e$ to
$e'$.
\end{definition}

In the sequel, we will almost always identify two isomorphic maps
$\bm,\bm'$. This of course implies that the (non-embedded,
combinatorial) graphs associated with $\bm,\bm'$ are isomorphic, but
this is stronger: For instance the two maps of Fig.\ref{fig:cartes}
are not isomorphic, since a map isomorphism preserves the
degrees of faces.

An {\em automorphism} of a map $\bm$ is an isomorphism of $\bm$ with
itself. It should be interpreted as a {\em symmetry} of the map. An
important fact is the following.

\begin{proposition}\label{sec:isom-autom-sous-2}
An automorphism of $\bm$ that fixes an oriented edge fixes all the
oriented edges.
\end{proposition}

Loosely speaking, the only automorphism of a rooted map is the
identity. This explains why rooting is an important tool in the
combinatorial study of maps, as it ``kills the symmetries''. The idea
of the proof of the previous statement is to see that if $e$ is fixed
by the automorphism, then all the edges incident to $e^-$ should also
be fixed (since an automorphism preserves the orientation). One can
thus progress in the graph (by connectedness) and show that all the
edges are fixed. 

In a rooted map, the face $f_e$ incident to the root edge $e$ is
often called the {\em external face}, or root face. The other faces
are called {\em internal}. The vertex $e^-$ is called the root vertex.

From now on, unless otherwise specified, all maps will be rooted. 

We end this presentation by introducing the notion of {\em graph
  distance} in a map $\bm$. A {\em chain} of length $k\geq 1$ is a
sequence $e_{(1)},\ldots,e_{(k)}$ of oriented edges in
$\overrightarrow{E}(\bm)$, such that $e_{(i)}^+=e_{(i+1)}^-$ for
$1\leq i\leq k-1$, and we say that the chain links the vertices
$e_{(1)}^-$ and $e_{(k)}^+$. We also allow, for every vertex $u\in
V(\bm)$, a chain with length $0$, starting and ending at $u$. The {\em
  graph distance} $d_\bm(u,v)$ between two vertices $u,v\in V(\bm)$
is the minimal $k$ such that there exists a chain with length $k$
linking $u$ and $v$. A chain with minimal length between two vertices
is called a {\em geodesic chain}.

\subsection{The Cori-Vauquelin-Schaeffer bijection}\label{cha:les-meth-biject}

Via the identification of maps up to isomorphisms the set of all maps
becomes a countable set. For instance, the set $\bM_n$ of all rooted maps
with $n$ edges is a finite set: The $2n$ oriented edges should be
organized around a finite family of polygons (the faces of the map),
and the number of ways to associate the boundary edges of these
polygons is finite. A natural question to ask is ``what is the
cardinality of $\bM_n$?''.

Tutte answered this question (and many other counting problems for
maps), motivated in part by the $4$-color problem. He developed a
powerful method, the ``quadratic method'', to solve the apparently
ill-defined equations for the generating functions of maps. For recent
developments in this direction, see the article by Bousquet-M\'elou and
Jehanne \cite{BoJe06}. The method, however, is a kind of ``black box''
which solves such counting problems without giving much extra
information about the structure of maps. One obtains
$$\#\bM_n=\frac{2}{n+2}3^n\cat_n\, ,$$ where
$\cat_n=\frac{1}{n+1}\binom{2n}{n}$ is the $n$-th Catalan number.  We
also mention the huge literature on the enumeration of maps using
matrix integrals, initiating in \cite{thooft,BrItPaZu}, which is
particularly popular in the physics literature. See \cite[Chapter
  4]{LaZv04} for an introduction to this approach.

Motivated by the very simple form of the formula enumerating $\bM_n$,
Cori and Vauquelin \cite{CoVa} gave in 1981 a bijective approach to
this formula. These approaches reached their full power with the
work of Schaeffer starting in his 1998 thesis \cite{schaeffer98}. We
now describe the bijective approach in the case of quadrangulations.

\subsubsection{Quadrangulations}\label{sec:quadrangulations}

A map $\bq$ is a quadrangulation if all its faces are of degree
$4$. We let $\bQ_n$ be the set of all (rooted) quadrangulations with $n$
faces. Quadrangulations are a very natural family of maps to consider,
in virtue of the fact that there exists a ``trivial'' bijection
between $\bM_n$ and $\bQ_n$, which can be described as follows. 

Let $\bm$ be a map with $n$ edges, and imagine that the vertices of
$\bm$ are colored in black. We then create a new map by adding inside
each face of $\bm$ a white vertex, and joining this white vertex to
every corner of the face $f$ it belongs to, by non-intersecting edges
inside the face $f$. In doing so, notice that some black vertices may
be joined to the same white vertex with several edges. Lastly, we
erase the interiors of the edges of the map $\bm$. We end up with a
map $\bq$, which is a plane quadrangulation with $n$ faces, each face
containing exactly one edge of the initial map. We adopt a rooting
convention, for instance, we root $\bq$ at the first edge coming after
$e$ in counterclockwise order around $e^-$, where $e$ is the root of
$\bm$.

Notice that $\bq$ also comes with a bicoloration of its vertices in
black and white, in which two adjacent vertices have different
colors. This says that $\bq$ is {\em bipartite}, and as a matter of
fact, every (planar!) quadrangulation is bipartite. So this coloring is
superfluous: One can recover it by declaring that the black vertices
are those at even distance from the root vertex of $\bq$, and the white
vertices are those at odd  distance from the root vertex.

Conversely, starting from a rooted quadrangulation $\bq$, we can
recover a bipartite coloration as above, by declaring that the
vertices at even distance from the root edge are black. Then, we draw
the diagonal linking the two black corners incident to every face of
$\bq$. Finally, we remove the interior of the edges of $\bq$ and root
the resulting map $\bm$ at the first outgoing diagonal from $e^-$ in
clockwise order from the root edge $e$ of $\bq$. One checks that this
is indeed a left- and right-inverse of the previous mapping from
$\bM_n$ to $\bQ_n$. See Fig.5 below for an illustration of these
bijections.

\begin{figure}

\begin{center}
\includegraphics[scale=.7]{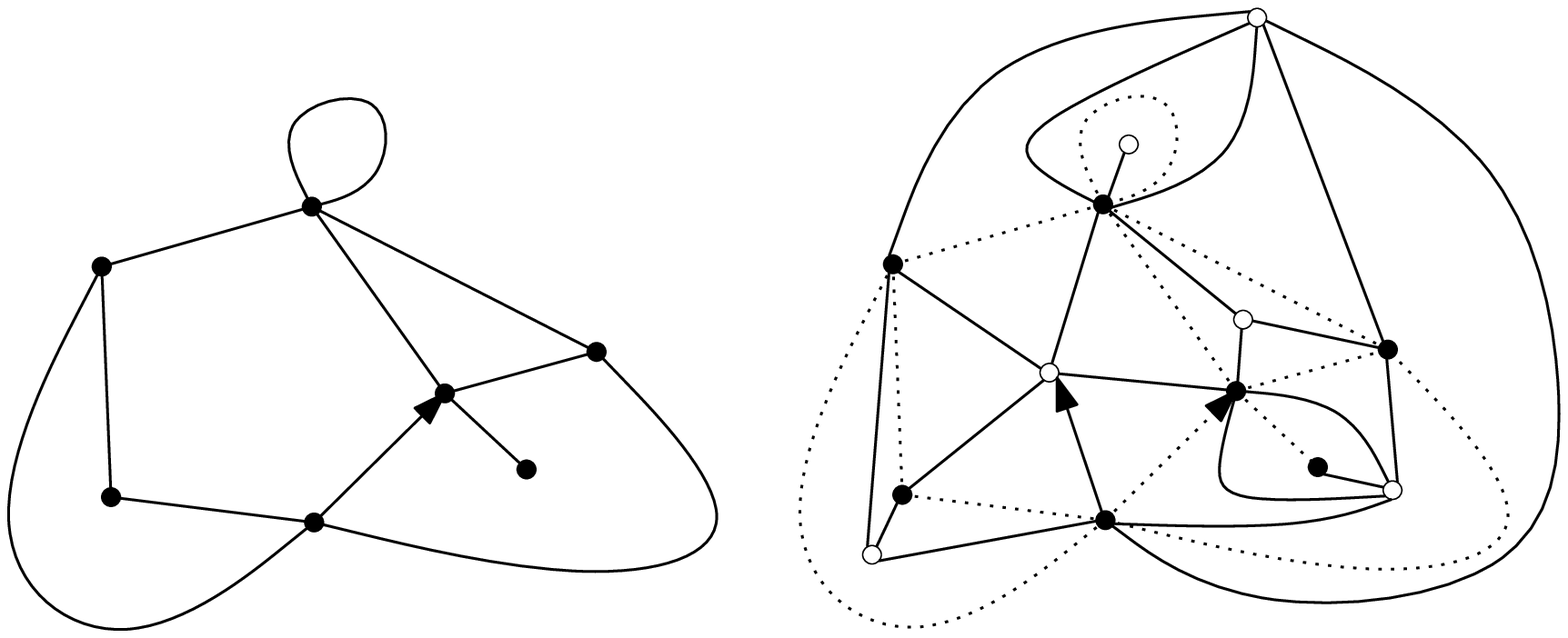}
\end{center}
\caption{The so-called ``trivial'' bijection}
\label{fig:bijtriv}
\end{figure}

For the record, we state the following useful fact. 

\begin{proposition}
A (planar) map is bipartite if and only if its faces all have even
degree.
\end{proposition}

\subsubsection{The CVS bijection}\label{sec:la-bijection-de}

Recall that $\bQ_n$ is the set of all rooted quadrangulations with $n$
faces. A simple application of Euler's formula shows that any element
of $\bQ_n$ has $2n$ edges ($4n$ oriented edges, $4$ for each face) and
$n+2$ vertices.

Let $\bT_n$ be the set of all labeled trees with $n$ edges, as defined
in Section 3.  If $(\tau,(\ell(u))_{u\in\tau})\in\bT_n$, then $\tau$
is a plane tree with $n$ edges, and $\ell:\tau\to\Z$ is a label
function on $\tau$, such that $\ell(\varnothing)=0$ and
$$|\ell(u)-\ell(\pi(u))|\leq 1\, ,\qquad \mbox{ for every }u\in
\tau\setminus\{\varnothing\}\, .$$

In order to avoid trivialities, we now assume that $n\geq 1$.  It will
be convenient here to view a plane tree $\tau$ as a planar map, by
embedding it in $\S^2$, and rooting it at the edge going from
$\varnothing$ to the vertex $1$. Let
$\varnothing=v_0,v_1,\ldots,v_{2n}=\varnothing$ be the contour
exploration of the vertices of the tree $\tau$ (see the
 end of subsection \ref{SSlabeledtree}). For $i\in
\{0,1,\ldots,2n-1\}$, we let $e_i$ be the oriented edge from $v_i$ to
$v_{i+1}$, and extend the sequences $(v_i)$ and $(e_i)$ to infinite
sequences by $2n$-periodicity. With each oriented edge $e_i$, we
can associate a corner around $e_i^-$, as explained 
in subsection 4.1. In the remaining part of Section 4, we will often identify 
the oriented edge $e_i$
with the associated corner, and we adopt the notation $\ell(e_i)=\ell(e_i^-)$. In
particular, note that $\ell(e_i)=V_i,0\leq i\leq 2n$ is the label
contour sequence as defined in Section 3.

For every $i\geq 0$, we define the {\em successor} of $i$ by
$$s(i)=\inf\{j> i:\ell(e_j)=\ell(e_i)-1\}\, ,$$ with the convention
that $\inf\emptyset=\infty$. Note that $s(i)=\infty$ if and only if
$\ell(e_i)$ equals $\min\{\ell(v):v\in \tau\}$. This is a simple
consequence of the fact that the integer-valued sequence
$(\ell(e_i),i\geq 0)$ can decrease only by taking unit steps.

Consider a point $v_*$ in $\S^2$ that does not belong to the support
of $\tau$, and denote by $e_\infty$ a corner around $v_*$, i.e.\ a
small neighborhood of $v_*$ with $v_*$ excluded, not intersecting the
corners $e_i,i\geq 0$. By convention, we set
$$\ell(v_*)=\ell(e_\infty)=\min\{\ell(u):u\in \tau\}-1.$$ 

For every $i\geq 0$, the successor
of the corner $e_i$ is then defined by
$$s(e_i)=e_{s(i)}\, .$$

\begin{figure}
\begin{center}
\includegraphics[scale=.8]{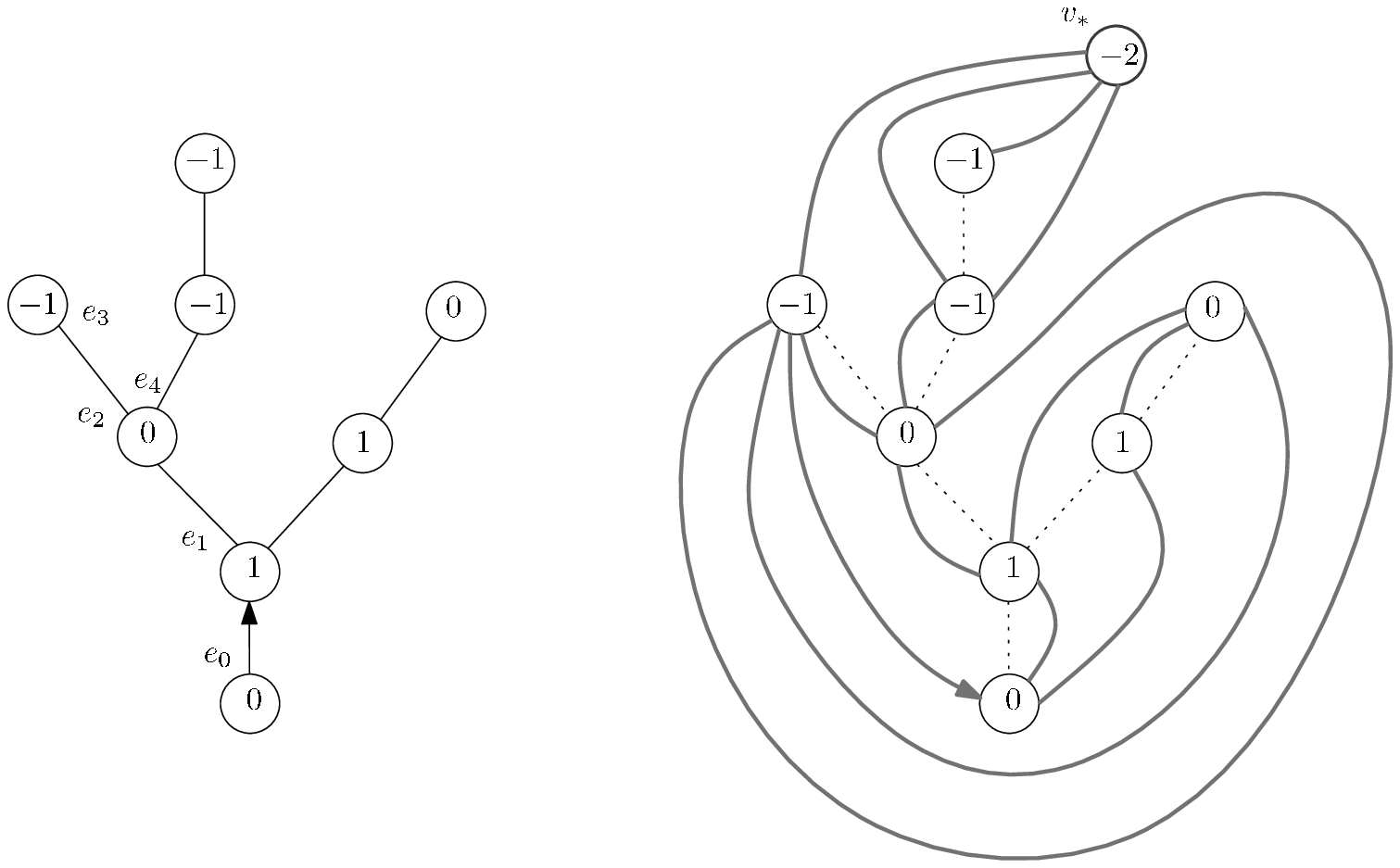}
\end{center}
\caption{Illustration of the Cori-Vauquelin-Schaeffer bijection, in
  the case $\epsilon=1$. For instance, $e_3$ is the successor of
  $e_0$, $e_2$ the successor of $e_1$, and so on.}
\label{fig:CVS}
\end{figure}

The CVS construction consists in drawing, for every $i\in
\{0,1,\ldots,2n-1\}$, an {\em arc}, which is an edge from the corner
$e_i$ to the corner $s(e_i)$ inside $\S^2\setminus(\{v_*\}\cup
\supp(\tau))$. See Fig.\ref{fig:CVS} for an illustration of the CVS construction.

\begin{lemma}\label{sec:cvs-bijection}
It is possible to draw the arcs in such a way that the graph with
vertex-set $\tau\cup \{v_*\}$ and edge-set consisting of the edges of $\tau$ and the
arcs is an embedded graph.
\end{lemma}

\proof
Since $\tau$ is a tree, we can see it as a map with a unique face
$\S^2\setminus\supp(\tau)$. The latter can in turn be seen as an open
polygon, bounded by the edges $e_0,e_1,\ldots,e_{2n-1}$ in
counterclockwise order. Hence, the result will follow if we can show
that the arcs do not cross, i.e.\ that it is not possible to find
pairwise distinct corners $e^{(1)},e^{(2)},e^{(3)},e^{(4)}$ that arise
in this order in the cyclic order induced by the contour exploration,
and such that $e^{(3)}=s(e^{(1)})$ and $e^{(4)}=s(e^{(2)})$.

If this were the case, then we would have $\ell(e^{(2)})\geq \ell(e^{(1)})$, as
otherwise the successor of $e^{(1)}$ would be between $e^{(1)}$ and
$e^{(2)}$. Similarly, $\ell(e^{(3)})\geq \ell(e^{(2)})$. But by
definition, $\ell(e^{(3)})=\ell(e^{(1)})-1$, giving $\ell(e^{(2)})\geq
\ell(e^{(3)})+1\geq \ell(e^{(2)})+1$, which is a contradiction.
\cq

\smallskip
We call $\bq$ the graph with vertex-set $V(\tau)\cup \{v_*\}$ and
edge-set formed by the arcs, now excluding the (interiors of the)
edges of $\tau$.

\begin{lemma}\label{sec:cvs-bijection-1}
The embedded graph $\bq$ is a quadrangulation with $n$ faces. 
\end{lemma}

\proof
First we check that $\bq$ is connected, and hence is a map. But this
is obvious since the consecutive successors of any given corner $e$,
given by $e,s(e),s(s(e)),\ldots$, form a finite sequence ending at
$e_\infty$. Hence, every vertex in $\bq$ can be joined by a chain to
$v_*$, and the graph is connected. 

To check that $\bq$ is a quadrangulation, let us consider an 
edge of $\tau$, corresponding to two oriented edges $e,\ov{e}$. Let
us first assume that $\ell(e^+)=\ell(e^-)-1$. Then, the successor
of $e$ is incident to $e^+$ and the preceding construction gives an arc 
starting from $e^-$ (more precisely from the corner associated with $e$)
and ending at $e^+$. Next,
let $e'$ be the corner following $\ov{e}$ in the contour exploration
around $\tau$. Then $\ell(e')=\ell(e^-)=\ell(\ov{e})+1$, giving that
$s(\ov{e})=s(s(e'))$. Indeed, $s(e')$ is the first corner coming after
$e'$ in contour order and with label $\ell(e')-1=\ell(e)-1$, while
$s(s(e'))$ is the first corner coming after $e'$ with label
$\ell(e)-2$. Therefore, it has to be the first corner coming after
$\ov{e}$, with label $\ell(e)-2=\ell(\ov{e})-1$.

We deduce that the arcs joining the corners $e$ to $s(e)$, resp. 
$\ov{e}$ to $s(\ov{e})$, resp. $e'$ to $s(e')$, resp.  $s(e')$ to $s(s(e'))=s(\ov{e})$, form a
quadrangle, that contains the edge $\{e,\ov{e}\}$, and no other edge
of $\tau$.

\begin{figure}
\begin{center}
\includegraphics[scale=.73]{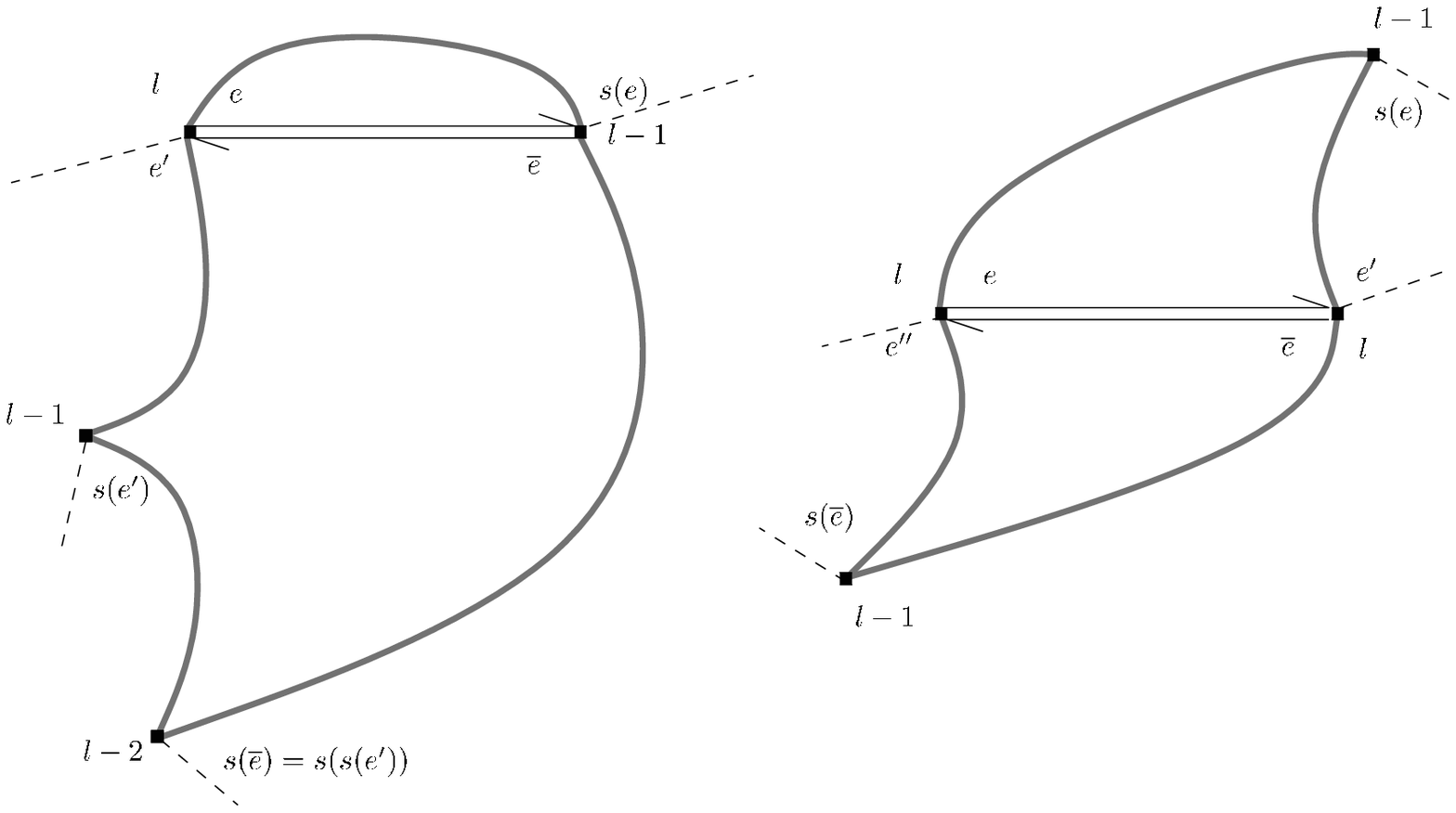}
\end{center}
\caption{Illustration of the proof of Lemma \ref{sec:cvs-bijection-1}.
In this figure, $l=\ell(e)$}
\label{fig:phipsi}
\end{figure}

If $\ell(e^+)=\ell(e^-)+1$, the situation is the same by interchanging
the roles of $e$ and $\ov{e}$. 

The only case that remains is when $\ell(e^+)=\ell(e^-)$. In this
case, if $e'$ and $e''$ are the corners following $e$ and $\ov{e}$
respectively in the contour exploration of $\tau$, then
$\ell(e)=\ell(e')=\ell(\ov{e})=\ell(e'')$, so that $s(e)=s(e')$ on the
one hand and $s(\ov{e})=s(e'')$ on the other hand.  We deduce that the
edge $\{e,\ov{e}\}$ is the diagonal of a quadrangle formed by the arcs
linking $e$ to $s(e)$, $e'$ to $s(e')=s(e)$, $\ov{e}$ to $s(\ov{e})$ and
$e''$ to $s(e'')=s(\ov{e})$.  The different cases are summed up in Fig.\ref{fig:phipsi}.

Now, notice that $\bq$ has $2n$ edges (one per corner of $\tau$) and
$n+2$ vertices, so it must have $n$ faces by Euler's formula. So all
the faces must be of the form described above. This completes the proof.
\cq

\smallskip
Note that the quadrangulation $\bq$ has a distinguished vertex $v_*$,
but for now it is not a rooted quadrangulation. To fix this root, we
will need an extra parameter $\epsilon\in \{-1,1\}$. If $\epsilon=1$
we let the root edge of $\bq$ be the arc linking $e_0$ with $s(e_0)$,
and oriented from $s(e_0)$ from $e_0$. If $\epsilon=-1$, the root edge
is this same arc, but oriented from $e_0$ to $s(e_0)$. 

In this way, we have defined a mapping $\Phi$, from
$\bT_n\times\{-1,1\}$ to the set $\bQ_n^\bullet$ of pairs $(\bq,v_*)$,
where $\bq\in \bQ_n$ and $v_*\in V(\bq)$. We call such pairs {\em
  pointed quadrangulations}. 
  
  \begin{theorem}\label{sec:la-bijection-de-3}
For every $n\geq 1$, the mapping $\Phi$ is a bijection from
$\bT_n\times\{-1,1\}$ onto $\bQ_n^\bullet$.
\end{theorem}

We omit the proof of this result. See Chassaing
and Schaeffer  \cite[Theorem 4]{CSise}. 

\begin{corollary}\label{sec:cvs-bijection-2}
We have the following formula for every $n\geq 1$: 
$$\#\bM_n=\#\bQ_n=\frac{2}{n+2}3^n\cat_n$$
\end{corollary}

\proof
We first notice that $\#\bQ_n^\bullet=(n+2)\#\bQ_n$, since every
quadrangulation $\bq\in \bQ_n$ has $n+2$ vertices, each of which
induces a distinct element of $\bQ_n^\bullet$. On the other hand, it
is obvious that 
$$\#\bT_n\times\{-1,1\}=2\cdot3^n\#\bA_n=2\cdot 3^n\cat_n\, .$$ The
result follows from Theorem \ref{sec:la-bijection-de-3}. 
\cq

\smallskip
The probabilistic counterpart of this can be stated as follows. 

\begin{corollary}\label{sec:cvs-bijection-3}
Let $Q_n$ be a uniform random element in $\bQ_n$, and conditionally
given $Q_n$, let $v_*$ be chosen uniformly at random in $V(Q_n)$. 
On the other hand, let $\theta_n$ be chosen uniformly at random in 
$\bT_n$, and let $\epsilon$ be independent of $\theta_n$ and uniformly distributed in
$\{-1,1\}$. Then $\Phi(\theta_n,\epsilon)$ has the same distribution as
$(Q_n,v_*)$.
\end{corollary}

The proof is obvious, since the probability that $(Q_n,v_*)$ equals
some particular $(\bq,v)\in \bQ_n^\bullet$ equals
$((n+2)\#\bQ_n)^{-1}=(\#\bQ_n^\bullet)^{-1}$. 

\subsubsection{Interpretation of the labels}\label{sec:interpr-labels}

The CVS bijection will be of crucial importance to us when we will
deal with metric properties of random elements of $\bQ_n$, because the
labels on $\bq$ that are inherited from a
labeled tree through the CVS construction turn out to measure certain distances in $\bq$. Recall
that the set $\tau$ is identified with $V(\bq)\setminus\{v_*\}$ if
$(\tau,\ell)$ and $\bq$ are associated through the CVS bijection (the
choice of $\epsilon$ is irrelevant here). Hence, the function $\ell$
is also a function on $V(\bq)\setminus\{v_*\}$, and we extend it by
letting, as previously, $\ell(v_*)=\min\{\ell(u): u\in \tau\}-1$. For
simplicity, we write
$$\min\ell=\min\{\ell(u):u\in \tau\}\, .$$

\begin{proposition}\label{sec:interpr-labels-1}
For every $v\in V(\bq)$, we have
\begin{equation}
\label{eq:28}
d_\bq(v,v_*)=\ell(v)-\min\ell +1\, ,
\end{equation}
where $d_\bq$ is the graph distance on $\bq$.
\end{proposition}

\proof
Let $v\in V(\bq)\setminus\{v_*\}=\tau$, and let $e$ be a corner (in
$\tau$) incident to $v$. Then the chain of arcs
$$e\to s(e)\to s^2(e)\to \ldots\to e_\infty$$ is a chain of length
$\ell(e)-\ell(e_\infty)=\ell(v)-\ell(v_*)$ between $v$ and
$v_*$. Therefore, $d_\bq(v,v_*)\leq \ell(v)-\ell(v_*)$. On the other
hand, if $v=v_0,v_1,\ldots,v_d=v_*$ are the consecutive vertices of
any chain linking $v$ to $v_*$, then since $|\ell(e)-\ell(s(e))|=1$ by
definition for any corner $e$ and since the edges of $\bq$ all connect a
corner to its successor, we get
$$d=\sum_{i=1}^d|\ell(v_i)-\ell(v_{i-1})|\geq
|\ell(v_0)-\ell(v_d)|=\ell(v)-\ell(v_*)\, ,$$
as desired. 
\cq

\smallskip
\noindent{\bf Remark.} The preceding proof also shows that the chain of arcs 
$e\to s(e)\to s^2(e)\to \ldots\to e_\infty$ is a geodesic chain linking $e^-$ to $v_*$.
Such a geodesic chain, or more generally a chain of the
form $e\to s(e)\to s^2(e)\to \ldots\to s^k(e)$, will be called a successor geodesic chain. 

\smallskip
The triangle inequality for $d_\bq$ (or the second part of the proof) gives the useful bound
\begin{equation}\label{eq:1}
d_\bq(u,v)\geq |\ell(u)-\ell(v)|\, ,
\end{equation}
This bound will be improved in the next subsection.

As a consequence of the proposition, we obtain for instance that the ``volume of
spheres'' around $v_*$ can be interpreted in terms of $\ell$:
for every $k\geq 0$,
$$|\{v\in
V(\bq):d_\bq(v,v_*)=k\}|=|\{u\in\tau:\ell(u)-\min\ell+1=k\}|\, .$$

\subsubsection{Two useful bounds}\label{sec:two-useful-bounds}

The general philosophy in the forthcoming study of random planar maps
is then the following: Information about labels in a random
labeled tree, which follows from the results of
subsection 3.3 if this tree is uniformly distributed over
$\bT_n$, allows one to obtain information about distances in
the associated quadrangulation. One major problem
with this approach is that exact information will only be
available for distances to a distinguished vertex $v_*$. There is no simple expression
for the distances between two vertices distinct from $v_*$ in terms of 
the labels in the tree. However, more advanced properties of the CVS
bijection allow to get useful bounds on these distances. Recall that
$e_0,e_1,e_2,\ldots$ is the contour sequence of corners (or oriented
edges) around a tree $\tau\in \bA_n$, starting from the root (see the beginning of
subsection \ref{sec:la-bijection-de}). We view
$(e_i,i\geq 0)$ as cyclically ordered, and for any two corners $e,e'$ of
$\tau$, we let $[e,e']$ be the set of all corners encountered when starting
from $e$, following the cyclic contour order, and stopping when visiting $e'$.

\begin{proposition}\label{sec:two-useful-bounds-1}
Let $((\tau,\ell),\epsilon)$ be an element in $\bT_n\times \{-1,1\}$,
and $(\bq,v_*)=\Phi(((\tau,\ell),\epsilon))$. Let $u,v$ be two vertices
in $V(\bq)\setminus\{v_*\}$, and let $e,e'$ be two corners of $\tau$
such that $e^-=u,(e')^-=v$. 

\noindent{\rm (i)} There holds that
$$d_\bq(u,v)\leq \ell(u)+\ell(v)-2\min_{e''\in [e,e']}\ell(e'')+2\, ,$$ 

\noindent{\rm (ii)} There holds that 
$$d_\bq(u,v)\geq \ell(u)+\ell(v)-2\min_{w\in [[u,v]]}\ell(w)\, ,$$
where $\llbracket u,v \rrbracket$ is the set of all vertices lying on the
geodesic path from $u$ to $v$ in the tree $\tau$.
\end{proposition}

\proof
For simplicity, let $m=\min_{e''\in [e,e']}\ell(e'')$. Let $e''$ be the
first corner in $[e,e']$ such that $\ell(e'')=m$. The corner
$s^k(e)$, whenever it is well defined (i.e.\ whenever
$d_\bq(e^-,v_*)\geq k$), is called the $k$-th successor of $e$.
Then $e''$ is the $(\ell(e)-m)$-th successor of
$e$. Moreover, by definition, $s(e'')$ does not belong to
$[e,e']$ since it has lesser label than $e''$, and necessarily,
$s(e'')$ is also the $(\ell(e')-m+1)$-st successor of $e'$. Hence,
the successor geodesic chain $e\to s(e)\to s^2(e)\to\cdots\to s(e'')$ from $u=e^-$ to $s(e'')^+$,
concatenated with the similar geodesic chain from $v$ to $s(e'')^+$ is a
path of length
$$\ell(u)+\ell(v)-2m+2\, ,$$ and the distance $d_\bq(u,v)$ is less
than or equal to this quantity. This proves (i).

Let us prove (ii). Let $w\in \llbracket u,v \rrbracket$ be such that 
$\ell(w)=\min\{\ell(w'):w'\in \llbracket u,v \rrbracket\}$. If 
$w=u$ or $w=v$ then the statement follows trivially from
(\ref{eq:1}). So we exclude this case. We can then write $\tau$ as
the union $\tau=\tau_1\cup\tau_2$ of two connected 
subgraphs of $\tau$ such that $\tau_1\cap\tau_2=\{w\}$,
$\tau_1$ contains $u$ but not $v$ and $\tau_2$ contains $v$
but not $u$. There may be several such decompositions, so
we just choose one.
We consider a geodesic path $\gamma$ from $u$ to $v$ in $\bq$.
If $v_*$ belongs to this path, then this means that $d_\bq(u,v)=
d_\bq(v_*,u)+d_\bq(v_*,v)$ and the desired lower bound 
immediately follows from (\ref{eq:28}). So we may assume 
that $v_*$ does not belong to $\gamma$. From our choice
of $\tau_1$ and $\tau_2$, we can then find two corners $e_{(1)}$
and $e_{(2)}$ of $\tau$ such that $e_{(1)}^-$ belongs to $\tau_1$ and $e_{(2)}^-$
belongs to $\tau_2$, $e^-_{(1)}$ and $e^-_{(2)}$ are consecutive points 
on $\gamma$, and the corners $e_{(1)}$ and $e_{(2)}$ are connected by an
edge of $\bq$. From the latter property, we must have $e_{(2)}=s(e_{(1)})$ or $e_{(1)}=s(e_{(2)})$. Consider
only the first case for definiteness (the other one is treated in
a similar fashion). Since the contour exploration of vertices of
$\tau$ must visit $w$ between any visit of $u=e_{(1)}^-$ and any
visit of $v=e_{(2)}^-$, the definition of the successor ensures that
$\ell(w)\geq \ell(e_{(2)})$ (with equality only possible if $w=e_{(2)}^-$).
Then, using (\ref{eq:1}) once again, we have
\begin{eqnarray*}
d_\bq(u,v)&=& d_\bq(u,e_{(2)}^-) + d_\bq(e_{(2)}^-,v)\\
&\geq& \ell(u)-\ell(e_{(2)}^-) + \ell(v)- \ell(e_{(2)}^-) \\
&\geq& \ell(u)+\ell(v) - 2\ell(w),
\end{eqnarray*}
giving the desired result.
\cq

\section{Basic convergence results for uniform quadrangulations}
\label{sec:limit-laws-radius}

For the remaining part of this course, our main goal will be to study
the scaling limits of random planar quadrangulations chosen according
to the uniform probability measure on $\bQ_n$. Thanks to Corollary
\ref{sec:cvs-bijection-3}, the CVS bijection and the study of scaling
limits of random labeled trees will turn out to be useful tools to study 
this problem. Ultimately,
the question we would like to address is to study the convergence in
distribution of an appropriately rescaled version of the random metric
space $(V(Q_n),d_{Q_n})$, in the sense of the Gromov-Hausdorff
topology.

One of the motivations for this problem comes from physics, and we
refer the interested reader to \cite{ADJ} for an extensive
discussion. In the past 15 years or so, physicists have been starting
to view random maps as possible discrete models for a continuum model
of random surfaces (called the Euclidean 2-dimensional quantum gravity
model), which is still ill-defined from a mathematical point of
view. We thus want to investigate whether the scaling limit of $Q_n$
exists in the above sense, and does define a certain random surface.
One can also ask the natural question of whether this limiting random
surface is {\em universal}, in the sense that it also arises as the
scaling limit of many different models of random maps, for instance,
maps chosen uniformly at random in the set of all $p$-angulations with $n$
faces: 
$$\bM^p_n=\{\bm:\deg(f)=p \mbox{ for every }f\in F(\bm),
\#F(\bm)=n\}\, , \qquad p\geq 3\, .$$ Indeed, most of the results that
we will describe in the sequel do have analogs in this more general
setting \cite{jfmgm05,mierinv,mierweill,legall06}, thanks to nice
generalizations of the CVS bijection that are due to Bouttier, Di
Francesco and Guitter \cite{BdFGmobiles}.

This is of course analogous to the celebrated Donsker Theorem,
according to which Brownian motion is the universal scaling limit of discrete
random walks, as well as to the fact that the Brownian CRT is the
scaling limit of many different models of random trees (see the 
remarks at the end of subsection \ref{convCRT}).

\subsection{Radius and profile}\label{sec:basic-conv-results}

We will first address a simpler question than the one raised above,
which is to determine by what factor we should rescale the distance
$d_{Q_n}$ in order to get an interesting scaling limit as
$n\to\infty$.

Let $\bq\in \bQ_n$ be a rooted planar quadrangulation, and $v$ be a
vertex of $\bq$.  As before, let $d_\bq$ denote the graph distance on
the vertex set of $\bq$. We define the {\em radius} of $\bq$ seen
from $v$ as
$$\mathcal{R}(\bq,v)=\max_{u\in V(\bq)} d_{\bq}(u,v)\, ,$$
and the {\em profile} of $\bq$ seen from $v$ as the sequence
$$I_{\bq,v}(k)=\card\{u\in V(\bq):d_{\bq}(u,v)=k\}\, ,\qquad k\geq 0$$
which measures the `volumes' of the spheres centered at $v$ in the
graph metric. The profile can be seen as a measure on $\Z_+$ with total
volume $n+2$. Our first limit theorem is the following.

\begin{theorem}
  \label{sec:limit-laws-radius-1}
Let $Q_n$ be uniformly distributed over $\bQ_n$,
and conditionally on $Q_n$, let $v_*$ be chosen uniformly among the
$n+2$ vertices of $Q_n$. Let also $(\ee,Z)$ be as in subsection 
\ref{convSnake}.

(i) We have
$$\left(\frac{9}{8n}\right)^{1/4}\mathcal{R}(Q_n,v_*)
\build\la_{n\to\infty}^{(d)}\sup Z-\inf Z\, .
$$

(ii) If $v_{**}$ is another vertex chosen uniformly in $V(Q_n)$ and independent of
$v_*$, 
$$\left(\frac{9}{8n}\right)^{1/4}d_{Q_n}(v_*,v_{**})
\build\la_{n\to\infty}^{(d)}\sup Z\, .$$ 

(iii) Finally, the following
convergence in distribution holds for the weak topology on probability
measures on $\R_+$:
$$\frac{I_{Q_n,v_*}((8n/9)^{1/4}\cdot)}{n+2}\build\la_{n\to\infty}^{(d)}
\mathcal{I}\, ,$$
where $\mathcal{I}$ is the occupation measure of $Z$ above its
infimum, defined as follows: For every non-negative, measurable 
$g:\R_+\to\R_+$, 
$$\langle \mathcal{I},g\rangle=\int_0^1\d t\, g(Z_t-\inf Z)\, .$$
\end{theorem}

The points (i) and (iii) are due to Chassaing and Schaeffer
\cite{CSise}, and (ii) is due to Le Gall \cite{legall05}, although
these references state these properties in a slightly different
context, namely, in the case where $v_*$ is the root vertex rather
than a uniformly chosen vertex. This indicates that as $n\to\infty$,
the root vertex plays no particular role. Some information about the
limiting distributions in (i) and (ii) can be found in Delmas \cite{Delmas}.

Property (ii) identifies the so-called $2$-point function of the Brownian map. 
An important generalization of this result has been obtained by
Bouttier and Guitter \cite{BG}, who were able to compute the $3$-point function,
namely the joint asymptotic distribution of the mutual distances between
three vertices chosen uniformly at random in $V(Q_n)$.

\smallskip
\proof
Let $((T_n,L_n),\epsilon)$ be a uniform random element in
$\bT_n\times\{-1,1\}$. Then by Corollary \ref{sec:cvs-bijection-3} we
may assume that $(Q_n,v_*)$ equals $\Phi(((T_n,L_n),\epsilon))$, where
$\Phi$ is the CVS bijection.

Let $C_n$ and $V_n$ be respectively the contour function and the label
contour function of $(T_n,L_n)$ (cf. subsections \ref{SSdiscretetree} and 
\ref{SSlabeledtree}), and 
let $u^n_i,0\leq i\leq
  2n$ be the contour exploration of vertices of $T_n$ as defined in 
  subsection \ref{SSlabeledtree}
  (so that $C_n(i)=|u^n_i|$ and $V_n(i)=L_n({u^n_i})$).
  
By Proposition \ref{sec:interpr-labels-1}, the radius of $Q_n$ viewed
from $v_*$ then equals 
$$\mathcal{R}(Q_n,v_*)=\max L_n-\min L_n+1=\max V_n-\min V_n+1\, . $$
Property (i) immediately follows from 
this equality and Theorem \ref{convsnake}.

As for (ii), we first observe that we may slightly change the hypothesis on the distribution of 
$v_{**}$. It clearly suffices to prove the desired convergence when $v_{**}$
is replaced by a vertex that is
uniformly chosen among the $n$ vertices of $Q_n$ that are distinct from both $v_*$
and the vertex $\varnothing$ of $T_n$ (recall that
$V(Q_n)\setminus\{v_*\}=V(T_n)$). 

Now, for $s\in [0,2n)$, we let $\langle s\rangle=\lceil s\rceil$
  if $C_n$ has slope $+1$ immediately after $s$, and $\langle
  s\rangle=\lf s\rf$ otherwise.  Then, if $u\in T_n$, we have
  $u^n_{\langle s\rangle}=u$ if and only if $u\not =\varnothing$ and $s$ is a time
  when the contour exploration around $T_n$ explores either of the two
  oriented edges between $u$ and its parent $\pi(u)$. Therefore, for
  every $u\in T_n\setminus\{\varnothing\}$, the Lebesgue measure of
  $\{s\in [0,2n):u^n_{\langle s\rangle}=u\}$ equals $2$.
    Consequently, if $U$ is a uniform random variable in $[0,1)$,
      independent of $(T_n,L_n)$, then $u^n_{\langle 2nU\rangle}$
      is uniform in $T_n\setminus\{\varnothing\}$. Hence, it suffices
      to prove the desired result with $u^n_{\langle 2nU\rangle}$ instead
      of $v_{**}$.

Since $|s-\langle s\rangle|\leq
1$, Theorem \ref{convsnake}
entails that
\begin{eqnarray*}
\Big(\frac{8n}{9}\Big)^{-1/4}d_{Q_n}(v_*,u^n_{\langle
  2nU\rangle})&=& \Big(\frac{8n}{9}\Big)^{-1/4} (L_n(u^n_{\langle
  2nU\rangle})-\min
L_n+1)\\
&=&\Big(\frac{8n}{9}\Big)^{-1/4}(V_n(\langle
2nU\rangle)-\min V_n+1)\, ,
\end{eqnarray*}
converges in distribution to $Z_U-\inf Z$ (here
$U$ is also assumed to be independent of $(\ee,Z)$). The fact that
$Z_U-\inf Z$ has the same distribution as $\sup Z$, or equivalently as $-\inf Z$,
can be derived from the invariance of the CRT under uniform re-rooting, see 
e.g. \cite{legweill}. This completes the proof of
(ii). 

Finally, for (iii) we just note that, for every 
bounded continuous $g:\R_+\to \R$,
\begin{eqnarray*}
\lefteqn{\frac{1}{n+2}\,\sum_{k\in\Z_+}
  I_{Q_n,v_*}(k)\,g(((8n/9)^{-1/4}k)}\\
&=&\frac{1}{n+2}\sum_{v\in Q_n} g((8n/9)^{-1/4}d_{Q_n}(v_*,v))\\
&=&E_{**}[g((8n/9)^{-1/4}d_{Q_n}(v_*,v_{**}))]\\
&\build\longrightarrow_{n\to\infty}^{} &E_U[g(Z_U-\inf Z)]\\
&=&\int_0^1\d t\, g(Z_t-\inf Z)\, ,
\end{eqnarray*}
where $E_{**}$ and $E_{U}$ means that we take the expectation only
with respect to $v_{**}$ and $U$ in the corresponding expressions
(these are conditional expectations given $(Q_n,v_*)$ and
$(\ee,Z)$ respectively). In the penultimate step, we used the convergence
established in the course of the proof of (ii). \hfill$\square$

\subsection{Convergence as metric spaces}\label{sec:conv-as-metr}

We would like to be able to understand the full scaling limit picture
for random maps, in a similar way as it was done for trees, where
we showed, using Theorem \ref{Aldous}, that the
distances in discrete trees, once rescaled by $\sqrt{2n}$, converge to
the distances in the CRT $(\TT_\ee,d_\ee)$. We thus
ask if there is an analog of the CRT that arises as the limit of the
properly rescaled metric spaces $(Q_n,d_{Q_n})$. In view of Theorem
\ref{sec:limit-laws-radius-1}, the correct normalization for the
distance should be $n^{1/4}$.

Assume that $(T_n,L_n)$ is uniformly distributed over $\T_n$, let $\epsilon$ be
uniform in $\{-1,1\}$ and independent of $(T_n,L_n)$, and let $Q_n$ be
the random uniform quadrangulation with $n$ faces and with a uniformly
chosen vertex $v_*$, which is obtained from $((T_n,L_n),\epsilon)$ via the CVS
bijection.  We now follow Le Gall \cite{legall06}\footnote{At this
  point, it should be noted that \cite{legall06,lgp,legall08} consider
  another version of Schaeffer's bijection, where no distinguished
  vertex $v_*$ has to be considered. This results in considering pairs
  $(T_n,L_n)$ in which $L_n$ is conditioned to be positive. The
  scaling limits of such pairs are still tractable, and in
  fact, are simple functionals of $(\ee,Z)$, as shown in
  \cite{legweill,legall05}. So there will be some differences from our
  exposition, but these turn out to be unimportant.}.  Recall our notation $u^n_0,u^n_1,\ldots,u^n_{2n}$
  for the contour exploration of the vertices of $T_n$, and recall that in the CVS 
  bijection these vertices are also viewed as elements of
$V(Q_n)\setminus\{v_*\}$. Define a pseudo-metric on $\{0,\ldots,2n\}$
by letting $d_n(i,j)=d_{Q_n}(u^n_i,u^n_j)$. A major problem comes from the fact
that $d_n(i,j)$ cannot be expressed as a simple functional of $(C_n,V_n)$. The only
distances that we are able to handle in an easy way are distances to
$v_*$, through the following rewriting of (\ref{eq:28}): 
\begin{equation}\label{eq:15}
d_{Q_n}(v_*,u^n_i)=V_n(i)-\min V_n+1\, .
\end{equation}
We also define, for $i,j\in \{0,1,\ldots,2n\}$, 
$$d_n^0(i,j)=V_n(i)+V_n(j)-2\max \Big(\min_{i\leq k\leq j}
V_n(k),\min_{j\leq k\leq i} V_n(k)\Big)+2\, .$$ 
Here, if $j<i$, the condition $i\leq k\leq j$  means that
$k\in\{i,i+1,\ldots,2n\}\cup\{0,1,\ldots,j\}$ and similarly 
for the condition $j\leq k\leq i$ if $i<j$.

As a consequence of
Proposition \ref{sec:two-useful-bounds-1}(i), we have the bound
$d_n\leq d_n^0$.

We now extend the function $d_n$ to $[0,2n]^2$ by letting
\begin{eqnarray}
 d_n(s,t)&=&(\lceil s\rceil-s)(\lceil t\rceil -t)d_n(\lf s \rf,\lf
 t\rf)+(\lceil s\rceil-s)(t-\lf t\rf)d_n(\lf s\rf,\lceil
 t\rceil)\nonumber
\\ &&+(s-\lf s\rf)(\lceil t\rceil
 -t)d_n(\lceil s \rceil,\lf t\rf)+(s-\lf s\rf)(t-\lf
 t\rf)d_n(\lceil s\rceil,\lceil t\rceil)\, ,\label{eq:8}
\end{eqnarray}
recalling that $\lfloor s\rfloor=\sup\{k\in \Z_+:k\leq s\}$ and $\lceil s\rceil
=\lfloor s\rfloor +1$. The function $d_n^0$ is extended to $[0,2n]^2$ by the obvious
similar formula.

It is easy to check that $d_n$ thus extended is continuous on
$[0,2n]^2$ and satisfies the triangle inequality (although this is not
the case for $d_n^0$), and that the bound $d_n\leq
d_n^0$ still holds. We define a rescaled version of these functions by letting
$$D_n(s,t)=\left(\frac{9}{8n}\right)^{1/4}d_n(2ns,2nt)\, ,\qquad 0\leq
s,t\leq 1\, .$$ 
We define similarly
the functions $D_n^0$ on $[0,1]^2$. Then, as a consequence of
Theorem \ref{convsnake}, we have
\begin{equation}
  \label{eq:5}
  (D_n^0(s,t),0\leq s,t\leq
  1)\build\la_{n\to\infty}^{(d)}(D^0(s,t),0\leq s,t\leq 1)\, ,
\end{equation}
for the uniform topology on ${C}([0,1]^2,\R)$, where by
definition
\begin{equation}\label{eq:2}
D^0(s,t)=Z_s+Z_t-2\max\Big(\min_{s\leq r\leq t}Z_r, \min_{t\leq r\leq
  s} Z_r\Big)\, ,
\end{equation} 
where if $t<s$ the condition $s\leq r\leq t$ means that $r\in[s,1]\cup[0,t]$.

We can now
state
\begin{proposition}
  \label{sec:conv-as-metr-2}
  The family of laws of $(D_n(s,t),0\leq s,t\leq 1)$, as $n$ varies,
  is relatively compact for the weak topology on probability
  measures on ${C}([0,1]^2,\R)$.
\end{proposition}

\proof
  Let $s,t,s',t'\in[0,1]$. Then by a simple use of the triangle
  inequality, and the fact that $D_n\leq D_n^0$,
$$|D_n(s,t)-D_n(s',t')|\leq D_n(s,s')+D_n(t,t')\leq
D_n^0(s,s')+D_n^0(t,t')\, ,$$ which allows one to estimate the modulus of
continuity at a fixed $\delta>0$:
\begin{equation}
  \label{eq:6}
\sup_{\substack{|s-s'|\leq \delta\\|t-t'|\leq \delta}}|D_n(s,t)-D_n(s',t')|\leq
2\sup_{|s-s'|\leq \delta}D_n^0(s,s')\, .
\end{equation}
However, the convergence in distribution (\ref{eq:5}) entails that
for every $\eps>0$, 
$$\limsup_{n\to\infty}P\left(\sup_{|s-s'|\leq \delta}D_n^0(s,s')\geq
  \eps\right) \leq P\left(\sup_{|s-s'|\leq \delta}D^0(s,s')\geq
  \eps\right)\, ,$$ and the latter quantity goes to $0$ when $\delta\to 0$
(for any fixed value of $\epsilon>0$) by the continuity of $D^0$ and the fact that
$D^0(s,s)=0$. Hence, taking $\eta>0$ and letting $\eps=\eps_k=2^{-k}$,
we can choose $\delta=\delta_k$ (tacitly depending also on $\eta$)
such that
$$\sup_{n\geq 1}P\left(\sup_{|s-s'|\leq \delta_k}D_n^0(s,s')\geq
  2^{-k}\right)\leq \eta 2^{-k}\, ,\qquad k\geq 1,$$ entailing
$$P\left(\bigcap_{k\geq 1}\left\{\sup_{|s-s'|\leq
    \delta_k}D_n^0(s,s')\leq 2^{-k}\right\} \right)\geq 1-\eta\, ,$$
  for all $n\geq 1$.  Together with (\ref{eq:6}), this shows that with
  probability at least $1-\eta$, the function $D_n$ belongs to the set of all
  functions $f$ from $[0,1]^2$ into $\R$ such that $f(0,0)=0$ and, for every
  $k\geq 1$,
$$\sup_{\substack{|s-s'|\leq \delta_k\\|t-t'|\leq
    \delta_k}}|f(s,t)-f(s',t')|\leq 2^{-k}\, .$$ The latter set is
compact by the Arzel\`a-Ascoli theorem. The conclusion then follows from
Prokhorov's theorem.
\hfill$\square$

\smallskip
At this point, we are allowed to say that the random distance
functions $D_n$ admit a limit in distribution, up to taking
$n\to\infty$ along a subsequence:
\begin{equation}\label{eq:10}
  (D_n(s,t),0\leq s,t\leq 1)\build\la_{}^{(d)}(
  D(s,t),0\leq s,t\leq 1)
\end{equation}
for the uniform topology on ${C}([0,1]^2,\R)$. In fact, we are going
to need a little more than the convergence of $D_n$. From the relative
compactness of the components, we see that the closure of the collection of laws of the
triplets $$((2n)^{-1}C_n(2n\cdot),(9/8n)^{1/4}V_n(2n\cdot),D_n),\quad n\geq 1$$
is compact in the space of all probability measures on
${C}([0,1],\R)^2\times {C}([0,1]^2,\R)$. Therefore, it is possible to			%CMI
choose a subsequence $(n_k,k\geq 1)$ so that this triplet converges in
distribution to a limit, which is denoted by $(\ee,Z,D)$ (from Theorem
\ref{convsnake}, this is of course consistent with the preceding
notation). The joint convergence to the triplet $(\ee,Z,D)$ gives a
coupling of $D,D^0$ such that $D\leq D^0$, since $D_n\leq D_n^0$ for
every $n$.

Define a random equivalence relation on $[0,1]$ by letting $s\approx
t$ if $D(s,t)=0$.  We let $M=[0,1]/\approx$ be the associated quotient space,
endowed with the quotient distance, which we still denote by $D$. The
canonical projection $[0,1]\to M$ is denoted by $\bp$. 

Finally, let $s_*\in [0,1]$ be such that $Z_{s_*}=\inf Z$. It can be proved
that $s_*$ is unique a.s., see \cite{MM05} or \cite{legweill}, and we will admit this fact
(although it is not really needed for the next statement). We set
$\rho_*=\bp(s_*)$.  We can now state the main result of this section.

\begin{theorem}
  \label{sec:conv-as-metr-3}
  The random pointed metric space $(M,D,\rho_*)$ is the limit in distribution of the
  spaces $(V(Q_n),(9/8n)^{1/4}d_{Q_n},v_*)$, for the
  Gromov-Hausdorff topology, along the subsequence $(n_k,k\geq
  1)$. Moreover, we have a.s.\, for every $x\in M$ and $s\in
  [0,1]$ such that $\bp(s)=x$,
$$D(\rho_*,x)=D(s_*,s)=Z_s-\inf Z\, .$$
\end{theorem}

Note that, in the discrete model, a point at which the minimal label
in $T_n$ is attained lies at distance $1$ from $v_*$. Therefore, the
point $\rho_*$ should be seen as the continuous analog of the
distinguished vertex $v_*$.  The last identity in the statement of the theorem is
then of course the continuous analog of (\ref{eq:28}) and
(\ref{eq:15}). 

\smallskip
\proof
For the purposes of this proof, it is useful to assume, using the
Skorokhod representation theorem, that the convergence 
$$((2n)^{-1/2}C_n(2n\cdot),(9/8n)^{1/4}V_n(2n\cdot),D_n)\la (\ee,Z,D)$$
holds a.s.\ along the subsequence $(n_k)$. In what follows we restrict our attention
to values of $n$ in this sequence.

For every $n$, let $i^{(n)}_*$ be
any index in $\{0,1,\ldots,2n\}$ such that $V_n(i^{(n)}_*)=\min
V_n$. Then for every $v\in V(Q_n)$, it holds that 
$$|d_{Q_n}(v_*,v)-d_{Q_n}(u^n_{i^{(n)}_*},v)|\leq 1$$ because
$d_{Q_n}(v_*,u^n_{i^{(n)}_*})=1$ ($v_*$ and $u^n_{i^{(n)}_*}$ are linked by an arc
in the CVS bijection). Moreover, since $(8n/9)^{-1/4}V_n(2n\cdot)$
converges to $Z$ uniformly on $[0,1]$, and since we know\footnote{We
  could also perform the proof without using this fact, but it makes
  things a little easier} that $Z$ attains its overall infimum at a
unique point $s_*$, it is easy to obtain that $i^{(n)}_*/2n$ converges
as $n\to\infty$ towards $s_*$. 

For every integer $n$, we construct a correspondence $\RR_n$ between
$V(Q_n)$ and $M_n$, by putting:
\begin{enumerate}
\item[$\bullet$] $(v_*,\rho_*)\in \RR_n$\,;
\item[$\bullet$] $(u^n_{\lfloor 2ns\rfloor},\bp(s))\in\RR_n$, for every $s\in[0,1]$.
\end{enumerate}

We then verify that the distortion of $\RR_n$ (with respect to the metrics 
$(9/8n)^{1/4}d_{Q_n}$ on $V(Q_n)$ and $D$ on $M$) converges to
$0$ 	a.s. as $n\to\infty$. We first observe that
\begin{align*}
&\sup_{s\in[0,1]} | (9/8n)^{1/4}d_{Q_n}(v_*,u^n_{\lfloor 2ns\rfloor})
- D(\rho_*,\bp(s))|\\
&\quad\leq 
(9/8n)^{1/4}+ \sup_{s\in[0,1]} | (9/8n)^{1/4}d_{Q_n}(u^n_{i_*^{(n)}},u^n_{\lfloor 2ns\rfloor})
- D(\rho_*,\bp(s))|\\
&\quad=(9/8n)^{1/4}+ \sup_{s\in[0,1]} |D_n(i_*^{(n)}/2n, \lfloor 2ns\rfloor/2n)
- D(s_*,s)|,
\end{align*}
which tends to $0$ as $n\to\infty$, by the a.s. uniform convergence 
of $D_n$ to $D$, and the fact that $i^{(n)}_*/2n$ converges
to $s_*$. Similarly, we have
\begin{align*}
&\sup_{s,t\in[0,1]} | (9/8n)^{1/4}d_{Q_n}(u^n_{\lfloor 2ns\rfloor},u^n_{\lfloor 2nt\rfloor})
- D(\bp(s),\bp(t))|\\
&\quad = \sup_{s,t\in[0,1]} |D_n(\lfloor 2ns\rfloor/2n,\lfloor
2nt\rfloor/2n)-D(s,t)|
\end{align*}
which tends to $0$ as $n\to\infty$. We conclude that the distortion of $\RR_n$
converges to $0$ a.s. and that the pointed metric spaces 
$(V(Q_n),(9/8n)^{-1/4}d_{Q_n},v_*)$ also converge a.s. 
to $(M,D,\rho_*)$ in the Gromov-Hausdorff topology. 

Let us prove the last statement of the theorem. 
Using once again the uniform
convergence of $D_n$ to $D$, we obtain that for every
$s\in [0,1]$,
\begin{eqnarray*}
D(s_*,s)&=&\lim_{n\to\infty} D_n(i_*^{(n)}/2n,\lfloor 2ns \rfloor/2n)
\\
&=&\lim_{n\to\infty}\left(\frac{8n}{9}\right)^{-1/4}d_{Q_n}(v_*,u^n_{\lfloor
  2ns\rfloor })\\
&=&\lim_{n\to\infty}\left(\frac{8n}{9}\right)^{-1/4}(V_n(\lfloor
  2ns\rfloor )-\min V_n+1)\\
&=&Z_s-\inf Z\, ,
\end{eqnarray*}
as desired. 
\cq

\medskip

It is tempting to call $(M,D)$ the ``Brownian map'', or the ``Brownian
continuum map'', by analogy with the fact that the ``Brownian
continuum random tree'' is the scaling limit of  uniformly distributed
plane trees with $n$ edges. However, the choice of the subsequence in
Theorem \ref{sec:conv-as-metr-3} poses a problem of uniqueness of
the limit. As we see in the previous statement, only the distances to
$\rho_*$ are {\em a priori} defined as simple functionals of the process
$Z$. Distances between other points in $M$ seem to be harder to
handle.
The following conjecture is however very appealing.

\begin{conjecture}
 \label{sec:conv-as-metr-4}
 The spaces $(V(Q_n),n^{-1/4}d_{Q_n})$ converge in distribution, for
 the Gromov-Hausdorff topology.
\end{conjecture}

Marckert and Mokkadem \cite{MM05} and Le Gall \cite{legall06} give a
natural candidate for the limit (called the Brownian map in
\cite{MM05}) but until now the convergence result in the above
conjecture has not been proved.

\section{Identifying the Brownian map}\label{sec:ident-brown-map}

\subsection{The Brownian map as a quotient of the CRT}\label{sec:brownian-map-as}

In the previous section, we wrote the scaling limit of rescaled 
random quadrangulations (along a suitable subsequence)  as a quotient space $M=[0,1]/\!\approx$
where the equivalence relation $\approx$
is defined by $s\approx t$ iff $D(s,t)=0$. In this section, we
provide a more explicit description of this quotient.

Recall the notation of the previous section. In particular,
$((T_n,L_n),\epsilon)$ is uniformly distributed over 
$\mathbf{T}_n\times\{-1,1\}$, and $(Q_n,v_*)$ is the pointed quadrangulation
that is the image of $((T_n,L_n),\epsilon)$ under the 
CVS bijection. For every $n\geq 1$, $u^n_0,u^n_1,\ldots, u^n_{2n}$ 
is the contour exploration of the vertices of $T_n$. Thus, $C_n(k)= |u^n_k|$
and $V_n(k)=L_n(u^n_k)$ for $0\leq k\leq 2n$. 

As in the proof 
of Theorem \ref{sec:conv-as-metr-3}, we may assume that,
along the sequence $(n_k)$ we have the almost sure
convergence
\begin{eqnarray}
\label{Groconv}
((2n)^{-1/2}C_n(2ns),(9/8n)^{1/4}V_n(2ns),D_n(s,t))_{s,t\in[0,1]}\\
\build{\la}_{n\to\infty}^{} (\ee_s,Z_s,D(s,t))_{s,t\in[0,1]}\nonumber
\end{eqnarray}
uniformly over $[0,1]^2$. 
Recall from the proof of Theorem \ref{sec:conv-as-metr-3}
that this implies the almost sure convergence
$$
\Big(V(Q_n),\Big(\frac{9}{8}\Big)^{1/4} d_{Q_n}\Big)
\build{\la}_{n\to\infty}^{} (M,D)
$$
in the Gromov-Hausdorff sense, along the sequence $(n_k)$.

As in Section 2 above, introduce the random equivalence 
relation $\sim_\ee$ on $[0,1]$ by
$$s\sim_\ee t \hbox{ iff }\ee_s=\ee_t=\min_{s\wedge t\leq r\leq s\vee t} \ee_r$$
and recall that the CRT $\t_\ee$ is defined as the quotient space
$[0,1]/\!\sim_\ee$ equipped with the distance $d_\ee$. 

\begin{lemma}
\label{quotientCRT}
We have almost surely for every $s,t\in[0,1]$,
$$s\sim_\ee t \Rightarrow D(s,t) = 0 \quad (\Leftrightarrow s\approx t).$$
\end{lemma}

\proof We can use the convergence of the first components in (\ref{Groconv})
to see that if $s\sim_\ee t$ and $s<t$ we can find integers $i_n<j_n$
such that $i_n/2n\la s$, $j_n/2n\la t$, and, for every 
sufficiently large $n$ (belonging to the sequence $(n_k)$),
$$C_n(i_n)=C_n(j_n) = \min_{i_n\leq k\leq j_n} C_n(k).$$
Then, from the definition of the contour function, we must have
$u^n_{i_n}=u^n_{j_n}$, and thus $d_n(i_n,j_n)=0$. Using the convergence  
(\ref{Groconv}) again, we conclude that $D(s,t)=0$.
\endproof

\medskip
\noindent{\bf Consequence.} Recall that $p_\ee: [0,1] \la \t_\ee$
denotes the canonical projection. Then $D(s,t)$ only depends 
on $p_\ee(s)$ and $p_\ee(t)$. We can therefore put
for every $a,b\in\t_\ee$,
$$D(a,b)=D(s,t)$$
where $s$, resp. $t$, is an arbitrary representative of $a$,
resp. of $b$, in $[0,1]$. Then $D$ is (again) a pseudo-distance
on $\t_\ee$. With a slight abuse of notation we keep 
writing $a\approx b$ iff $D(a,b)=0$, for $a,b\in \t_\ee$. 
Then the Brownian map $M$ can be written as
$$M= [0,1]/\! \approx \;= \t_\ee /\! \approx$$
where the first equality was the definition of $M$ and the second 
one corresponds to the fact that there is an obvious canonical isometry
between the two quotient spaces.

One may wonder why it is more interesting to write the Brownian 
map $M$ as a quotient space of the CRT $\t_\ee$ rather than as
a quotient space of $[0,1]$. The point is that it will be possible
to give a simple intuitive description of $\approx$ viewed 
as an equivalence relation on $\t_\ee$. This is indeed the main
goal of the next section.

\subsection{Identifying the equivalence relation $\approx$}

We noticed in subsection \ref{snakedeter} that the process $Z$ (the Brownian snake
driven by $\ee$) can be viewed as indexed by $\t_\ee$. This will
be important in what follows: For $a\in\t_\ee$, we will write 
$Z_a=Z_t$ for any choice of $t$ such that $a=p_\ee(t)$. 
We also set $a_*=p_\ee(s_*)$: $a_*$ is thus the unique vertex of $\t_\ee$
such that
$$Z_{a_*}=\min_{a\in\t_\ee} Z_a.$$

We first need to define intervals on the tree $\t_\ee$. For simplicity
we consider only leaves of $\t_\ee$. Recall that a point $a$ of $\t_\ee$
is a leaf if $\t_\ee\backslash \{a\}$ is connected. 
Equivalently a vertex $a$ distinct from the root $\rho$ is a leaf if and only $p_\ee^{-1}(a)$ is a singleton. Note in
particular that $a_*$ is a leaf of $\t_\ee$. 

Let $a$ and $b$ be two (distinct) leaves of $\t_\ee$, and let $s$ and $t$ be the 
unique elements of $[0,1)$ such that $p_\ee(s)=a$ and $p_\ee(t)=b$. 
Assume that $s<t$ for definiteness. We then set
\begin{align*}
  &[a,b]=p_\ee([s,t])\\
  &[b,a]=p_\ee([t,1]\cup [0,s]).
\end{align*}
It is easy to verify that $[a,b]\cap[b,a]= \llbracket a,b\rrbracket$
is the line segment between $a$ and $b$ in $\t_\ee$. 

\begin{theorem}
\label{equival-relation}
Almost surely, for every distinct $a,b\in\t_\ee$,
$$a\approx b\quad \Leftrightarrow
\left\{\begin{array}{l}
a,b\hbox{ are leaves of }\t_\ee\hbox{ and}\\
\noalign{\smallskip}
Z_a=Z_b=\max\Big(\min_{c\in[a,b]} Z_c, \min_{c\in[b,a]} Z_c\Big)
\end{array}
\right.
$$
\end{theorem}

\begin{remark} We know that the minimum of $Z$ over $\t_\ee$ is attained at the
unique vertex $a_*$. If $a$ and $b$ are (distinct) leaves of $\t_\ee\backslash\{a_*\}$,
exactly one of the two intervals $[a,b]$ and $[b,a]$ contains the vertex $a_*$. Obviously
the minimum of $Z$ over this interval is equal to $Z_{a_*}$ and thus cannot be equal 
to $Z_a$ or $Z_b$.
\end{remark}

The proof of the implication $\Leftarrow$ in the theorem is easy. Suppose that $a=p_\ee(s)$
and $b=p_\ee(t)$ with $s<t$ (for definiteness). If 
$$Z_a=Z_b=\max\Big(\min_{c\in[a,b]} Z_c, \min_{c\in[b,a]} Z_c\Big)$$
this means that
$$Z_s=Z_t=\max\Big(\min_{r\in[s,t]} Z_r, \min_{r\in[t,1]\cup[0,s]} Z_r\Big).$$
The last identity is equivalent to saying that $D^0(s,t)=0$, and since $D\leq D^0$
we have also $D(s,t)=0$, or equivalently $a\approx b$.

Unfortunately, the proof of the converse implication is much harder, and we will
only give some key ideas of the proof, referring to \cite{legall06} for additional details.

We start with a preliminary lemma. We denote by ${\rm vol}(\cdot)$ the mass measure
on $\t_\ee$, which is simply the image of the Lebesgue measure on $[0,1]$
under the projection $p_\ee:[0,1]\la \t_\ee$.

\begin{lemma}
\label{volume-ball}
Almost surely, for every $\delta\in(0,1)$, there exists a (random)
constant $C_\delta(\omega)$ such that, for every $r>0$ and every
$a\in\t_\ee$,
$${\rm vol}(\{b\in\t_\ee: D(a,b)\leq r\}) \leq C_\delta\,r^{4-\delta}.$$
\end{lemma}

We omit the proof of this lemma. The first ingredient of the proof is a
``re-rooting invariance property'' of random planar maps, which makes it possible
to reduce the proof to the case $a=a_*$. In that case we can use
the formula $D(a_*,b)=Z_b - \min Z$ and explicit moment
calculations for the Brownian snake (see Corollary 6.2 in
\cite{legall08} for a detailed proof).

Let us come to the proof of the implication $\Rightarrow$ in Theorem
\ref{equival-relation}.  For simplicity we consider only the case when
$a$ and $b$ are leaves of $\t_\ee$ (it would be necessary to show also
that the equivalence class of any vertex of $\t_\ee$ that is not a
leaf is a singleton -- this essentially follows from Lemma 2.2 in \cite{legall06}). We let $s,t\in[0,1]$ be such that $a=p_\ee(s)$
and $b=p_\ee(t)$, and assume for definiteness that $0\leq s_*<s<t\leq
1$.

We assume that $a\approx b$, and our goal is to prove that
$$Z_a=Z_b=\min_{c\in[a,b]} Z_c.$$
We already know that $Z_a=Z_b$, because 
$$Z_a- \min Z = D(a_*,a)=D(a_*,b)= Z_b - \min Z.$$

\noindent{\bf First step.} We first establish that
\begin{equation}
\label{claim1}
Z_a=Z_b=\min_{c\in\llbracket a,b\rrbracket} Z_c.
\end{equation}
To see this, we go back to the discrete picture. We can find $a_n,b_n\in T_n$ such 
that $a_n \la a$ and $b_n \la b$ as $n\to \infty$ (strictly speaking these convergences 
make no sense: What we mean is that $a_n=u^n_{i_n}$, $b_n=u^n_{j_n}$ with
${i_n}/{2n} \la s$ and ${j_n}/{2n} \la t$). Then the condition $D(a,b)=0$
implies that
\begin{equation}
\label{conv-dista}
n^{-1/4}\, d_{Q_n}(a_n,b_n) \la 0 .
\end{equation}

Recall, from Proposition \ref{sec:two-useful-bounds-1},
the notation  $\llbracket a_n,b_n\rrbracket$ for the set of vertices lying
on the geodesic path from $a_n$ to $b_n$ in the tree $T_n$.
By Proposition \ref{sec:two-useful-bounds-1}(ii), we have
$$d_{Q_n}(a_n,b_n)\geq L_n(a_n)+L_n(b_n)- 2\,\min_{c\in\llbracket a_n,b_n\rrbracket} L_n(c).$$
We multiply both sides of this inequality by $n^{-1/4}$ and let $n$
tend to $\infty$, using (\ref{conv-dista}). Modulo some technical details that we omit
(essentially one needs to check that any vertex of $\t_\ee$ belonging to
$\llbracket a,b\rrbracket$ is of the form $p_\ee(r)$, where $r=\lim k_n/2n$ and the integers
$k_n$ are such that $u^n_{k_n}$ belongs to $\llbracket a_n,b_n\rrbracket$),
we get that 
$$Z_a+Z_b-2\min_{c\in\llbracket a,b \rrbracket} Z_c\leq 0$$
from which (\ref{claim1}) immediately follows.
\smallskip

\noindent{\bf Second step.} 
We argue by contradiction, assuming that
$$\min_{c\in[a,b]} Z_c < Z_a=Z_b.$$

Let $\gamma_n$ be a discrete geodesic from $a_n$ to $b_n$ in the quadrangulation $Q_n$
(here we view $a_n$ and $b_n$ as vertices of the quadrangulation $Q_n$,
and this geodesic is of course different from the geodesic from $a_n$ to $b_n$
in the tree $T_n$). From (\ref{conv-dista}) the maximal distance between $a_n$
(or $b_n$) and a vertex visited by $\gamma_n$ is $o(n^{1/4})$ as $n\to\infty$. 
As a consequence, using the triangle inequality and (\ref{eq:28}), we have
$$\sup_{u\in\gamma_n} |L_n(u) - L_n(a_n)|=o(n^{1/4})$$
as $n\to\infty$.

To simplify the presentation of the argument, we assume
that, for infinitely many values of $n$, the geodesic path $\gamma_n$
from $a_n$ to $b_n$ stays in the lexicographical interval $[a_n,b_n]$.
This lexicographical interval is defined, analogously to the 
continuous setting, as the set of all vertices visited by the contour exploration
sequence $(u^n_i)_{0\leq i\leq 2n}$ between its last visit of $a_n$
and its first visit of $b_n$. Note that the preceding assumption may not
hold, and so the real argument is slightly more complicated
than what follows.

We use the previous assumption to prove the following claim.
If $x\in [a,b]$, we denote by $\phi_{a,b}(x)$ the last ancestor of 
$x$ that belongs to $\llbracket a,b\rrbracket$ (the condition 
$x\in [a,b]$ ensures that the ancestral line $\llbracket\rho,x\rrbracket$
intersects $\llbracket a,b\rrbracket$). Alternatively, $\phi_{a,b}(x)$ is
the point of $\llbracket a,b\rrbracket$ at minimal $d_\ee$-distance of $x$
in the tree $\t_\ee$. 

\medskip
\noindent{\bf Claim.} Let $\eps>0$. For every $c\in[a,b]$ such that
$$\left\{ \begin{array}{ll}
Z_c<Z_a +\eps&\\
Z_x>Z_a +\eps/2&\qquad\forall x\in \llbracket \phi_{a,b}(c),c\rrbracket
\end{array}\right.$$
we have $D(a,c)\leq \eps$. 

\medskip
The claim eventually leads to the desired contradiction: Using the first step of the proof (which ensures that 
$Z_c\geq Z_a$ for $c\in\llbracket a,b\rrbracket$) and the properties
of the Brownian snake, one can check that, under the condition
$$\min_{c\in[a,b]} Z_c < Z_a=Z_b,$$
 the volume of the
set of all vertices $c$ that satisfy the assumptions of the claim is bounded 
below by a (random) positive constant times $\eps^2$, at least
for sufficiently small $\eps>0$ (see Lemma 2.4 in \cite{legall06}
for a closely related statement). The desired contradiction follows
since Lemma \ref{volume-ball} implies that, for every $\delta\in(0,1)$, 
$${\rm vol}(\{c: D(a,c)\leq \eps\})\leq C_\delta\,\eps^{4-\delta}.$$

\begin{figure}
\begin{center}
\ifx\JPicScale\undefined\def\JPicScale{1}\fi
\unitlength \JPicScale mm
\begin{picture}(114,92)(0,0)
\linethickness{0.25mm}
\multiput(20,0)(0.12,0.12){750}{\line(1,0){0.12}}
\linethickness{0.25mm}
\multiput(10,30)(0.12,-0.12){167}{\line(1,0){0.12}}
\linethickness{0.25mm}
\multiput(20,20)(0.12,0.12){83}{\line(1,0){0.12}}
\linethickness{0.25mm}
\put(20,20){\line(0,1){10}}
\linethickness{0.25mm}
\put(40,20){\line(0,1){20}}
\linethickness{0.25mm}
\linethickness{0.25mm}
\multiput(20,70)(0.12,-0.12){292}{\line(1,0){0.12}}
\linethickness{0.25mm}
\multiput(40,50)(0.12,0.12){83}{\line(1,0){0.12}}
\linethickness{0.25mm}
\multiput(45,65)(0.12,-0.12){42}{\line(1,0){0.12}}
\linethickness{0.25mm}
\multiput(50,60)(0.12,0.12){42}{\line(1,0){0.12}}
\linethickness{0.25mm}
\multiput(25,15)(0.12,0.12){42}{\line(1,0){0.12}}
\linethickness{0.25mm}
\linethickness{0.25mm}
\put(30,60){\line(0,1){10}}
\linethickness{0.25mm}
\put(35,55){\line(0,1){5}}
\linethickness{0.25mm}
\multiput(65,65)(0.12,-0.12){83}{\line(1,0){0.12}}
\linethickness{0.25mm}
\multiput(60,70)(0.12,-0.12){42}{\line(1,0){0.12}}
\linethickness{0.25mm}
\multiput(55,85)(0.12,-0.12){83}{\line(1,0){0.12}}
\linethickness{0.25mm}
\put(65,75){\line(0,1){5}}
\linethickness{0.25mm}
\multiput(65,75)(0.12,0.12){83}{\line(1,0){0.12}}
\linethickness{0.25mm}
\put(70,80){\line(0,1){5}}
\linethickness{0.25mm}
\put(60,80){\line(0,1){10}}
\linethickness{0.25mm}
\multiput(95,85)(0.12,-0.12){42}{\line(1,0){0.12}}
\linethickness{0.25mm}
\put(85,65){\line(1,0){15}}
\linethickness{0.25mm}
\multiput(95,65)(0.12,0.12){83}{\line(1,0){0.12}}
\linethickness{0.25mm}
\put(100,70){\line(1,0){5}}
\linethickness{0.25mm}
\multiput(80,45)(0.12,0.12){42}{\line(1,0){0.12}}
\linethickness{0.25mm}
\multiput(80,45)(0.12,-0.12){42}{\line(1,0){0.12}}
\linethickness{0.25mm}
\put(45,25){\line(1,0){10}}
\linethickness{0.25mm}
\put(35,15){\line(1,0){5}}
\linethickness{0.25mm}
\put(30,50){\line(1,0){10}}
\linethickness{0.25mm}
\put(10,25){\line(1,0){5}}
\linethickness{0.25mm}
\put(15,15){\line(1,0){10}}
\linethickness{0.25mm}
\multiput(70,60)(0.12,0.12){125}{\line(1,0){0.12}}
\linethickness{0.25mm}
\put(80,70){\line(0,1){10}}
\linethickness{0.25mm}
\put(60,40){\line(0,1){10}}
\linethickness{0.75mm}
\qbezier(20,70)(23.91,73.28)(26.88,71.88)
\qbezier(26.88,71.88)(29.84,70.47)(30,70)
\linethickness{0.75mm}
\qbezier(30,70)(33.59,68.06)(34.38,64.25)
\qbezier(34.38,64.25)(35.16,60.44)(35,60)
\linethickness{0.25mm}
\put(65,65){\line(0,1){10}}
\linethickness{0.75mm}
\qbezier(60,70)(56.72,73.91)(58.12,76.88)
\qbezier(58.12,76.88)(59.53,79.84)(60,80)
\linethickness{0.75mm}
\qbezier(60,80)(61.94,83.59)(65.75,84.38)
\qbezier(65.75,84.38)(69.56,85.16)(70,85)
\linethickness{0.75mm}
\qbezier(70,85)(79.81,89.28)(80.25,85.12)
\qbezier(80.25,85.12)(80.69,80.97)(80,80)
\linethickness{0.75mm}
\qbezier(95,85)(99.22,89.91)(104.38,90.12)
\qbezier(104.38,90.12)(109.53,90.34)(110,90)
\linethickness{0.25mm}
\put(100,80){\line(0,1){5}}
\linethickness{0.25mm}
\put(65,45){\line(1,0){15}}
\linethickness{0.25mm}
\multiput(70,45)(0.12,-0.12){42}{\line(1,0){0.12}}
\linethickness{0.25mm}
\multiput(65,55)(0.12,-0.12){42}{\line(1,0){0.12}}
\linethickness{0.25mm}
\put(65,60){\line(1,0){5}}
\linethickness{0.5mm}
\qbezier(60,50)(63.59,49.69)(64.38,47.5)
\qbezier(64.38,47.5)(65.16,45.31)(65,45)
\linethickness{0.5mm}
\qbezier(65,45)(66.31,48.91)(65.75,51.88)
\qbezier(65.75,51.88)(65.19,54.84)(65,55)
\linethickness{0.5mm}
\qbezier(65,55)(61.72,58.59)(63.12,59.38)
\qbezier(63.12,59.38)(64.53,60.16)(65,60)
\linethickness{0.5mm}
\qbezier(65,60)(60.75,62.59)(60.25,66.12)
\qbezier(60.25,66.12)(59.75,69.66)(60,70)
\linethickness{0.25mm}
\multiput(75,75)(0.12,-0.12){42}{\line(1,0){0.12}}
\linethickness{0.75mm}
\qbezier(80,80)(76.41,79.69)(75.62,77.5)
\qbezier(75.62,77.5)(74.84,75.31)(75,75)
\linethickness{0.75mm}
\qbezier(75,75)(75.31,72.06)(77.5,71)
\qbezier(77.5,71)(79.69,69.94)(80,70)
\linethickness{0.75mm}
\qbezier(80,70)(83.56,79.47)(89,82.38)
\qbezier(89,82.38)(94.44,85.28)(95,85)
\linethickness{0.75mm}
\qbezier(35,60)(38.91,64.91)(41.88,65.12)
\qbezier(41.88,65.12)(44.84,65.34)(45,65)
\linethickness{0.75mm}
\qbezier(45,65)(42.38,61.09)(43.5,58.12)
\qbezier(43.5,58.12)(44.62,55.16)(45,55)
\linethickness{0.75mm}
\qbezier(45,55)(52.5,55.94)(56.25,62.5)
\qbezier(56.25,62.5)(60,69.06)(60,70)
\linethickness{0.25mm}
\multiput(35,35)(0.12,-0.12){42}{\line(1,0){0.12}}
\put(18,3){\makebox(0,0)[cc]{$\varnothing$}}

\put(48,10){\makebox(0,0)[cc]{tree $T_n$}}

\put(17,72){\makebox(0,0)[cc]{$a_n$}}

\put(114,90){\makebox(0,0)[cc]{$b_n$}}

\linethickness{0.25mm}
\put(105,85){\line(1,0){5}}
\linethickness{0.25mm}
\put(100,80){\line(1,0){5}}
\put(58,50){\makebox(0,0)[cc]{$u$}}

\put(57,70){\makebox(0,0)[cc]{$w$}}

\put(72,89){\makebox(0,0)[cc]{$\gamma_n$}}

\end{picture}

\caption{Illustration of the proof: The geodesic path $\gamma_n$
from $a_n$ to $b_n$ is represented by the thick curves.
The thin curves correspond to the beginning of the successor geodesic chain starting from $u$. 
This chain does not cross the line segment $\llbracket a_n,b_n \rrbracket$
and thus has to meet the path $\gamma_n$ at some point $w$.}
\end{center}
\end{figure}
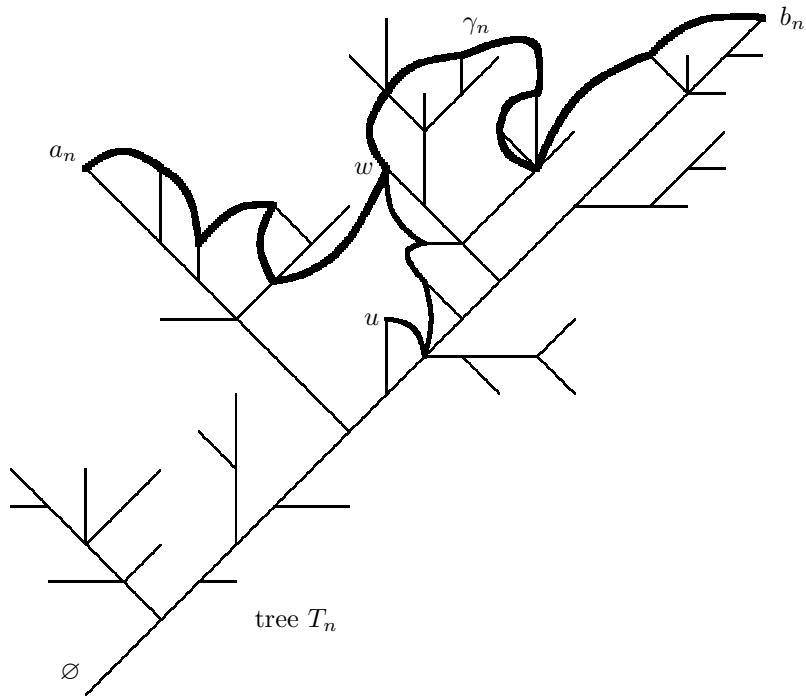

To complete this sketch, we explain why the claim holds. Again, we need to
go back to the discrete setting. We consider a vertex $u\in[a_n,b_n]$ such that
\begin{enumerate}
\item[(i)] $L_n(u) < L_n(a_n) + \eps n^{1/4}$;
\item[(ii)] $L_n(v) > L_n(a_n) + \frac{\eps}{2}\,n^{1/4}\ ,\qquad \forall v\in\llbracket \phi^n_{a_n,b_n}(u),u\rrbracket$
\end{enumerate}
where $\phi^n_{a_n,b_n}(u)$ is the last ancestor of $u$ in the tree $T_n$ that belongs to $\llbracket a_n,b_n\rrbracket$. 

Condition (ii) guarantees that the vertex $u$ lies ``between'' $\llbracket a_n, b_n\rrbracket$
and the geodesic $\gamma_n$: If this were not the case, the geodesic $\gamma_n$
would contain a point in $\llbracket \phi^n_{a_n,b_n}(u),u\rrbracket$, which is impossible
by (ii) (we already noticed that the label of a vertex of the geodesic $\gamma_n$
must be $L_n(a_n)+o(n^{1/4}$). 

Consider the geodesic path from $u$ to $v_*$ in $Q_n$ that is obtained from
the successor geodesic chain $e\to s(e)\to s^2(e)\to\cdots$ starting from any corner $e$ of
$u$ in $T_n$. Since arcs in the CVS bijection do not cross 
edges of the tree and since we know that
the vertex $u$ lies in the area between $\llbracket a_n, b_n\rrbracket$
and the geodesic $\gamma_n$, the geodesic we have just constructed 
cannot ``cross'' $\llbracket a_n,b_n\rrbracket $ and so it must intersect $\gamma_n$
at a vertex $w$. This vertex $w$ is such that
$$L_n(u)-L_n(w)= d_{Q_n}(u,w).$$
Since $w$ belongs to $\gamma_n$, we have $d_{Q_n}(w,a_n)=o(n^{1/4})$, and therefore
$$L_n(u)-L_n(a_n)= d_{Q_n}(u,a_n) + o(n^{1/4}).$$
By (i), we now get
$$d_{Q_n}(u,a_n)\leq \eps n^{1/4} + o(n^{1/4}).$$
We have thus obtained a discrete analog of the claim. To get the continuous version 
as stated above, we just need to do a careful passage to the limit $n\to\infty$.
\hfill$\square$

\subsection{Hausdorff dimension}\label{sec:hausdorff-dimension-1}

The limit in distribution (along a suitable subsequence) in Theorem
\ref{sec:conv-as-metr-3} can be written as $(\t_\ee/\!\approx,D)$, and
the space $\t_\ee/\!\approx$ is completely identified: Roughly
speaking two vertices $a$ and $b$ of the CRT $\t_\ee$ are identified
if and only if they have the same label $Z_a=Z_b$ and if one can go
from $a$ to $b$ following the ``contour'' of the tree $\t_\ee$ and
visiting only vertices with larger label.  In order to prove
Conjecture \ref{sec:conv-as-metr-4}, it would be necessary to
characterize the distance $D$. Much is known about $D$ (in particular
Theorem \ref{sec:limit-laws-radius-1} characterizes the distribution
of the profile of distances from the distinguished point $\rho_*$, and
one can show that this profile has the same distribution if one
replaces $\rho_*$ by a ``typical'' point of $M$). Still the
characterization of $D$ remains an open problem.

Nevertheless, one can show that the ``Brownian map''
$(\t_\ee/\!\approx,D)$, that is, any of the random metric spaces that
can arise as the limit in Theorem \ref{sec:conv-as-metr-3}, has
Hausdorff dimension $4$ and is homeomorphic to the $2$-sphere. This
was proved in \cite{legall06} and \cite{lgp}. The remainder of these
notes will be devoted to the proof of these two results.

\begin{theorem}\label{sec:hausdorff-dimension-2}
  Almost surely, the space $(M,D)$ has Hausdorff dimension $4$.
\end{theorem}

The lower bound is an easy consequence of Lemma
\ref{volume-ball}. Recall that $\mathrm{vol}$ is the image measure of
Lebesgue measure on $[0,1]$ under $p_\ee$. We let $\mathrm{Vol}$ be
the induced measure on $(M,D)$, that is, the image of Lebesgue
measure on $[0,1]$ under the projection $\bp:[0,1]\to M$. Then Lemma
\ref{volume-ball} implies that a.s., for every $\delta\in (0,1)$, and
every $x\in M$, it holds that 
$$\limsup_{r\downarrow
  0}\frac{\mathrm{Vol}(B_D(x,r))}{r^{4-\delta}}=0\, ,$$ where
$B_D(x,r)=\{y\in M:D(x,y)<r\}$ is the open ball centered at $x$ with
radius $r$. This last fact, combined with standard density theorems
for Hausdorff measures, implies that a.s.\ the Hausdorff dimension of
$(M,D)$ is greater than or equal to $4-\delta$, for every $\delta\in
(0,1)$.

For the upper bound, we rely on the following easy lemma. 

\begin{lemma}\label{sec:hausdorff-dimension-3}
  Almost surely, for every $\alpha\in (0,1/4)$, the label process $Z$
  is H{\"o}lder continuous with exponent $\alpha$. 
\end{lemma}

\proof
This is obtained by the classical Kolmogorov continuity criterion, and
moment estimates for $Z$. Let $s,t$ be such that $0\leq s<t\leq 1$,
and recall that conditionally given $\ee$, $Z_s-Z_t$ is a Gaussian
random variable with variance $d_\ee(s,t)$. Consequently,
for every $p>0$, there exists $C_p\in (0,\infty)$ such that 
$$E[|Z_s-Z_t|^p\, |\, \ee]=C_p d_\ee(s,t)^{p/2}\, ,$$
and since $\ee$ is a.s.\ H{\"o}lder continuous with exponent
$2\alpha$, we deduce the existence of a (random) $C'_p\in (0,\infty)$
such that 
$$E[|Z_s-Z_t|^p\, |\, \ee]\leq C'_p |s-t|^{p\alpha}\, .$$
The desired H\"older continuity property then follows from 
an application of the classical Kolmogorov lemma. \cq

\medskip

From this, we deduce that the projection $\bp:[0,1]\to M$ is a.s.\
H{\"o}lder continuous with index $\alpha\in (0,1/4)$ as well. Indeed,
using the fact that $D\leq D^0$, where $D^0$ is defined in
\eqref{eq:2}, we get
\begin{eqnarray*}
  D(\bp(s),\bp(t))&=&D(s,t)\\
  &\leq& Z_s+Z_t-2\inf_{s\wedge t\leq u\leq s\vee
    t}Z_u\\
  &\leq& 2\sup_{s\wedge t\leq u,v\leq s\vee t}|Z_u-Z_v|\\
  &\leq&
  C''_p|s-t|^\alpha\, ,
\end{eqnarray*}
for some $C''_p\in (0,\infty)$.  The fact that the Hausdorff dimension
of $(M,D)$ is bounded above by $1/\alpha$ is then a classical consequence
of this last property. This
completes the proof of the theorem.

\section{The homeomorphism theorem}\label{sec:home-theor}

\begin{theorem}\label{sec:home-theor-1}
Almost-surely, the Brownian map $(M,D)$ is homeomorphic to the
$2$-sphere $\mathbb{S}^2$. 
\end{theorem}

This result was first obtained by Le Gall and Paulin \cite{lgp}, by
arguing directly on the quotient space $M=\TT_\ee/\approx$. More
precisely, Le Gall and Paulin observe that the equivalence relations $\sim_\ee$
and $\approx$ may be viewed as equivalence relations on the sphere
$\mathbb{S}^2$. Upon showing that the associated classes are closed,
arcwise connected, and have connected complements, one can then apply
a theorem due to Moore \cite{moore25}, showing that under these
hypotheses, the quotient $\mathbb{S}^2/\approx$ is itself homeomorphic
to $\mathbb{S}^2$. Here, we will adopt a different approach,
introduced in Miermont \cite{miermont08}, which relies more on the discrete
approximations described in these notes. The idea is roughly as
follows: Even though the property of being homeomorphic to
$\mathbb{S}^2$ is not preserved under Gromov-Hausdorff convergence,
this preservation can be deduced under an additional property, called
{\em regular convergence}, introduced by Whyburn. This property says
heuristically that the spaces under consideration do not have small
bottlenecks, i.e.\ cycles of vanishing diameters that separate the
spaces into two macroscopic components.

In this section, when dealing with elements of the space $\K$ of
isometry classes of pointed compact metric spaces, we will often omit
to mention the distinguished point, as its role is less crucial than
it was in Sections \ref{sec:limit-laws-radius} and
\ref{sec:ident-brown-map}.

\subsection{Geodesic spaces and regular
  convergence}\label{sec:geod-spac-regul}

A metric space $(X,d)$ is said to be a {\em geodesic} metric space if
for every $x,y\in X$, there exists an isometry $f:[0,d(x,y)]\to X$
such that $f(0)=x$ and $f(d(x,y))=y$. Any such $f$ is called a
geodesic path between $x$ and $y$. For instance, real trees are
geodesic metric spaces by Definition \ref{sec:real-trees}. The set
$\K_{\mathrm{geo}}$ of isometry classes of (rooted) compact geodesic
metric spaces is closed in $(\K,d_{GH})$, as shown in \cite{BBI}.

\begin{definition}\label{sec:geod-spac-regul-2}
  Let $((X_n,d_n),n\geq 1)$ be a sequence of compact geodesic metric
  spaces, converging to $(X,d)$ in $(\K,d_{GH})$. We say that the
  convergence is {\em regular} if for every $\eps>0$, one can find
  $\delta>0$ and $N\in \N$ such that, for every $n>N$, every closed
  path $\gamma$ in $X_n$ with diameter at most $\delta$ is homotopic
  to $0$ in its $\eps$-neighborhood.
\end{definition}

For instance, let $Y_n$ be the complement in the unit sphere $\mathbb{S}^2\subset\R^3$
of the open $1/n$-neighborhood of the North pole, and endow $Y_n$ with the
intrinsic distance induced from the usual Euclidean metric on $\R^3$
(so that the distance between $x,y\in Y_n$ is the minimal length of a
path from $x$ to $y$ in $Y_n$).  Let $X_n$ be obtained by gluing two
(disjoint) copies of $Y_n$ along their boundaries, and endow it with
the natural intrinsic distance. Then $X_n$ converges in the
Gromov-Hausdorff sense to a bouquet of two spheres, i.e.\ two
(disjoint) copies of $\mathbb{S}^2$ whose North poles have been
identified. However, the convergence is not regular, because the path
$\gamma$ that consists in the boundary of (either copy of) $Y_n$ viewed as a subset of
$X_n$ has vanishing diameter as $n\to\infty$, but is not homotopic to $0$
in its $\eps$-neighborhood for any $\eps\in (0,1)$ and for any $n$. 
Indeed, such an $\eps$-neighborhood is a cylinder, around which
$\gamma$ makes one turn. 

\begin{theorem}\label{sec:geod-spac-regul-1}
  Let $((X_n,d_n),n\geq 1)$ be a sequence of $\K_{\mathrm{geo}}$ that
  converges regularly to a limit $(X,d)$ that is not reduced to a
  point. If $(X_n,d_n)$ is homeomorphic to $\mathbb{S}^2$ for every
  $n\geq1$, then so is $(X,d)$.
\end{theorem}

This theorem is an easy reformulation of a result of Whyburn in the
context of Gromov-Hausdorff convergence; see the paper by Begle
\cite{begle44}.  In the latter, it is assumed that every $X_n$ should
be a compact subset of a compact metric space $(Z,\delta)$,
independent of $n$, and that $X_n$ converges in the Hausdorff sense to
$X$. This transfers to our setting, because, if $(X_n,d_n)$ converges
to $(X,d)$ in the Gromov-Hausdorff sense, then one can find a compact
metric space $(Z,\delta)$ containing isometric copies $X'_n,n\geq 1$
and $X'$ of $X_n,n\geq 1$ and $X$, such that $X'_n$ converges in the
Hausdorff sense to $X'$, see for instance \cite[Lemma
A.1]{GrPfWi08}. In \cite{begle44}, it is also assumed in the
definition of regular convergence that for every $\eps>0$, there
exist $\delta>0$ and $N\in\N$ such that, for every $n\geq N$, any two points of $X_n$ that lie at distance
$\leq \delta$ are in a connected subset of $X_n$ of diameter $\leq
\eps$. This condition is tautologically satisfied for geodesic
metric spaces, which is the reason why we work in this context.

\subsection{Quadrangulations seen as geodesic
  spaces}\label{sec:quadr-seen-as}

Theorem \ref{sec:geod-spac-regul-1} gives a natural method to prove
Theorem \ref{sec:home-theor-1}, using the convergence of
quadrangulations to the Brownian map, as stated in Theorem
\ref{sec:conv-as-metr-3}. However, the finite space $(V(Q_n),d_{Q_n})$
is certainly not a geodesic space, nor homeomorphic to the
$2$-sphere. Hence, we have to modify a little these spaces so that
they satisfy the hypotheses of Theorem \ref{sec:geod-spac-regul-1}.  We will
achieve this by constructing a particular\footnote{The way we do this
  is by no means canonical. For instance, the emptied cubes $X_f$ used
  to fill the faces of $\bq$ below could be replaced by
  unit squares for the $l^1$ metric. However, our choice avoids the
  existence of too many geodesic paths between vertices of the map in
  the surface where it is embedded.} graphical representation of $\bq$.

Let $(X_f,d_f),f\in F(\bq)$ be disjoint copies of the emptied unit cube ``with
bottom removed''
$$\mathcal{C}=[0,1]^3\setminus \left((0,1)^2\times [0,1)\right)\, ,$$ endowed
  with the intrinsic metric $d_f$ inherited from the Euclidean metric
  (the distance between two points of $X_f$ is the minimal
  Euclidean length {\bf of a path in $X_f$}).  Obviously each
  $(X_f,d_f)$ is a geodesic metric space homeomorphic to a closed disk of
  $\R^2$. We will write elements of $X_f$ in the form $(s,t,r)_f$, where 
  $(s,t,r)\in\mathcal{C}$ and the subscript $f$ is used to differentiate points
of the different spaces $X_f$. The boundary $\partial X_f$ is then the collection 
of all points $(s,t,r)_f$ for $(s,t,r)\in ([0,1]^2\setminus
(0,1)^2)\times \{0\}$.
  
Let $f\in F(\bq)$ and let  $e_1,e_2,e_3,e_4$ be the four oriented
edges incident to $f$ enumerated in a way consistent with the counterclockwise order
on the boundary
 (here the labeling of these edges is chosen arbitrarily among the $4$
  possible labelings preserving the cyclic order). We then define
$$\begin{array}{lll}
c_{e_1}(t)=(t,0,0)_f&\, ,&\qquad 0\leq t\leq
1\\ c_{e_2}(t)=(1,t,0)_f&\, ,&\qquad 0\leq t\leq
1\\ c_{e_3}(t)=(1-t,1,0)_f&\, ,&\qquad 0\leq t\leq
1\\ c_{e_4}(t)=(0,1-t,0)_f&\, ,&\qquad 0\leq t\leq 1\, .
\end{array}$$
In this way, for every oriented edge $e$ of
the map $\bq$, we have defined a path $c_e$ which goes along one of
the four edges of the square $\partial X_f$, where $f$ is the face located to the left of
$e$.

We define an equivalence relation $\equiv$ on the disjoint union $\amalg_{f\in
  F(\bq)}X_f$, as the coarsest equivalence relation such that, for
every oriented edge $e$ of $\bq$, and every $t\in[0,1]$, we have
$c_e(t)\equiv c_{\ov{e}}(1-t)$. By identifying points of the same equivalence
class, we glue the oriented sides of the squares $\partial X_f$
pairwise, in a way that is consistent with the map structure. More
precisely, the topological quotient $\s_\bq:=\amalg_{f\in
  F(\bq)}X_f/\equiv$ is a surface which has a $2$-dimensional cell
complex structure, whose $1$-skeleton $\e_\bq:=\amalg_{f\in
  F(\bq)}\partial X_f/\equiv$ is a representative of the map $\bq$,
with faces ($2$-cells) $X_f\setminus\partial X_f$. In particular,
$\s_\bq$ is homeomorphic to $\mathbb{S}^2$ by \cite[Lemma
3.1.4]{MoTh01}. With an oriented edge $e$ of $\bq$ one associates an
edge of the graph drawing $\e_\bq$ in $\s_\bq$, more simply called an
edge of $\s_\bq$, made of the equivalence classes of points in
$c_e([0,1])$ (or $c_{\ov{e}}([0,1])$). We also let $\v_\bq$ be the
$0$-skeleton of this complex, i.e.\ the vertices of the graph ---
these are the equivalent classes of the corners of the squares
$\partial X_f$. We call them the vertices of $\s_\bq$ for simplicity.

We next endow the disjoint union $\amalg_{f\in F(\bq)} X_f$ with the
largest pseudo-metric $D_\bq$ that is compatible with $d_f,f\in
F(\bq)$ and with $\equiv$, in the sense that $D_\bq(x,y)\leq d_f(x,y)$
for $x,y\in X_f$, and $D_\bq(x,y)=0$ for $x\equiv y$. Therefore, the
function $D_\bq:\amalg_{f\in F(\bq)} X_f\times \amalg_{f\in F(\bq)}
X_f\to \R_+$ is compatible with the equivalence relation $\equiv$, and its
quotient mapping defines a pseudo-metric on the quotient space
$\s_\bq$, which is still denoted by $D_\bq$. 

\begin{proposition}\label{sec:turn-quadr-into-1}
The space $(\s_\bq,D_\bq)$ is a geodesic metric space homeomorphic to
$\mathbb{S}^2$. Moreover, the space
$(\v_\bq,D_\bq)$ is isometric to $(V(\bq),d_\bq)$, and any geodesic path in
$\s_\bq$ between two elements of $\v_\bq$ is a concatenation of edges
of $\s_\bq$. Last,
$$d_{GH}((V(\bq),d_\bq),(\s_\bq,D_\bq))\leq 3\, .$$
\end{proposition}

\proof We first  check that $D_\bq$ is a true metric on
$\s_\bq$, i.e.\ that it separates points. To see this, we use the fact
\cite[Theorem 3.1.27]{BBI} that $D_\bq$ admits the constructive
expression:
\begin{align*}
&D_\bq(a,b)\\&\;=\inf\left\{ \sum_{i=0}^nd(x_i,y_i):n\geq
0,x_0=a,y_n=b,y_i\equiv x_{i+1}\hbox{ for }0\leq i\leq n-1\right\},
\end{align*}
where we have set
$d(x,y)=d_f(x,y)$ if $x,y\in X_f$ for some $f$, and $d(x,y)=\infty$
otherwise. It follows that, for $a\in X_f\setminus \partial X_f$ and
$b\neq a$, $D_\bq(a,b)>\min(d(a,b),d_f(a,\partial X_f))$ \hfill $>0$, so $a$			%CMI
and $b$ are separated.

To verify that $D_\bq$ is a a true metric on $\s_\bq$, it remains to treat the case 
where $a\in \partial X_f,b\in\partial X_{f'}$
for some $f,f'\in F(\bq)$.  The crucial observation is that a shortest path in
$X_f$ between two points of $\partial X_f$ is entirely contained in
$\partial X_f$. It is then a simple exercise to check that if $a,b$ are 
in distinct equivalence classes, the distance $D_\bq(a,b)$ will be larger than the
length of some fixed non-trivial path with values in $\e_\bq$. More
precisely, if (the equivalence classes of) $a,b$ belong to the same edge of
$\s_\bq$, then we can find representatives $a',b'$ in the same $X_f$
and we will have $D_\bq(a,b)\geq d_f(a',b')$. If the equivalence class of $a$ is
not a vertex of $\s_\bq$ but that of $b$ is, then $D_\bq(a,b)$ is at
least equal to the distance of $a\in X_f$ to the closest corner of the
square $\partial X_f$. Finally, if the (distinct) equivalence classes
of $a,b$ are both vertices, then $D_\bq(a,b)\geq 1$.  One deduces that
$D_\bq$ is a true distance on $\s_\bq$, which makes it a geodesic metric
space by \cite[Corollary 3.1.24]{BBI}. Since $\s_\bq$ is a
compact topological space, the metric $D_\bq$ induces the quotient
topology on $\s_\bq$ by \cite[Exercise 3.1.14]{BBI}, hence
$(\s_\bq,D_\bq)$ is homeomorphic to $\mathbb{S}^2$.

From the observations in the last paragraph, a shortest path between
vertices of $\s_\bq$ takes values in $\e_\bq$. Since an edge
of $\s_\bq$ is easily checked to have length $1$ for the distance
$D_\bq$, such a shortest path will have the same length as a geodesic path
for the (combinatorial) graph distance between the two vertices. Hence
$(\v_\bq,D_\bq)$ is indeed isometric to $(V(\bq),d_\bq)$. The last
statement follows immediately from this and the fact that
$\diam(X_f,d_f)\leq 3$, entailing that $\v_\bq$ is $3$-dense in
$(\s_\bq,D_\bq)$, i.e.\ its $3$-neighborhood in $(\s_\bq,D_\bq)$
equals $\s_\bq$. \cq

\medskip

In view of the proposition, we can view $D_\bq$ as an extension to
$\mathcal{S}_\bq$ of the graph distance $d_\bq$ on $V(\bq)$. For this
reason, we will denote $D_\bq$ by $d_\bq$ from now on, which should
not set any ambiguity.

%Note that if $\gamma$ is a loop of diameter $\leq D$ in $S_\bq$, then
%we can find a simple loop $\gamma'$ of diameter $\leq D+4$ in $S_\bq$
%that takes its values only in edges of $S_\bq$, i.e.\ only in the
%points of $\amalg_{f\in F(\bq)}\partial X_f/\sim$. Such a simple path
%is obtained as the boundary of the faces of $S_\bq$ that contain some
%point of $\gamma$. Consequently, to prove $1$-regular convergence for
%quadrangulations, it suffices to prove that any simple loop of
%diameter $\leq \delta$ and made of edges of $S_\bq$, splits $S_\bq$
%into two connected components, one of which is of diameter $\leq
%\eps$.

\subsection{Proof of the homeomorphism
  theorem}\label{sec:proof-theorem-reft1-1}

We now work in the setting of the beginning of subsection
\ref{sec:brownian-map-as}. Recall that the uniform pointed
quadrangulation $(Q_n,v_*)$ is encoded by a uniform random element
$(T_n,L_n)$ of $\bT_n$ via the CVS bijection (the parameter
$\epsilon\in \{-1,1\}$ will play no role here), and that $C_n$ and
$V_n$ are the contour and label processes of $(T_n,L_n)$. We assume that
the amost sure convergence  \eqref{Groconv} holds uniformly on
$[0,1]^2$, along the sequence $(n_k)$, which is fixed. In what
follows, all  convergences as $n\to\infty$ hold
along this sequence, or along some further subsequence.

We can also assume that $(V(Q_n),d_{Q_n})$ is actually the (isometric)
space $(\v_{Q_n},$ \hfill $d_{Q_n})$, i.e.\ the subspace of vertices of the				%CMI
space $(\s_{Q_n},d_{Q_n})$ constructed in the previous
subsection. Recalling from subsection \ref{sec:la-bijection-de} that, in
the CVS bijection, each edge of the tree $T_n$ lies in exactly one
face of $Q_n$, we may and will assume that $T_n$ is also embedded in
the surface $\s_{Q_n}$, in such a way that the set of its vertices is
$\v_{Q_n}\setminus\{v_*\}$, where $v_*\in V(Q_n)$ is identified with
its counterpart in $\v_{Q_n}$, and that each edge of $T_n$ lies
entirely in the corresponding face of $\s_{Q_n}$ via the CVS
bijection.

We will rely on the following lemma. Let
$\mathrm{Sk}(\TT_\ee)$ be the complement of the set of leaves in the
CRT $\TT_\ee$. Equivalently, $\mathrm{Sk}(\TT_\ee)$ is the set of all points $a\in\TT_\ee$
such that $\TT_\ee\setminus\{a\}$ is disconnected, and it also coincides with
the set of all $a\in \TT_\ee$ that can be written $a=p_\ee(s)=p_\ee(s')$
for some $0\leq s< s'<1$. The set $\mathrm{Sk}(\TT_\ee)$ is called the {\em
  skeleton} of $\TT_\ee$.

\begin{lemma}\label{sec:estim-lengths-geod-2}
  The following property is true with probability $1$. Let $a\in
  \mathrm{Sk}(\TT_\ee)$, and let $s\in (0,1)$ be such that
  $a=p_\ee(s)$. Then for every $\eps>0$, there exists $t\in
  (s,(s+\eps)\wedge 1)$ such that $Z_t<Z_s$. 
\end{lemma}

This lemma is a consequence of \cite[Lemma 3.2]{lgp} (see also \cite[Lemma 2.2]{legall06}
for a slightly weaker statement). The proof relies on a precise
study of the label function $Z$, and we refer the interested reader to
\cite{lgp}.  Note that this result (and the analogous statement
derived by time-reversal) implies that a.s.,
if $a\in \mathrm{Sk}(\TT_\ee)$, then in each component of
$\TT_\ee\setminus\{a\}$, one can find points $b$ that are arbitrarily
close to $a$ and such that $Z_b<Z_a$.

\begin{lemma}\label{sec:proof-theorem-reft1}
  Almost surely, for every $\eps>0$, there exists $\delta\in(0,\eps)$
  such that, for $n$ large enough, any simple loop $\gamma_n$ made of
  edges of $\s_{Q_n}$, with diameter $\leq n^{1/4}\delta$, splits $\s_{Q_n}$
  in two Jordan domains, one of which has diameter $\leq n^{1/4}\eps$.
\end{lemma}

\proof We argue by contradiction. Assume that, with positive
probability, along some (random) subsequence of $(n_k)$ there exist
simple loops $\gamma_n$ made of edges of $\s_{Q_n}$, with diameters
$o(n^{1/4})$ as $n\to\infty$, such that the two Jordan domains bounded
by $\gamma_n$ are of diameters $\geq n^{1/4}\eps$, where $\eps>0$ is
some fixed constant. From now on we argue on this event. By abuse of
notation we will sometimes identify the chain $\gamma_n$ with the
set of vertices it visits, or with the union of its edges, in a way that should be
clear from the context.

By the Jordan curve theorem, the path $\gamma_n$ splits $\s_{Q_n}$
into two Jordan domains, which we denote by
$\mathcal{D}_n$ and $\mathcal{D}_n'$.  Since the diameters of both these
domains are at least $n^{1/4}\eps$, and since every point in
$\s_{Q_n}$ is at distance at most $3$ from some vertex, we can find
vertices $y_n$ and $y'_n$ belonging to $\mathcal{D}_n$ and $\mathcal{D}_n'$
respectively, and which lie at distance at
least $n^{1/4}\eps/4$ from $\gamma_n$. Since $V(Q_n)=T_n\cup\{v_*\}$,
we can always assume that $y_n$ and $y_n'$ are distinct from
$v_*$. Now, consider the geodesic path from $y_n$ to
$y'_n$ in $T_n$, and let $x_n$ be the first vertex of this path that belongs to
$\gamma_n$. 

In the contour exploration around $T_n$, the vertex $x_n$ is visited at
least once in the interval between $y_n$ and $y'_n$, and another time in the
interval between $y'_n$ and $y_n$. More precisely, let $j_n$ and
$j'_n$ be such that $y_n=u^n_{j_n},y'_n=u^n_{j'_n}$, and assume first
that $j_n<j'_n$ for infinitely many $n$. For such $n$, we can find
integers $i_n\in (j_n,j'_n)$ and $i'_n\in (0,j_n)\cup(j'_n,2n)$ such
that $x_n=u^n_{i_n}=u^n_{i'_n}$.  Up to further extraction, we may and
will assume that
\begin{equation}\label{eq:4}
  \frac{i_n}{2n}\to s\, ,\qquad \frac{i'_n}{2n}\to s'\, ,\qquad \frac{j_n}{2n}\to t\, ,\qquad
  \frac{j'_n}{2n}\to t'\, ,
\end{equation}
for some $s,s',t,t'\in [0,1]$ such that $t\leq s\leq t'$ and $s'\in
[0,t]\cup [t,1]$. Since
$$d_{Q_n}(x_n,y_n)\wedge d_{Q_n}(x_n,y'_n)\geq n^{1/4}\eps/4\, ,$$ we
deduce from \eqref{Groconv} that $D(s,t),D(s',t),D(s,t'),D(s',t')>0$,
and in particular, $s,s',t,t'$ are all distinct. Since
$u^n_{i_n}=u^n_{i'_n}$, we conclude that $s\sim_\ee s'$, so that
$p_\ee(s)\in \mathrm{Sk}(\TT_\ee)$. One obtains the same conclusion by
a similar argument if $j_n>j_n'$ for every $n$ large. We let
$x=p_\ee(s)$ and $y=p_\ee(t)$. Note that $y\not =x$ because $D(s,t)>0$ (recall
Lemma \ref{quotientCRT}).

Since $x\in \mathrm{Sk}(\TT_\ee)$, by Theorem
\ref{equival-relation}  we deduce 
that $D(a_*,x)=D(s_*,s)>0$, where
$a_*=p_\ee(s_*)$ is as before the a.s.\ unique leaf of $\TT_\ee$
where $Z$ attains its minimum. In particular, we obtain by \eqref{eq:15},
\eqref{Groconv} and the fact that $\diam(\gamma_n)=o(n^{1/4})$ that
$$\liminf_{n\to\infty}n^{-1/4}d_{Q_n}(v_*,\gamma_n)=
\liminf_{n\to\infty} n^{-1/4}d_{Q_n}(v_*,x_n)>0\, .$$ Therefore, for
$n$ large enough, $v_*$ does not belong to $\gamma_n$, and for
definiteness, we will assume that for such $n$, $\mathcal{D}_n$ is the
component of $\s_{Q_n}\setminus\gamma_n$ that does not contain
$v_*$.

Now, we let $L_n^+=L_n-\min L_n+1$, and in the rest of this proof, we
call $L_n^+(v)=d_{Q_n}(v_*,v)$ the {\em label of the vertex $v$} in
$Q_n$.  Let $l_n=d_{Q_n}(v_*,\gamma_n)=\min_{v\in \gamma_n} L_n^+(v)$
be the minimal distance from $v_*$ to a point visited by $\gamma_n$.
Note that, for every vertex $v\in \mathcal{D}_n$, the property
$L_n^+(v)\geq l_n$ holds, since any geodesic chain from $v_*$ to $v$ in $Q_n$ has to
cross $\gamma_n$. 

Recalling that the vertex $x_n$ was chosen so that the simple path in
$T_n$ from $x_n$ to $y_n$ lies entirely in $\mathcal{D}_n$, we
conclude that the labels of vertices on this path are all greater than
or equal to $l_n$. By passing to the limit, one concludes that for
every $c$ in the path $\llbracket x,y \rrbracket$ in $\TT_\ee$, there holds that $Z_c\geq
Z_x$. Since the process $Z$ evolves like Brownian
motion along line segments of the tree $\t_\ee$, we deduce that for every $c\in \llbracket x,y \rrbracket$ close enough to $x$, we have in
fact $Z_c>Z_x$. From the interpretation of line segments 
in $\TT_\ee$ in terms of the coding function $\ee$ (see the end of
subsection \ref{SScodingrealtree}), we can find $\ov{s}\in (0,1)$ such that
$p_\ee(\ov{s})=x$, and such that, for every $u>\ov{s}$ sufficiently close to $\ov{s}$, the
intersection of $\llbracket x,p_\ee(u)\rrbracket$ with $\llbracket x,y\rrbracket$
will be of the form $\llbracket x,p_\ee(r)\rrbracket$ for some $r\in(\ov{s},u]$. By
Lemma \ref{sec:estim-lengths-geod-2}, and the fact that
$Z_c\geq
Z_x$ for every $c\in \llbracket x,y \rrbracket$ close enough to $x$,
we can find $u>\ov{s}$ encoding a point
$a=p_\ee(u)$ and some $\eta>0$ such that $Z_a\leq
Z_x-(9/8)^{1/4}\eta$, and such that $\llbracket x,a \rrbracket \cap \llbracket x,y \rrbracket
=\llbracket x,b\rrbracket$ for some
$b\neq x$ such that $Z_b\geq Z_x+(9/8)^{1/4}\eta$.

\begin{figure}[h!]
\begin{center}
\includegraphics{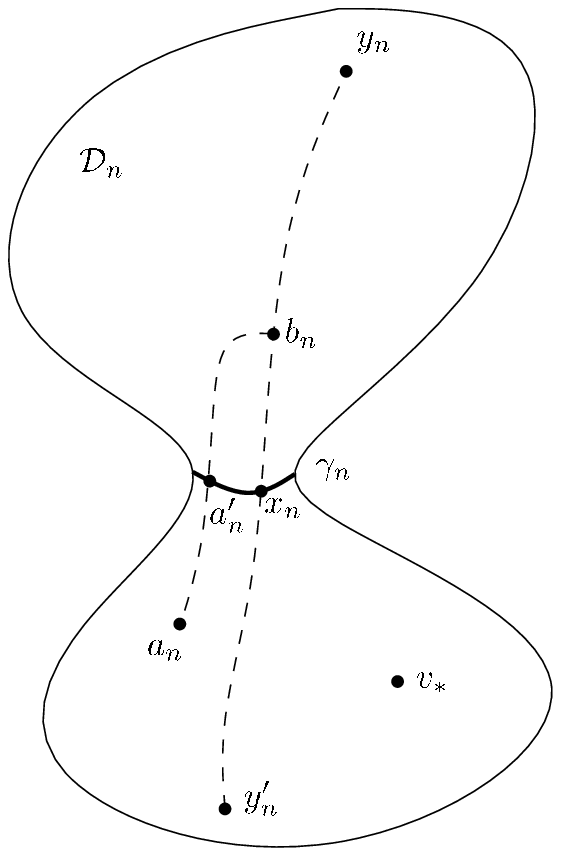}
\caption{Illustration of the proof. The surface $\s_{Q_n}$ is depicted
  as a sphere with a bottleneck circled by $\gamma_n$ (thick line).
  The dashed lines represent paths of $T_n$ that are useful in the
  proof: One enters the component $\mathcal{D}_n$, and the other goes
  out after entering, identifying in the limit a point of the skeleton
  with another.} \label{fig:1}
\end{center}
\end{figure}

We then go back once again to the discrete approximations of the
Brownian map, by considering $k_n$ such that $k_n/2n$ converges to
$u$. From the fact that $Z_a<Z_x$, we deduce that the vertex
$a_n=u^n_{k_n}$ has label $L^+_n(a_n)<l_n$ for every $n$ large
enough. Indeed, the convergence \eqref{Groconv} and the fact that
$\diam(\gamma_n)=o(n^{1/4})$ imply that $(9/8n)^{1/4}l_n\to Z_x-\inf
Z$. Consequently, the point $a_n$ does not belong to
$\mathcal{D}_n$. Moreover, the path in $T_n$ from $a_n$ to $x_n$ meets
the path from $x_n$ to $y_n$ at a point $b_n$ such that
$L_n^+(b_n)\geq l_n+\eta n^{1/4}$. The path from $a_n$ to $b_n$ has to
cross the loop $\gamma_n$ at some vertex, and we let $a'_n$ be the first such
vertex. By letting $n\to\infty$ one last time, we find a vertex
$a'\in \TT_\ee$, which in the appropriate sense is the limit of $a'_n$
as $n\to\infty$, such that $\llbracket a',x \rrbracket$ meets $\llbracket x,y \rrbracket$ at
$b$. In particular, $a'\neq x$. But since $a'_n$ and $x_n$ are both on
$\gamma_n$, we deduce that $D(a',x)=0$. This contradicts Theorem
\ref{equival-relation} because $x$ is not a leaf of $\TT_\ee$. This
contradiction completes the proof of the lemma. \cq

\medskip

We claim that Lemma \ref{sec:proof-theorem-reft1} suffices to verify
that the convergence of $(V(Q_n),$ \hfill $(9/8n)^{1/4}d_{Q_n})$ to $(M,D)$ is					%CMI
regular, and hence to conclude by Theorem \ref{sec:geod-spac-regul-1}
that the limit $(M,D)$ is a topological sphere. To see this, we first choose
$\eps<\diam(M)/3$ to avoid trivialities. Let $\gamma_n$ be a loop in
$\s_{Q_n}$ with diameter $\leq n^{1/4}\delta$. Consider the union of
the closures of faces of $\s_{Q_n}$ that are visited by $\gamma_n$.
The boundary of this union is a collection $\mathcal{L}$ of pairwise disjoint simple loops
made of edges of $\s_{Q_n}$. If $x,y$ belong to the preceding union of faces,
the fact that a face of $\s_{Q_n}$ has diameter less than
$3$ implies that there exist points $x'$ and $y'$ of $\gamma_n$ at distance at most $3$
from $x$ and $y$ respectively. Therefore, the diameters of the loops in $\mathcal{L}$ all are $\leq n^{1/4}\delta+6$.

By the Jordan Curve Theorem, each of these loops splits $\s_{Q_n}$ into
two simply connected components. By definition, one of these two
components contains $\gamma_n$ entirely. By Lemma
\ref{sec:proof-theorem-reft1}, one of the two components has diameter
$\leq n^{1/4}\eps$.  If we show that the last two properties hold
simultaneously for one of the two components associated with (at least) one of the loops
in $\mathcal{L}$, then obviously $\gamma_n$ will be homotopic to $0$ in
its $\eps$-neighborhood in $(\s_{Q_n},n^{-1/4}d_{Q_n})$. So assume the
contrary: The component not containing $\gamma_n$ associated with
every loop of $\mathcal{L}$ is of diameter $\leq n^{1/4}\eps$. If this
holds, then any point in $\s_{Q_n}$ must be at distance at most
$n^{1/4}\eps+3$ from some point in $\gamma_n$. Take $x,y$ such that
$d_{Q_n}(x,y)=\diam(\s_{Q_n})$. Then there exist points $x'$ and $y'$ in $\gamma_n$
at distance at most $n^{1/4}\eps+3$ respectively from $x$ and $y$, and we
conclude that $d_{Q_n}(x',y')\geq
\diam(\s_{Q_n})-6-2n^{1/4}\eps>n^{1/4}\delta\geq \diam(\gamma_n)$ for $n$
large enough by our choice of $\eps$. This contradiction completes the proof.

\medskip
{\it Note added in proof}. The uniqueness problem for the Brownian map
has been solved in two very recent papers of the authors: See
the preprints arxiv:1104.1606 and arxiv:1105.4842. Consequently, Conjecture 
\ref{sec:conv-as-metr-4}
is now a theorem, and analogs of this result hold for more general
random planar maps such as triangulations.

\end{document}